\documentclass[12pt]{degm}

\usepackage[latin1]{inputenc}
\usepackage{amsfonts}
\usepackage{amsmath}
\usepackage{amssymb}

\usepackage[all]{xy}

\newtheorem{theo}{Th{\'e}or{\`e}me}[section]
\newtheorem{prop}[theo]{Proposition}
\newtheorem{cor}[theo]{Corollaire}
\newtheorem{lemme}[theo]{Lemme}

\newtheorem{rque}[theo]{Remarque}
\newtheorem{defin}[theo]{D{\'e}finition}
\newtheorem{ex}[theo]{Exemple}

\renewcommand\refname{R{\'e}f{\'e}rences}

\DeclareMathOperator{\bs}{\backslash}

\DeclareMathOperator{\an}{an}
\DeclareMathOperator{\Aut}{Aut}

\DeclareMathOperator{\can}{can}
\DeclareMathOperator{\Conv}{Conv}
\DeclareMathOperator{\cusp}{cusp}
\DeclareMathOperator{\dR}{dR}
\DeclareMathOperator{\dlog}{dlog}
\DeclareMathOperator{\Det}{Det}
\DeclareMathOperator{\ev}{ev}
\DeclareMathOperator{\Ens}{Ens}
\DeclareMathOperator{\End}{End}
\DeclareMathOperator{\Fil}{Fil}

\DeclareMathOperator{\Gal}{Gal}
\DeclareMathOperator{\gal}{gal}
\DeclareMathOperator{\geom}{geom}
\DeclareMathOperator{\GL}{GL}
\DeclareMathOperator{\gr}{gr}
\DeclareMathOperator{\G}{G}
\DeclareMathOperator{\Hom}{Hom}

\DeclareMathOperator{\Isom}{Isom}
\DeclareMathOperator{\im}{im}
\DeclareMathOperator{\id}{id}
\DeclareMathOperator{\Lie}{Lie}
\DeclareMathOperator{\limproj}{\underset{\leftarrow}{lim}} 
\DeclareMathOperator{\limind}{\underset{\rightarrow}{lim}} 
\DeclareMathOperator{\Ker}{Ker}
\DeclareMathOperator{\N}{N}

\DeclareMathOperator{\ord}{ord}
\DeclareMathOperator{\Proj}{Proj}
\DeclareMathOperator{\pr}{pr}

\DeclareMathOperator{\Res}{Res}

\DeclareMathOperator{\rig}{rig}
\DeclareMathOperator{\Sch}{Sch}
\DeclareMathOperator{\SL}{SL}
\DeclareMathOperator{\Spec}{Spec}

\DeclareMathOperator{\Stab}{Stab}
\DeclareMathOperator{\Tate}{Tate}
\DeclareMathOperator{\Sym}{Sym}
\DeclareMathOperator{\Tr}{Tr}

\DeclareMathOperator{\A}{\mathbb{A}}
\DeclareMathOperator{\D}{\mathbb{D}}
\DeclareMathOperator{\p}{\mathbb{P}}

\DeclareMathOperator{\F}{\mathbb{F}}
\DeclareMathOperator{\Z}{\mathbb{Z}}
\DeclareMathOperator{\n}{\mathbb{N}}
\DeclareMathOperator{\Gm}{\mathbb{G}_m}
\DeclareMathOperator{\Q}{\mathbb{Q}}
\DeclareMathOperator{\R}{\mathbb{R}}
\DeclareMathOperator{\C}{\mathbb{C}}
\DeclareMathOperator{\cA}{\mathcal{A}}
\DeclareMathOperator{\bV}{\mathbb{V}}

\title{\Large Vari{\'e}t{\'e}s et formes modulaires de Hilbert  arithm{\'e}tiques pour 
$\Gamma_1(\mathfrak{c},\mathfrak{n})$}

\author{\large Mladen Dimitrov, Jacques Tilouine}

\begin{document}

\maketitle

\vspace{-1cm}
Cet article a pour origine deux expos{\'e}s donn{\'e}s par 
le second auteur {\`a} Varenna. Son contenu est cependant 
assez diff{\'e}rent de celui des expos{\'e}s.
Il traite de certains aspects arithm{\'e}tiques du lien entre
les vari{\'e}t{\'e}s et les formes modulaires de
Hilbert. Plusieurs points classiques ne sont pas abord{\'e}s faute de temps 
: op{\'e}rateurs de Hecke, th{\'e}orie de Serre-Tate. Par contre, cet article 
  donne des d{\'e}tails sur les compactifications toro{\"\i}dales
de la vari{\'e}t{\'e} ab{\'e}lienne 
de Hilbert-Blumenthal universelle (voir la partie \ref{kugasato}), ainsi
que plusieurs applications; la plupart, si elles sont peut-{\^e}tre
connues des experts, ne semblent pas figurer
dans les publications sur ce th{\`e}me; en particulier, celles 
relatives {\`a} la th{\'e}orie de Hodge,
aux formes de poids demi-entier et de Hilbert-Jacobi 
nous paraissent nouvelles. 

Nous avons grandement b{\'e}n{\'e}fici{\'e} d'un 
s{\'e}minaire sur ce sujet que nous avons organis{\'e} au premier semestre 
2001-02 {\`a} Paris 13.  Nous tenons {\`a} en
remercier tous les participants et en particulier G. Chenevier,
Y. Henrio, A. Mokrane, S. Morel et  S. Rozensztajn.
Nous voulons {\'e}galement remercier H. Hida, 
qui nous a {\'e}clair{\'e}s sur plusieurs points de
ce travail. Une partie de cet article a {\'e}t{\'e} r{\'e}dig{\'e}e 
alors que le second auteur
s{\'e}journait {\`a} l'Institut de Math{\'e}matiques de l'Universit{\'e} de M{\"u}nster
dans le cadre de la SFB 478 sur l'invitation de C. Deninger.
Il a appr{\'e}ci{\'e} les excellentes conditions de travail et
l'atmosph{\`e}re cordiale qui y r{\`e}gnent.

Nous avons divis{\'e} notre travail en deux articles; 
ainsi ce texte est muni d'un compagnon
 \cite{dimdg} qui donne les d{\'e}tails sur les compactifications 
(toro{\"\i}dales et minimale) des vari{\'e}t{\'e}s de Hilbert-Blumenthal
en niveau $\Gamma_1(\mathfrak{c},\mathfrak{n})$, en particulier 
aux pointes ramifi{\'e}es.
L'organisation du pr{\'e}sent article est la suivante :

\vspace{-1cm}
\tableofcontents

\section[Vari{\'e}t{\'e}s modulaires de Hilbert analytiques.]{Vari{\'e}t{\'e}s modulaires de Hilbert analytiques}

R{\'e}f{\'e}rences :  \cite{freitag}\cite{geer}

\vspace{-.4cm}
\subsubsection*{Notations.}

Soit $F$  un   corps  de nombres totalement r{\'e}el de degr{\'e} $d=d_F$, 
d'anneau des entiers $\mathfrak{o}$, de diff{\'e}rente $\mathfrak{d}$
et de discriminant $\Delta_F=\N_{F\!/\!\Q}(\mathfrak{d})$.
On note $J_F=\Hom_{\Q\!-\mathrm{alg.}}(F,{\C})$ l'ensemble de ses
plongements (r{\'e}els). On abr{\'e}gera $\N=\N_{F\!/\!\Q}$.

\smallskip

On se donne un groupe alg{\'e}brique $D_{/\Q}$, interm{\'e}diaire entre $\Gm_{/\Q}$ et 
$\Res^F_{\Q} \Gm$, connexe :

\vspace{-3mm}
$$ \Gm\hookrightarrow D\hookrightarrow \Res^F_{\Q} \Gm.$$

On d{\'e}finit le groupe alg{\'e}brique 
$G_{/\Q}$ (resp. $G^*_{/\Q}$) comme le produit fibr{\'e} de 
$D$ (resp. $\Gm$) et de  $\Res^F_{\Q} \GL_2$ au-dessus de  $\Res^F_{\Q} \Gm$.
On a le diagramme cart{\'e}sien suivant :

\vspace{-3mm}
$$\xymatrix@R=10pt{ \Res^F_{\Q}\SL_2 \ar[d]\ar@{^{(}->}[r] & G^* \ar[d]\ar@{^{(}->}[r] & 
G \ar[d]\ar@{^{(}->}[r] & \Res^F_{\Q} \GL_2\ar[d]^{\nu} \\
1 \ar@{^{(}->}[r]& \Gm \ar@{^{(}->}[r]& D \ar@{^{(}->}[r]& \Res^F_{\Q} \Gm,}$$

\vspace{-3mm}
\noindent o{\`u} la fl{\`e}che $\nu:\Res^F_{\Q} \GL_2\rightarrow \Res^F_{\Q} \Gm$ 
est  donn{\'e}e par la norme r{\'e}duite.

\smallskip

Le sous-groupe de Borel standard de $G$, 
son radical unipotent et  son tore maximal standard sont not{\'e}s 
$B$, $U$ et  $T$, respectivement. On pose $T_1=T\cap \ker(\nu)$. 

Pour toute $\Q$-alg{\`e}bre $R$ et pour tout groupe alg{\'e}brique $H$ sur $\Q$, on 
note $H_R$ le groupe de ses $R$-points. 

\smallskip

\begin{rque}
Dans toutes les applications le groupe $G$ sera  soit $\Res^F_{\Q} \GL_2$, soit 
$G^*$. Nous avons pr{\'e}f{\'e}r{\'e} ne pas fixer $G$ 
d{\`e}s le d{\'e}part, car $G^*$ intervient dans l'{\'e}tude g{\'e}om{\'e}trique 
des formes modulaires de Hilbert (le probl{\`e}me  de modules de vari{\'e}t{\'e}s 
ab{\'e}liennes de Hilbert associ{\'e} {\`a} $G^*$ est repr{\'e}sentable  : voir la partie
\ref{vahb}), 
alors que  $\Res^F_{\Q} \GL_2$ intervient dans l'{\'e}tude arithm{\'e}tique des
formes modulaires de Hilbert (les vari{\'e}t{\'e}s de Shimura 
associ{\'e}es {\`a}  $\Res^F_{\Q} \GL_2$ ne sont en g{\'e}n{\'e}ral que des
espaces de modules grossiers, mais  on conna{\^\i}t l'existence de 
repr{\'e}sentations galoisiennes associ{\'e}es aux  formes  modulaires de Hilbert propres
pour $\Res^F_{\Q} \GL_2$). Cette pr{\'e}sentation a {\'e}t{\'e} inspir{\'e}e par \cite{BL}.
\end{rque}

\vspace{-.4cm}
\subsubsection*{Le domaine sym{\'e}trique hermitien  $\mathfrak{H}_F$.}
Soit $(F\otimes \R)_+$ (resp. $G_{\R}^+$) la composante neutre de  $(F\otimes \R)^\times$  (resp. de $G_{\R}$).
Le groupe $G_{\R}^+$ agit par homographies 
sur l'espace 
$\mathfrak{H}_F=\{z\in F\otimes\C \enspace |\enspace \im(z)\in 
(F\otimes \R)_+\}\cong \mathfrak{H}^{J_F}$, o{\`u} $\mathfrak{H}$ d{\'e}signe le 
demi-plan de Poincar{\'e} (l'isomorphisme {\'e}tant donn{\'e} par 
$\xi\otimes w \mapsto (\tau(\xi)w)_{\tau\in J_F}$).

Posons $\underline{i}=(\sqrt{-1},...,\sqrt{-1})\in \mathfrak{H}_F$ et 
$K_{\infty}^+=\Stab_{G_{\R}^+}(\underline{i})$. Alors  $G_{\R}^+/K_{\infty}^+ \cong
\mathfrak{H}_F$, par $g\mapsto g(\underline{i})$ d'inverse
$x+\underline{i}y \mapsto 
\begin{pmatrix} y^{1/2} & xy^{-1/2} \\ 0 & y^{-1/2}\end{pmatrix}K_{\infty}^+$.

Via l'inclusion  $\mathfrak{H}_F\hookrightarrow \p^1(F\otimes\C)
\text{, }z\mapsto \Big[\begin{matrix} z \\ 1\end{matrix}\Big]$,
l'action de  $G_{\R}^+$ sur $\mathfrak{H}_F$ 
est compatible avec l'action naturelle de $G_{\C}$ sur $\p^1(F\otimes\C)$.

Les points rationnels 
$\p^1(F)$ du bord $\p^1(F\otimes\R)$ de $\mathfrak{H}_F$
sont appel{\'e}s les {\it pointes}. 
On pose $\infty= \Big[\begin{matrix} 1 \\ 0\end{matrix}\Big]$. 
Le groupe $G_{\Q}$ agit transitivement 
sur l'ensemble des pointes. On a  $B_{\Q}=\Stab_{G_{\Q}}(\infty)$ et
 $\p^1(F)\cong G_{\Q}/B_{\Q}$.

\medskip
 On munit l'espace $\mathfrak{H}_F^*=\mathfrak{H}_F\sqcup\p^1(F)$
de la topologie de Satake, donn{\'e}e par :

$\relbar$ la topologie usuelle sur  $\mathfrak{H}_F$,

$\relbar$ pour toute pointe 
 $\mathcal{C}\in \p^1(F)$, s'{\'e}crivant $\mathcal{C}=\gamma \infty$
avec $\gamma\in G_{\Q}$,
un syst{\`e}me fondamental de voisinages ouverts de $\mathcal{C}$ est donn{\'e} par
les $\{\gamma W_H\}_{H\in\R_+^*}$, o{\`u} 
$$W_H=\{z\in \mathfrak{H}_F
\enspace | \enspace \prod_\tau \im(z_\tau)>H\}.$$

%\noindent {\bf Fait : }
 
L'espace $\mathfrak{H}_F^*$
est s{\'e}par{\'e} (mais pas localement compact!) pour cette topologie 
(voir  \cite{freitag} I.2.9.).

\vspace{-.4cm}

\subsubsection*{Action de $G_{\Q}^+$ sur les $\mathfrak{o}$-r{\'e}seaux.}

Le groupe $G_{\Q}$ agit {\`a} gauche sur $F^2$, par
$\gamma\cdot(m,n)=(m,n)\gamma^{-1}$, o{\`u} 
$\gamma\in G_{\Q}$ et  $m,n \in F$.
Soit $G_{\Q}^+$ le sous-groupe de $G_{\Q}$ form{\'e} des {\'e}l{\'e}ments dont le 
d{\'e}terminant appartient au sous-groupe des {\'e}l{\'e}ments 
totalement positifs $F_+^\times$ de $F^\times$.
Posons $\mathfrak{o}_+^\times=\mathfrak{o}^\times\cap F_+^\times$.

\smallskip

Pour tout id{\'e}al fractionnaire $\mathfrak{f}$ de $F$
on pose $\mathfrak{f}^*=\mathfrak{f}^{-1}\mathfrak{d}^{-1}$.
On a un accouplement parfait $\Tr_{F\!/\!\Q} :
\mathfrak{f}\times\mathfrak{f}^* \rightarrow\Z$.
 
\smallskip

 Soit $L$ un $\mathfrak{o}$-r{\'e}seau de  $F^2$; c'est un 
$\mathfrak{o}$-module projectif de rang deux, donc il s'{\'e}crit, quitte {\`a} changer
la base de $F^2$, comme
$L=\mathfrak{e}\oplus\mathfrak{f}^*$,
avec $\mathfrak{e}$ et  $\mathfrak{f}$ des id{\'e}aux fractionnaires de $F$. 
Le stabilisateur 
du r{\'e}seau $\mathfrak{e}\oplus\mathfrak{f}^*$ dans $ G_{\Q}^+$ est {\'e}gal {\`a} :
$$ G^+(\mathfrak{e}\oplus\mathfrak{f}^*):=\{\gamma\in G_{\Q}^+| \det(\gamma)\in 
\mathfrak{o}_+^\times\} \cap
\begin{pmatrix} \mathfrak{o} & (\mathfrak{ef})^* \\ 
\mathfrak{efd} & \mathfrak{o}\end{pmatrix}$$

Lorsque $G=G^*$ (resp. $G=\Res^F_{\Q} \GL_2$), on {\'e}crit 
$\SL(\mathfrak{e}\oplus\mathfrak{f}^*)$ 
(resp. $\GL^+(\mathfrak{e}\oplus\mathfrak{f}^*)$), 
{\`a} la place de $ G^+(\mathfrak{e}\oplus\mathfrak{f}^*)$.
Notons que $\mathfrak{o}_+^\times\cap\Q=\{1\}$, et donc 
$\SL(\mathfrak{e}\oplus\mathfrak{f}^*)$ est form{\'e} d'{\'e}l{\'e}ments dont
le d{\'e}terminant vaut $1$.

\vspace{-.2cm}
\begin{lemme}
Dans la $\SL_2(F)$-orbite de  tout $\mathfrak{o}$-r{\'e}seau $L$ de $F^2$,
il existe un $\mathfrak{o}$-r{\'e}seau de la forme  $\mathfrak{o}\oplus \mathfrak{c}^*$,
avec $\mathfrak{c}$  un id{\'e}al fractionnaire de $F$.
\end{lemme}

\noindent {\bf D{\'e}monstration : } 
Il est clair que la  $\SL_2(F)$-orbite de  $L$  
 contient au moins un r{\'e}seau de la forme 
$\mathfrak{e}\oplus\mathfrak{f}^*$.
Prenons $a\in\mathfrak{e}$ et $c\in\mathfrak{fd}$ satisfaisant
$a\mathfrak{o}+c\mathfrak{ef}^*=\mathfrak{e}$. Par le th{\'e}or{\`e}me de
Bezout, il existe alors une matrice unimodulaire 
$\begin{pmatrix} a& b\\ c & d\end{pmatrix}\in \SL_2(F)\cap
\begin{pmatrix} \mathfrak{e} & \mathfrak{f}^* \\ 
\mathfrak{fd} & \mathfrak{e}^{-1}\end{pmatrix}$ et l'image de
$\mathfrak{e}\oplus\mathfrak{f}^*$ par cette matrice vaut 
 $\mathfrak{o}\oplus \mathfrak{c}^*$, avec $\mathfrak{c}^*=\mathfrak{ef}^*$. 
L'id{\'e}al $\mathfrak{c}^*$ est canoniquement isomorphe {\`a} $\wedge^2_{\mathfrak{o}}L$ et donc ne
d{\'e}pend pas  de la matrice de passage unimodulaire particuli{\`e}re choisie.
\hfill $\square$

En vertu de ce lemme, nous ne consid{\'e}rerons par la suite que des
$\mathfrak{o}$-r{\'e}seau de la forme $\mathfrak{o}\oplus \mathfrak{c}^*$,
avec $\mathfrak{c}$  un id{\'e}al fractionnaire de $F$.

\vspace{-.2cm}
\begin{rque} Consid{\'e}rons le cas o{\`u} $G=\Res^F_{\Q} \GL_2$. Alors
l'application $L\mapsto \wedge^2_{\mathfrak{o}}L$ induit une bijection 
entre l'ensemble des $G_{\Q}^+$-orbites de $\mathfrak{o}$-r{\'e}seaux $L$ de $F^2$
et le groupe de classes strictes d'id{\'e}aux $\mathrm{Cl}_F^+$. 

Notons que deux groupes  $\GL^+(\mathfrak{o}\oplus \mathfrak{c}^*)$
et  $\GL^+(\mathfrak{o}\oplus \mathfrak{c'}^*)$ sont conjugu{\'e}s dans 
$G_{\Q}^+$, si et seulement, si les id{\'e}aux  $\mathfrak{c}$ et $\mathfrak{c'}$
appartiennent au m{\^e}me genre 
(i.e. $\mathfrak{c'}=\xi \mathfrak{e}^2\mathfrak{c}$, 
avec $\xi\in F_+^\times$ et $\mathfrak{e}$ id{\'e}al de $F$). 

\end{rque}

\vspace{-.4cm}

\subsubsection*{Sous-groupes de congruence de $ G_{\Q}^+$.}
\label{sgc}
On fixe dans la suite un id{\'e}al fractionnaire $\mathfrak{c}$ et on consid{\`e}re 
le r{\'e}seau $L_0=\mathfrak{o}\oplus \mathfrak{c}^*$.

On se donne  aussi  un id{\'e}al $\mathfrak{n}\subsetneq \mathfrak{o}$.
Le  $\mathfrak{o}/\mathfrak{n}$-module  $\mathfrak{n}^{-1}L_0/L_0$
est  libre de rang $2$. Prenons $x_0\in F$, avec 
$\mathfrak{o}=\mathfrak{n}+x_0\mathfrak{cd}$. La multiplication par
$x_0$ induit alors les isomorphismes 
  $\mathfrak{o}/\mathfrak{n}\overset{\sim}{\longrightarrow}
\mathfrak{c}^*/\mathfrak{c}^*\mathfrak{n}$ et  
$\mathfrak{cd}/\mathfrak{cdn}\overset{\sim}{\longrightarrow}
\mathfrak{o}/\mathfrak{n}$
ce qui nous permet d'identifier  l'image de $\SL(\mathfrak{o}\oplus\mathfrak{c}^*)$
 dans $\Aut(\mathfrak{n}^{-1}L_0/L_0)$ avec  $\SL_2(\mathfrak{o}/\mathfrak{n})$ par
l'application $\begin{pmatrix} a & bx_0 \\ c & d\end{pmatrix}\mapsto
\begin{pmatrix} a & b \\ cx_0 & d\end{pmatrix}$, o{\`u} $a,b,d\in
\mathfrak{o}/\mathfrak{n}$ et $c \in \mathfrak{cd}/\mathfrak{cdn}$.
Faisons l'hypoth{\`e}se :

\medskip

\noindent {\bf (NT)}  $\mathfrak{n}$ est premier {\`a} $\N(\mathfrak{cd})$ 
et $\mathfrak{n}$ ne divise ni $2$, ni $3$.

\smallskip

\noindent Soit $ \Gamma_0^1(\mathfrak{c},\mathfrak{n})=
 \SL_2(F)\cap
\begin{pmatrix} \mathfrak{o} & \mathfrak{c}^* \\ 
\mathfrak{cdn} & \mathfrak{o}\end{pmatrix}$,
$ \Gamma^1(\mathfrak{c},\mathfrak{n})=\Ker(\SL(\mathfrak{o}\oplus \mathfrak{c}^*)\rightarrow \Aut(\mathfrak{n}^{-1}L_0/L_0))  $

\noindent $ \text{et } \Gamma_1^1(\mathfrak{c},\mathfrak{n})=\Big{\{}\gamma=
\begin{pmatrix}a &b \\c &d\end{pmatrix}
\in \Gamma_0^1(\mathfrak{c},\mathfrak{n})
\enspace \Big{|} \enspace  d\equiv 1 \pmod{\mathfrak{n}} \Big{\}}.$

La r{\'e}duction modulo $\mathfrak{n}$ induit un diagramme  cart{\'e}sien :
$$\xymatrix@R=10pt{ \Gamma^1(\mathfrak{c},\mathfrak{n})\ar[d]\ar@{^{(}->}[r] & 
\Gamma_1^1(\mathfrak{c},\mathfrak{n})\ar[d]\ar@{^{(}->}[r]& 
\Gamma_0^1(\mathfrak{c},\mathfrak{n})\ar[d]\ar@{^{(}->}[r]& 
\SL(\mathfrak{o}\oplus \mathfrak{c}^*) \ar[d]\\
{\begin{pmatrix}1 &0 \\0 &1\end{pmatrix}} \ar@{^{(}->}[r]& 
 {\begin{pmatrix}* &* \\0 &1\end{pmatrix}} \ar@{^{(}->}[r]& 
 {\begin{pmatrix}* &* \\0 &*\end{pmatrix}} \ar@{^{(}->}[r]& 
 \SL_2(\mathfrak{o}/\mathfrak{n})}$$

Les groupes $\Gamma^1(\mathfrak{c},\mathfrak{n})\subset
\Gamma_1^1(\mathfrak{c},\mathfrak{n})\subset
\Gamma_0^1(\mathfrak{c},\mathfrak{n})\subset
\SL(\mathfrak{o}\oplus \mathfrak{c}^*)$ sont des 
sous-groupes de congruence de $\SL_2(F)$.
De m{\^e}me on d{\'e}finit les sous-groupes de congruence : 
$$\Gamma(\mathfrak{c},\mathfrak{n})\subset
\Gamma_1(\mathfrak{c},\mathfrak{n})\subset
\Gamma_0(\mathfrak{c},\mathfrak{n})\subset
\GL^+(\mathfrak{o}\oplus \mathfrak{c}^*)\subset \GL_2^+(F),$$
$$\Gamma^D(\mathfrak{c},\mathfrak{n})\subset
\Gamma_1^D(\mathfrak{c},\mathfrak{n})\subset
\Gamma_0^D(\mathfrak{c},\mathfrak{n})\subset
G^+(\mathfrak{o}\oplus \mathfrak{c}^*)\subset G_{\Q}^+,$$
\noindent en rempla{\c c}ant la condition 
d'unimodularit{\'e} par celle d'avoir son d{\'e}terminant appartenant {\`a} 
$\mathfrak{o}_+^\times$ (resp. {\`a} $
\mathfrak{o}_{D+}^\times:=D_{\Q}\cap\mathfrak{o}_+^\times$).

\vspace{-.3cm} 
\begin{lemme}\label{NT}
 Sous l'hypoth{\`e}se  {\bf (NT)} le groupe 
 $\Gamma_1(\mathfrak{c},\mathfrak{n})$ est sans torsion.
\end{lemme}

\noindent{\bf D{\'e}monstration : }Par l'absurde. Supposons qu'il existe un {\'e}l{\'e}ment 
de $\Gamma_1(\mathfrak{c},\mathfrak{n})$ d'ordre premier $p$. Le 
d{\'e}terminant de cet {\'e}l{\'e}ment est une racine de l'unit{\'e} totalement positive, donc
{\'e}gale {\`a} $1$.  Cet {\'e}l{\'e}ment admet  comme valeur 
propre une racine de l'unit{\'e} $\zeta_p\neq 1$, ainsi que son inverse $\zeta_p^{-1}$.
En prenant sa trace on trouve que $\zeta_p$ est quadratique sur $F$, 
i.e.  $[F(\zeta_p):F]\leq 2$. Par ailleurs, $\zeta_p+\zeta_p^{-1}-2\in
\mathfrak{n}$, donc $\N(\mathfrak{n})$ est une puissance de $p$.
D'apr{\`e}s  {\bf (NT)}  on a alors que $F$ et 
$\Q(\zeta_p)$ sont lin{\'e}airement disjoints, d'o{\`u} 
$[F(\zeta_p):F]=p-1$. On en d{\'e}duit que $p=2$ ou $p=3$, 
ce qui implique, par un calcul facile, 
que $\mathfrak{n}$ divise $2$ ou  $3$. Contradiction.\hfill $\square$

\begin{rque} La condition {\bf (NT)} est optimale pour  que
$\Gamma_1(\mathfrak{c},\mathfrak{n})$ soit sans torsion. En effet, 
comme les matrices $\begin{pmatrix} -1 & 0  \\ 0 & -1\end{pmatrix}\in 
\Gamma_1^1(\mathfrak{d}^{-1},(2))$ 
et $\begin{pmatrix} -2 & 1  \\ -3 & 1\end{pmatrix}\in 
\Gamma_1^1(\mathfrak{d}^{-1},(3))$ sont d'ordre fini, 
$\mathfrak{n}$ ne peut diviser ni $2$   ni $3$.
Par ailleurs, la condition que $\mathfrak{n}$ soit  premier {\`a} 
$\Delta_F$ est aussi n{\'e}cessaire, comme le montre 
la matrice  d'ordre fini
$\begin{pmatrix} 1 & 1  \\ 
\frac{\sqrt{5}-5}{2} & \frac{\sqrt{5}-3}{2} \end{pmatrix}\in 
\Gamma_1^1(\mathfrak{d}^{-1},(\sqrt{5}))$ (ici $F=\Q(\sqrt{5})$).
Enfin, la condition que $\mathfrak{n}$ soit  premier {\`a} 
$\N(\mathfrak{c})$ est b{\'e}nigne, car par le th{\'e}or{\`e}me d'approximation 
faible, chaque classe de $\mathrm{Cl}_F^+$ contient  des id{\'e}aux  $\mathfrak{c}$ 
premiers {\`a} $\N(\mathfrak{n})$.
\end{rque}

{\bf Dans toute la suite du texte on suppose que l'hypoth{\`e}se 
(NT) est satisfaite, de sorte que 
$\Gamma_1(\mathfrak{c},\mathfrak{n})$ soit sans torsion.}

\vspace{-.4cm} 
\subsubsection*{Pointes pour les  sous-groupes de congruence.}

Soit $\Gamma$ un sous-groupe de congruence. Comme 
$F^\times \Gamma/F^\times$ est commensurable avec $\mathrm{PSL}_2(\mathfrak{o})$, 
l'ensemble de ses pointes est aussi  $\p^1(F)$ et l'ensemble  
$\Gamma\bs\p^1(F)$ est fini. 

Les deux lemmes suivants d{\'e}crivent les classes d'isomorphisme 
de $\Gamma_0^{\bullet}(\mathfrak{c},\mathfrak{n})$-pointes 
($\bullet=\varnothing,D,1$).
Cette description sera utilis{\'e}e 
dans la partie \ref{vahb}, o{\`u} nous {\'e}tudions les 
$\Gamma_1^\bullet(\mathfrak{c},\mathfrak{n})$-pointes.
 Notons que 
$\Gamma_0^D(\mathfrak{c},\mathfrak{n})=
G^+(\mathfrak{o}\oplus\mathfrak{c}^*)\cap 
G^+(\mathfrak{o}\oplus(\mathfrak{cn})^*)$.

\vspace{-.4cm} 
\begin{lemme}
Soient $\Big(\begin{matrix} a \\ c\end{matrix}\Big),
\Big(\begin{matrix} a' \\ c'\end{matrix}\Big)\in F^2-\{0\}$ et
 soit $\mathfrak{f}$ un id{\'e}al fractionnaire de $F$.\\
Si  $a\mathfrak{o}+c\mathfrak{f}^*=a'\mathfrak{o}+c'\mathfrak{f}^*$, alors
 il existe $\gamma\in 
 \SL(\mathfrak{o}\oplus\mathfrak{f}^*)$  tel que
$\Big(\begin{matrix} a' \\ c'\end{matrix}\Big)=\gamma \Big(\begin{matrix} a \\ c\end{matrix}\Big).$ 
\end{lemme}

\noindent{\bf D{\'e}monstration : } 
Posons $\mathfrak{b}=a\mathfrak{o}+c\mathfrak{f}^*$. Il existent 
$\gamma,\gamma'\in  G^1_{\Q}\cap 
\begin{pmatrix} \mathfrak{b} & (\mathfrak{bf})^* \\ 
\mathfrak{bfd} & \mathfrak{b}^{-1}\end{pmatrix}$ tels que 
 $\Big(\begin{matrix} a \\ c\end{matrix}\Big)=\gamma\infty$ et 
$\Big(\begin{matrix} a' \\ c'\end{matrix}\Big)=\gamma'\infty$.
Comme $\gamma'\gamma^{-1}\in\SL_2(\mathfrak{o}\oplus\mathfrak{f}^*)$
 on a le lemme. \hfill $\square$

En notant $\mathcal{I}_F$ l'ensemble des id{\'e}aux fractionnaires
et $\mathrm{Cl}_F$ le groupe des classes de $F$, on en d{\'e}duit :

\vspace{-.2cm} 
\begin{lemme}\label{classebiendef} On a deux bijections :

$il_{\mathfrak{c}}: G^+(\mathfrak{o}\oplus\mathfrak{c}^*)\bs 
\left(F^2-\{0\}\right)\overset{\sim}{\longrightarrow} \mathcal{I}_F\text{, }
\Big(\begin{matrix} a \\ c\end{matrix}\Big)\mapsto 
\mathfrak{b}=a\mathfrak{o}+c\mathfrak{c}^*$, et 

$il_{\mathfrak{cn}}: G^+(\mathfrak{o}\oplus(\mathfrak{cn})^*)\bs 
\left(F^2-\{0\}\right)\overset{\sim}{\longrightarrow} \mathcal{I}_F\text{, } 
\Big(\begin{matrix} a \\ c\end{matrix}\Big)\mapsto 
\mathfrak{b}'=a\mathfrak{o}+c(\mathfrak{cn})^*$,

qui induisent  deux bijections d'ensembles finis :

\medskip
$cl_{\mathfrak{c}}: G^+(\mathfrak{o}\oplus\mathfrak{c}^*)
\bs \mathbb{P}^1(F)\overset{\sim}{\longrightarrow} \mathrm{Cl}_F\text{, et }
cl_{\mathfrak{cn}}: G^+(\mathfrak{o}\oplus (\mathfrak{cn})^*)
\bs \mathbb{P}^1(F)\overset{\sim}{\longrightarrow} \mathrm{Cl}_F.$
\end{lemme}

\noindent{\bf D{\'e}monstration : } Les fl{\`e}ches sont bien d{\'e}finie (voir 
 l'action de $G_{\Q}^+$ sur $F^2$ d{\'e}finie plus haut). Le lemme pr{\'e}c{\'e}dent donne 
l'injectivit{\'e}. La surjectivit{\'e} d{\'e}coule du fait que tout id{\'e}al 
dans un corps de nombres peut {\^e}tre  engendr{\'e} par  deux de ses  {\'e}l{\'e}ments.

\vspace{-.4cm} 
\subsubsection*{Vari{\'e}t{\'e}s modulaires de Hilbert analytiques.}
{\'E}tant donn{\'e} un sous-groupe de congruence  $\Gamma$ 
on d{\'e}finit {\it la vari{\'e}t{\'e} modulaire de Hilbert analytique }
$M^{\an }=\Gamma\bs\mathfrak{H}_F$. La vari{\'e}t{\'e} $M^{\an }$ est lisse, 
si et seulement, si $\Gamma$ est sans torsion. En revanche 
$M^{\an }$ n'est jamais compacte.

Les vari{\'e}t{\'e}s modulaires de Hilbert, dont  nous {\'e}tudierons en d{\'e}tail la 
g{\'e}om{\'e}trie, sont celles
correspondant aux groupes de congruence 
$\Gamma_1^D(\mathfrak{c},\mathfrak{n})$.

\vspace{-.4cm} 
\subsubsection*{Compactification de Satake.}
L'espace quotient $M^{\an *}=\Gamma\bs\mathfrak{H}_F^*$
est compact pour la topologie de Satake. Il est l'union de 
$M^{\an }$ et d'un nombre fini de points, appel{\'e}s les pointes 
(voir  \cite{freitag} Sect.I).
Il est muni d'une structure de vari{\'e}t{\'e}
analytique complexe normale pour laquelle les pointes sont des 
points singuliers si $d_F>1$ (voir  \cite{freitag} II.4).

\vspace{-.4cm}

\section[Vari{\'e}t{\'e}s ab{\'e}liennes de Hilbert-Blumenthal et
formes de Hilbert.]{Vari{\'e}t{\'e}s ab{\'e}liennes de Hilbert-Blumenthal et
formes de Hilbert}\label{vahb}

  Dans la suite  $\Gamma$ (resp.  $\Gamma^1$)  d{\'e}signe le 
sous-groupe de congruence $\Gamma_1^D(\mathfrak{c},\mathfrak{n})$ 
 (resp.  $\Gamma_1^1(\mathfrak{c},\mathfrak{n})$) et on pose 
 $M^{\an}=
\Gamma_1^D(\mathfrak{c},\mathfrak{n})\!\bs\!\mathfrak{H}_F$
(resp.  $M^{1,\an}=
\Gamma_1^1(\mathfrak{c},\mathfrak{n})\!\bs\!\mathfrak{H}_F$).

\begin{defin} Une vari{\'e}t{\'e} ab{\'e}lienne {\`a} multiplication r{\'e}elle par $\mathfrak{o}$
sur un sch{\'e}ma $S$ est la donn{\'e}e d'un sch{\'e}ma   
ab{\'e}lien  $f:A\rightarrow S$ de dimension relative $d_F$ et 
d'une injection  $\iota:\mathfrak{o}\hookrightarrow \End(A/S)$.
\end{defin}

Soit $\mathfrak{c}$ un id{\'e}al fractionnaire.
Pour chaque  vari{\'e}t{\'e} ab{\'e}lienne {\`a} multiplication r{\'e}elle $A/S$, on d{\'e}finit
 un faisceau en $\mathfrak{o}$-modules
sur le gros site  {\'e}tale de $S$ en associant {\`a} un $S$-sch{\'e}ma $Y$
le $\mathfrak{o}$-module $A(Y)\otimes_{\mathfrak{o}}\mathfrak{c}$.
Ce foncteur 
est repr{\'e}sentable par une vari{\'e}t{\'e} ab{\'e}lienne {\`a} 
multiplication r{\'e}elle sur $S$, not{\'e}e 
$A\otimes_{\mathfrak{o}}\mathfrak{c}$ (voir \cite{DePa}); elle est caract{\'e}ris{\'e}e par:
$$A\otimes_{\mathfrak{o}}\mathfrak{c}=
{\begin{cases}A/A[\mathfrak{c}^{-1}]\text{, si }
 \mathfrak{c}^{-1} \text{ entier}. \\ 
(A^t\otimes\mathfrak{c}^{-1})^t\text{, si }
 \mathfrak{c}\text{ entier.} 
\end{cases}}$$

La premi{\`e}re formule s'obtient en tensorisant par $A$ sur $\mathfrak{o}$ 
la suite exacte  courte 
$0\rightarrow \mathfrak{o}\rightarrow \mathfrak{c}\rightarrow \mathfrak{c}/\mathfrak{o}\rightarrow 0$.
La seconde en r{\'e}sulte par dualit{\'e}.

{\`A} partir de $\iota:\mathfrak{o}\hookrightarrow \End(A/S)$
on obtient $\mathfrak{c}\hookrightarrow
\Hom_{\mathfrak{o}}(A,A\otimes_{\mathfrak{o}}\mathfrak{c})$.
Soit $\mathfrak{c}_+=\mathfrak{c}\cap (F\otimes \R)_+$.
Soit $\Sym_{\mathfrak{o}}(A,A^t)$ le $\mathfrak{o}$-module des homomorphismes sym{\'e}triques de $A$ vers $A^t$ et 
$\mathcal{P}(A)\subset \Sym_{\mathfrak{o}}(A,A^t)$
le c{\^o}ne des polarisations.

\vspace{-.2cm}
\begin{defin}\cite{DePa}
Une vari{\'e}t{\'e} ab{\'e}lienne $A$ de Hilbert-Blumenthal (abr{\'e}g{\'e}e en VAHB) sur un sch{\'e}ma $S$
est une vari{\'e}t{\'e} ab{\'e}lienne {\`a} multiplication r{\'e}elle par $\mathfrak{o}$,
v{\'e}rifiant la condition de Deligne-Pappas suivante:

 {\bf (DP)} il existe un isomorphisme $\mathfrak{o}$-{\'e}quivariant 
$\lambda: A\otimes\mathfrak{c} \overset{\sim}{\longrightarrow}A^t$ 
tel que via $\lambda$ on a 
$(\mathfrak{c},\mathfrak{c}_+)\cong (\Sym_{\mathfrak{o}}(A,A^t),\mathcal{P}(A)).$

Un tel isomorphisme $\lambda$ est appel{\'e} une $\mathfrak{c}$-polarisation.
\end{defin}

Le groupe $\mathfrak{o}_+^{\times}$ agit sur l'ensemble des
$\mathfrak{c}$-polarisations d'une VAHB $A/S$.

\vspace{-.2cm}
\begin{defin}\label{classe}
On appelle une classe de $\mathfrak{c}$-polarisation, une 
orbite $\overline{\lambda}$ de  $\mathfrak{c}$-polarisations
sous $\mathfrak{o}_{D+}^{\times}=\mathfrak{o}_+^{\times}\cap D_{\Q}$. 
\end{defin}

\vspace{-.2cm}
\begin{rque}
Si $\Delta_F$ est inversible dans $S$,
alors  la condition  {\bf (DP)} est {\'e}quivalente {\`a} la condition suivante 
de Rapoport \cite{rapoport} (voir \cite{DePa} Cor.2.9 et \cite{goren} Chap.3.5)

 {\bf (R)} le faisceau $\underline{\omega}=f_*\Omega^1_{A/S}$ est   
localement libre de rang $1$ sur $\mathfrak{o}\otimes\mathcal{O}_S$ pour la topologie de Zariski.
\end{rque}

\vspace{-.2cm}
\begin{defin}
Une $\mu_{\mathfrak{n}}$-structure
de niveau sur une VAHB $A/S$ est la donn{\'e}e d'une immersion ferm{\'e}e 
$\mathfrak{o}$-lin{\'e}aire de $S$-sch{\'e}mas en groupes finis
$\alpha : (\mathfrak{o}/\mathfrak{n})(1)\hookrightarrow A[\mathfrak{n}]$,
o{\`u} $(\mathfrak{o}/\mathfrak{n})(1)=(\Gm\otimes \mathfrak{d}^{-1})
[\mathfrak{n}]$  
d{\'e}signe le  dual de Cartier du $S$-sch{\'e}ma constant $\mathfrak{o}/\mathfrak{n}$.
\end{defin}

\begin{rque}
Comme $\mathfrak{c}$ est premier {\`a} $\mathfrak{n}$, 
la $\mathfrak{c}$-polarisation $\lambda$, combin{\'e}e avec l'accouplement 
de Weil $A[\mathfrak{n}]\times A^t[\mathfrak{n}]\rightarrow 
(\Gm\otimes \mathfrak{d}^{-1})[\mathfrak{n}]$, donne un accouplement 
$\mathfrak{o}$-{\'e}quivariant parfait 
$A[\mathfrak{n}]\times A[\mathfrak{n}]\rightarrow
(\Gm\otimes \mathfrak{c}^*)[\mathfrak{n}]$. {\'E}tant donn{\'e} une $\mu_{\mathfrak{n}}$-structure
de niveau $\alpha : (\mathfrak{o}/\mathfrak{n})(1)\hookrightarrow A[\mathfrak{n}]$,
{\`a} l'aide de ce dernier, on lui associe de mani{\`e}re canonique un morphisme 
$\mathfrak{o}$-lin{\'e}aire surjectif de $S$-sch{\'e}mas en groupes finis 
$\alpha^*: A[\mathfrak{n}]\rightarrow  \mathfrak{c}^{-1}/\mathfrak{nc}^{-1}$,
appel{\'e} le $\lambda$-dual de Cartier  de $\alpha$. On a une suite exacte  :
$$0 \rightarrow (\mathfrak{o}/\mathfrak{n})(1) 
\overset{\alpha}{\longrightarrow} A[\mathfrak{n}]
\overset{\alpha^*}{\longrightarrow} 
\mathfrak{c}^{-1}/\mathfrak{nc}^{-1}\rightarrow  0$$
\end{rque}

\vspace{-.4cm} 
\subsubsection*{Construction analytique de la VAHB universelle sur $M^{1,\an}$.}

Pour tout $z\in \mathfrak{H}_F$ et $\gamma\in G_{\R}$
on pose $j(\gamma,z)=c\cdot z+d\in(F\otimes\C)^\times$. 
D'apr{\`e}s l'identit{\'e} $j(\gamma\gamma',z)=j(\gamma,\gamma'(z))j(\gamma',z)$ 
on a  un $1$-cocycle :
$$ G_{\R}\longrightarrow 
(\mathfrak{o}\otimes\mathcal{O}_{\mathfrak{H}_F})^\times
\text{, } \gamma\mapsto(z\mapsto j(\gamma,z)).$$

\medskip
On pose $\mathcal{A}^{\an }=\Gamma^1\bs
(\mathfrak{H}_F\times(F\otimes\C))/\mathfrak{o}\oplus\mathfrak{c}^*$, 
o{\`u} le groupe produit semi-direct 
$(\mathfrak{o}\oplus\mathfrak{c}^*)\rtimes\Gamma^1$ 
(pour $\gamma\cdot(m,n)=(m,n)\gamma^{-1}$)
agit {\`a} gauche sur $\mathfrak{H}_F\times(F\otimes\C)$ par :
\begin{equation}\label{vahb-univ}
\begin{cases} \gamma(z,v)=(\gamma(z),j(\gamma,z)^{-1}v) \\
(z,v)(m,n)=(z,v+m\cdot z+n)\end{cases}.
\end{equation}

La fibre du point $\Gamma^1 z\in M^{1,\an }$ est  
la vari{\'e}t{\'e} ab{\'e}lienne $\mathcal{A}^{\an }_z:=(F\otimes\C)/L_z$, 
o{\`u}  $L_z=(\mathfrak{o} z\oplus\mathfrak{c}^*)$.
La fl{\`e}che $\iota(\xi):(z,v)\mapsto
(z,\xi v)$ induit une action de $\xi\in \mathfrak{o}$ sur $\mathcal{A}^{\an }$,
d'o{\`u} une injection 
$\iota : \mathfrak{o}\hookrightarrow\End(\mathcal{A}^{\an }/M^{1,\an })$.

Pour tout fibr{\'e} vectoriel $E$ sur $M^{1,\an }$, soit $E^\vee$ le fibr{\'e} dual.
Il est facile de voir que   $\Lie(\mathcal{A}^{\an }/M^{1,\an })=\Gamma^1\bs
(\mathfrak{H}_F\times(F\otimes\C))$ et 
$\underline{\omega}=\Lie(\mathcal{A}^{\an }/M^{1,\an })^\vee$ sont localement libres de rang $1$
sur $\mathfrak{o}\otimes\mathcal{O}_{M^{1,\an }}$.

\smallskip
Pour tout $\mathfrak{o}$-module $L$ on a 
un isomorphisme entre  $\Hom_{\mathfrak{o}}(L,\mathfrak{d}^{-1})$ 
et $L^*=\Hom_{\Z}(L,\Z)$, obtenu en composant avec $\Tr_{F\!/\!\Q}$.

\smallskip
On a un isomorphisme  $\mathfrak{o}$-lin{\'e}aire
$\wedge^2_{\mathfrak{o}}L_z\cong \mathfrak{c}^*$,
venant de l'accouplement parfait $\lambda_z :  L_z \times L_z\rightarrow F$ 
\enspace $(u,v)\mapsto \frac{uv^c-u^cv}{2i\im(z)}$. 
L'application $\Tr_{F\!/\!\Q}\circ \lambda_z$ nous fournit 
un isomorphisme  $ L_z\otimes_{\mathfrak{o}}\mathfrak{c}\cong L_z^*$,
d'o{\`u} une  $\mathfrak{c}$-polarisation 
$\mathcal{A}_z\otimes\mathfrak{c} \cong \mathcal{A}^{t}_z $.

Si  $\mathfrak{o}=\mathfrak{n}+y_0\mathfrak{c}$,
 la fl{\`e}che $M^{1,\an }\times \mathfrak{n}^{-1}\mathfrak{d}^{-1}/\mathfrak{d}^{-1}
\rightarrow \mathcal{A}^{\an }[\mathfrak{n}] $, 
$(z,v) \mapsto (z,y_0v)$ munit $\mathcal{A}^{\an }$ d'une 
$\mu_{\mathfrak{n}}$-structure de niveau.

\begin{prop}
 $(\mathcal{A}^{\an },\iota,\lambda,\alpha)/ M^{1,\an }$ est 
une VAHB $\mathfrak{c}$-polaris{\'e}e analytique, munie d'une 
$\mu_{\mathfrak{n}}$-structure de niveau.
 
La fl{\`e}che $\mathcal{A}^{\an }\rightarrow M^{1,\an }$
est universelle, i.e. pour toute VAHB analytique
$A/S$ munie d'une $\mu_{\mathfrak{n}}$-structure de niveau et d'une 
 $\mathfrak{c}$-polarisation, il existe une unique 
fl{\`e}che $\varphi: S\rightarrow M^{1,\an }$ et un unique isomorphisme de VAHB 
munies de $\mu_{\mathfrak{n}}$-structure de niveau et de 
$\mathfrak{c}$-polarisation
$A\cong \mathcal{A}^{\an }\times_{M^{1,\an }} S$. En particulier, si 
$A$ est une VAHB complexe 
 munie d'une $\mu_{\mathfrak{n}}$-structure de niveau et d'une 
 $\mathfrak{c}$-polarisation, il existe un unique
point $z\in M^{1,\an }$ et un unique isomorphisme 
$A\cong \mathcal{A}^{\an }_z$.
\end{prop}

\noindent{\bf Id{\'e}e de la d{\'e}monstration : } Il est clair que toute VAHB complexe est 
isomorphe {\`a} une VAHB de la forme  $\mathcal{A}^{\an }_z$ et que les 
deux  VAHB  analytiques 
 $\mathcal{A}^{\an }_z$ et $\mathcal{A}^{\an }_{z'}$ sont isomorphes
comme  VAHB munies de leurs $\mu_{\mathfrak{n}}$-structures de niveau et 
$\mathfrak{c}$-polarisations   si et seulement si $z'\in \Gamma^1 z$. 

Soit $A/S$ comme dans l'{\'e}nonc{\'e}. Par ce qui pr{\'e}c{\`e}de, il existe une unique 
fl{\`e}che ensembliste $\varphi: S\rightarrow M^{1,\an }$ telle que 
$A\cong \mathcal{A}^{\an }\times_{M^{1,\an }} S$.
L'analyticit{\'e} de $\varphi$ se v{\'e}rifie localement, car $\varphi(s)=
\int_{\gamma_1}\omega(s)/\int_{\gamma_2}\omega(s)$,  o{\`u} 
$(\gamma_1,\gamma_2)$ est une $\mathfrak{o}$-base locale convenable 
de l'homologie de $A/S$ et $\omega$ est une $\mathfrak{o}\otimes 
\mathcal{O}_S$-base 
locale de   $\underline{\omega}$.\hfill $\square$

\begin{rque} 
1) Notons qu'en g{\'e}n{\'e}ral pour $G\neq G^*$ la vari{\'e}t{\'e} 
$M^{\an }=\Gamma\bs \mathfrak{H}_F$ n'est qu'un espace de modules 
grossier pour le probl{\`e}me de modules de classes d'isomorphismes
de VAHB munies d'une classe de $\mathfrak{c}$-polarisation (voir la 
d{\'e}finition \ref{classe}) et d'une $\mu_{\mathfrak{n}}$-structure de niveau.

 Comme $\Gamma^1$ est  un sous-groupe distingu{\'e} de $\Gamma$,
le quotient $\mathfrak{o}_{D+}^\times$ agit sur   $M^{1,\an}$. Sur 
les $S$-points  $\epsilon\in \mathfrak{o}_{D+}^\times$
envoie $(A,\iota,\lambda,\alpha)/S$ sur $(A,\iota,\epsilon\lambda,\alpha)/S$. 
On a $M^{\an}=\mathfrak{o}_{D+}^\times\bs
M^{1,\an}$.

En fait, le sous-groupe $\mathfrak{o}_{D+}^{\times}\cap 
\mathfrak{o}_{\mathfrak{n}}^{\times 2}$ agit trivialement, car 
la multiplication par $\epsilon\in \mathfrak{o}^\times$
induit un isomorphisme $(A,\iota,\lambda,\alpha)\cong
(A,\iota,\epsilon^2\lambda,\epsilon\alpha)$. Donc 
$M^{1,\an}$ est un rev{\^e}tement fini {\'e}tale de
$M^{\an}$, de groupe $\mathfrak{o}_{D+}^\times/ \mathfrak{o}_{D+}^{\times}\cap 
\mathfrak{o}_{\mathfrak{n}}^{\times 2}$.

Pour toute VAHB $A/S$ munie d'une classe de $\mathfrak{c}$-polarisation 
et d'une $\mu_{\mathfrak{n}}$-structure de niveau on a des fl{\`e}ches
$S\rightarrow M^{1,\an }$ dont les compos{\'e}es avec la 
projection $M^{1,\an }\rightarrow M^{\an }$ co{\"\i}ncident et telles
que $A/S$ avec sa  classe de $\mathfrak{c}$-polarisation  soit
le pull back de $A^{\an }/M^{1,\an }$ munie 
de la   classe de sa $\mathfrak{c}$-polarisation universelle.

2) Lorsque $G=\Res^F_{\Q}\GL_2$, Hida,
dans son livre \cite{hida-padic} Chap.4 Sect. 4.1.2, a donn{\'e} une autre description de
$M^{\an}$ comme espace de modules grossier des VAHB avec classes 
de $F^{\times}_{+}$-polarisation.  
Dans sa description,
$$M^{\an}=M^{\an}_1(\mathfrak{c},\mathfrak{n})=F_+^\times\bs
\underset{\mathfrak{c'}}{\coprod}M_1^1(\mathfrak{c'},\mathfrak{n})^{\an},$$
 o{\`u} $\mathfrak{c'}$ d{\'e}crit les id{\'e}aux de $F$ qui 
appartiennent {\`a} la m{\^e}me classe stricte que $\mathfrak{c}$. 
\end{rque}

\vspace{-.4cm} 
\subsubsection*{VAHB analytique de Tate.}

Soit $\mathcal{C}=\gamma\infty$ une $\Gamma$-pointe ($\gamma\in G^*_{\Q}$).
On commence par {\'e}tudier la forme d'un voisinage de $\mathcal{C}$ dans $M^{\an }$,
puis on va d{\'e}crire celle de $\mathcal{A}^{\an }$, au-dessus d'un tel 
voisinage dans $M^{1,\an }$.

\underline{1 cas} $\mathcal{C}=
\Big[\begin{matrix}1\\0\end{matrix}\Big]=\infty$. $\Stab_\Gamma(\infty)=\bigg{\{}
\begin{pmatrix} u\epsilon & b \\ 0 & u^{-1} \end{pmatrix}\enspace \Big|
\enspace u\in\mathfrak{o}_{\mathfrak{n}}^\times, 
\epsilon\in\mathfrak{o}_{D+}^\times
\enspace b\in\mathfrak{c}^*\bigg{\}}=
\mathfrak{c}^*\rtimes (\mathfrak{o}_{\mathfrak{n}}^\times
\times \mathfrak{o}_{D+}^\times)$, o{\`u}
$\mathfrak{o}_{D+}^\times=\mathfrak{o}_+^\times\cap D_{\Q}$ et 
pour tout id{\'e}al $\mathfrak{f}$ de $F$,
 $\mathfrak{o}_{\mathfrak{f}}^\times$ d{\'e}signe 
le groupe des unit{\'e}s de $ \mathfrak{o}$ 
congrues {\`a} $1$ modulo $\mathfrak{f}$.

Soit $\phi$ l'inclusion naturelle $F\rightarrow F\otimes\C$; 
on a une suite exacte courte: 
\begin{equation}\label{qexp}
0\rightarrow \mathfrak{c}^*
\overset{\phi}{\longrightarrow}
 F\otimes\C \overset{q}{\longrightarrow}
\Gm\otimes\mathfrak{c}^*\rightarrow 1,
\end{equation}
obtenue par produit tensoriel par $\mathfrak{c}^* $ de
$0\rightarrow \Z \rightarrow \C\overset{e^{2i\pi\cdot}}
{\longrightarrow} \C^*\rightarrow 1$.

Pour $m\in F$ et $z\in F\otimes\C$, 
on pose  $q_z^m=q(\phi(m)z)$ ($=q(\phi(m)z+\phi(n))$ pour tout 
 $n\in \mathfrak{c}^*$). On voit facilement:

\medskip
\noindent {\bf Fait : } 
Pour $H>0$ assez grand 
$\Stab_\Gamma(W_{H})=B_\Gamma:=\Gamma\cap B_{\Q}$. 

\medskip
Le groupe 
$ \mathfrak{o}_{\infty}^\times:=
\mathfrak{o}_{\mathfrak{n}}^\times\times \mathfrak{o}_{D+}^\times$ agit sur 
le quotient $D_{H}=\mathfrak{c}^*\bs W_{H}$ par 

$(u,\epsilon)\cdot(z+\mathfrak{c}^*)=u^2\epsilon z+\mathfrak{c}^*$, et on  a
le diagramme suivant :
$$ \xymatrix@R=15pt@C=17pt{ \mathfrak{H}_F\ar[d] & W_{H} 
\ar@{_{(}->}[l]\ar[d] & & & & \\
M^{\an }  &  
B_\Gamma\bs W_{H} \ar@{_{(}->}[l] & &
D_{H}=\mathfrak{c}^*\bs W_{H} 
\ar@{^{(}->}[r]\ar_{{ }\hspace{-4mm}\mod{\mathfrak{o}_{\infty}^\times}}[ll] &
F\otimes\C/\phi(\mathfrak{c}^*) 
\ar^{\hspace{-1mm}\underset{\sim}{q}}[r] &
\Gm\otimes\mathfrak{c}^*=:S_\infty}$$

Le diagramme suivant d{\'e}crit la structure de la VAHB universelle 
$\mathcal{A}^{\an }$
sur  le voisinage $B_{\Gamma^1} \bs W_{H}$ de la pointe $\infty$ dans 
$M^{1,\an}$ :
$$\xymatrix@R=15pt@C=15pt{
 B_{\Gamma^1}\bs(W_{H}\times F\otimes\C)/\mathfrak{o}\oplus\mathfrak{c}^*\ar[d] &
\mathfrak{c}^*\bs
(W_{H}\times F\otimes\C )/\mathfrak{o}\oplus\mathfrak{c}^* \ar[l]\ar[d]\ar@{^{(}->}[r] &
(\Gm\otimes\mathfrak{c}^*\times \Gm\otimes\mathfrak{c}^*)/q_z^{\mathfrak{o}} \ar[d]\\
B_{\Gamma^1}\bs W_{H}& D_{H} 
\ar@{^{(}->}[r]\ar_{\mod{\mathfrak{o}_{\mathfrak{n}}^\times}}[l] & 
\Gm\otimes\mathfrak{c}^*=:S_\infty }$$

\noindent{\bf Commentaires :} 1) La notation 
$q_z^{\mathfrak{o}}$ exprime que $m\in {\mathfrak{o}}$ agit
sur $\Gm\otimes\mathfrak{c}^*\times \Gm\otimes\mathfrak{c}^*$ par 
la formule $(q_z,q_v)\cdot m=(q_z,q_v q_z^m)$.

2) Le groupe $\mathfrak{o}_{\mathfrak{n}}^\times$ agit sur 
$S_\infty\times \Gm\otimes\mathfrak{c}^*$ par 
$u \cdot (q_z,q_v)=(q_z^{u^2},q_v^{u})$.

\begin{defin}
 La VAHB $\mathfrak{c}$-polaris{\'e}e au-dessus de $S_\infty$ 
ainsi obtenue s'appelle la VAHB analytique 
de Tate, not{\'e}e $\Tate_{\mathfrak{c},\mathfrak{o}}(q_z)$. 
Sa fibre au point $q_z\in S_\infty$ est {\'e}gale {\`a} 
$\Gm\otimes\mathfrak{c}^*/q_z^{\mathfrak{o}}$.
\end{defin}

\underline{2 cas} $\mathcal{C}=\Big[\begin{matrix}a\\c\end{matrix}\Big]=
\gamma\infty$, $\gamma=\begin{pmatrix}a & b \\c & d \end{pmatrix}\in G^*_{\Q}$.
$\Stab_\Gamma(\mathcal{C})=B_{\Gamma,\mathcal{C}}:=\Gamma\cap \gamma B_{\Q}\gamma^{-1}$. Un syst{\`e}me fondamental de voisinages  de la  pointe $\mathcal{C}$ est donn{\'e} 
par les $B_{\Gamma,\mathcal{C}} \bs \gamma W_{H}$.
Notons que pour tout sous-groupe $\Gamma'$ de $ G_{\Q}$
on a la suite exacte suivante :
$$1\rightarrow \Gamma'\cap U_{\Q}
\rightarrow \Gamma'\cap B_{\Q}
\rightarrow \pr(\Gamma'\cap B_{\Q})
\rightarrow 1,$$ 
o{\`u} $\pr: B_{\Q}\rightarrow T_{\Q}$ est la projection canonique.
Le diagramme suivant :
$$ \xymatrix@R=15pt{ \mathfrak{H}_F \ar[r]& 
M^{\an }=\Gamma\bs\mathfrak{H}_F\ar^{\sim}[r] & 
\gamma^{-1}\Gamma\gamma\bs\mathfrak{H}_F  & & \\
\gamma W_{H} \ar[r]\ar@{^{(}->}[u] & 
B_{\Gamma,\mathcal{C}}\bs \gamma W_{H} \ar^{{}\hspace{-5mm}\sim}[r]\ar@{^{(}->}[u]&
\gamma^{-1}\Gamma\gamma\cap B_{\Q}\bs  W_{H} \ar@{^{(}->}[u]
& & \gamma^{-1}\Gamma\gamma\cap U_{\Q}\bs  W_{H}, 
\ar_{{}\hspace{-1mm} \pr(\gamma^{-1}\Gamma\gamma\cap B_{\Q})}[ll]}$$

\noindent  permet de nous ramener au cas de la pointe 
$\infty$, pour  le groupe $\gamma^{-1}\Gamma\gamma$.

\bigskip

\noindent$\bullet$ Calcul de $\gamma^{-1}\Gamma\gamma\cap U_{\Q}$.
$\begin{pmatrix}1 & \xi^* \\0 & 1 \end{pmatrix}\in \gamma^{-1}\Gamma\gamma\iff
\begin{pmatrix}1+ac\xi^* & a^2\xi^* \\-c^2\xi^* & 1-ac\xi^* \end{pmatrix}\in \Gamma$

\medskip
\noindent$\iff \xi^*\in  
a^{-2}\mathfrak{c}^{*}\cap (ac)^{-1}\mathfrak{n}\cap c^{-2}\mathfrak{c}^{*-1}\mathfrak{n}
=(a^2 \mathfrak{c}^{*-1}+ac\mathfrak{n}^{-1}+c^2(\mathfrak{cn})^{*})^{-1}=$

\medskip
\noindent$=\mathfrak{c}^{*}(a^2\mathfrak{o}+ac(\mathfrak{cn})^{*}+
c^2\mathfrak{c}^{*2}\mathfrak{n}^{-1})^{-1}=\mathfrak{c}^{*}
(a\mathfrak{o}+c\mathfrak{c}^{*})^{-1}
(a\mathfrak{o}+c(\mathfrak{cn})^{*})^{-1}.$

\smallskip

Donc $\gamma^{-1}\Gamma\gamma\cap U_{\Q}=\mathfrak{(cbb')}^*$,
avec  $\mathfrak{b}=a\mathfrak{o}+c\mathfrak{c}^{*}\quad et\quad 
\mathfrak{b'}=a\mathfrak{o}+c(\mathfrak{cn})^{*}$.

\medskip

Posons  $X:=\mathfrak{cbb'}$ (sa classe 
est bien d{\'e}finie, d'apr{\`e}s le lemme \ref{classebiendef}).

\bigskip

\noindent$\bullet$ 
Calcul de $\mathfrak{o}_{\mathcal{C}}^\times:=
\pr(\gamma^{-1}\Gamma\gamma\cap B_{\Q})$. 
Posons $\mathfrak{o}_{\mathcal{C},1}^\times:=
\mathfrak{o}_{\mathcal{C}}^\times\cap T_{1,\Q}$. Alors l'on a une suite
exacte courte $1 \rightarrow \mathfrak{o}_{\mathcal{C},1}^\times
\rightarrow \mathfrak{o}_{\mathcal{C}}^\times
\overset{\nu}{\rightarrow} \nu(\mathfrak{o}_{\mathcal{C}}^\times)
\rightarrow 1$.
En prenant  $\gamma={\begin{pmatrix} a &0 \\ c & a^{-1}\end{pmatrix}}$ on a 
$$\mathfrak{o}_{\mathcal{C}}^\times=
\{ (u,\epsilon)\in F^\times\times \mathfrak{o}_{D+}^\times  |
\exists \xi^*\in a^{-2}\mathfrak{c}^{*},
 u^{-1}-1+ac\xi^*\in \mathfrak{n}, 
u\epsilon-u^{-1}-ac\xi^*\in
 \frac{a}{c}\mathfrak{c}^{*-1}\mathfrak{n}\}.$$

 Le groupe $\mathfrak{o}_{\mathcal{C}}^\times$
 ne d{\'e}pend que de la classe de $\gamma$
dans $\Gamma\bs  G_{\Q}/B_{\Q}$.
Un calcul d{\'e}montre que l'on a  
 $\mathfrak{o}^\times
\supset \mathfrak{o}_{\mathcal{C},1}^\times\supset 
\mathfrak{o}_{\mathfrak{n}}^\times$. 
Si l'id{\'e}al  $\mathfrak{n}$ est sans  facteurs carr{\'e}s, 
alors $\mathfrak{o}_{\mathcal{C},1}^\times=
\mathfrak{o}_{\mathfrak{n}}^\times$. 
Le calcul explicite du groupe 
$\mathfrak{o}_{\mathcal{C}}^\times$ dans le cas g{\'e}n{\'e}ral, est  un corollaire
d'une autre description des $\Gamma$-pointes, donn{\'e}e dans \cite{dimdg} Prop 3.3.

\begin{rque} En g{\'e}n{\'e}ral l'inclusion $\gamma^{-1}\Gamma\gamma\cap T_{1,\Q}
\subset \mathfrak{o}_{\mathcal{C},1}^\times$ est stricte, bien que ce soit
une {\'e}galit{\'e} pour la pointe $\infty$. N{\'e}anmoins le groupe 
$\gamma^{-1}\Gamma\gamma\cap T_{1,\Q}$
est d'indice fini dans $\mathfrak{o}^\times$.
\end{rque}

\noindent Le type de la pointe $\mathcal{C}$ est d{\'e}termin{\'e} par :
$\begin{cases}
\text{l'ideal } X^*, \\
\text{le groupe }  \mathfrak{o}_{\mathcal{C}}^\times, \\
\text{l'action de } \mathfrak{o}_{\mathcal{C}}^\times 
 \text{ sur }  X^* \bs W_{H}.
\end{cases}$

 Le fait de remplacer $\gamma$  par $\gamma
\begin{pmatrix} a' & c' \\ 0 & a'{}^{-1}\end{pmatrix}$,
multiplie  $X^*$ par $a'{}^{-2}$ et 
conjugue l'action de $\mathfrak{o}_{\mathcal{C}}^\times$
sur $X^* \bs W_{H}$, par l'isomorphisme
 $W_{H}\rightarrow W_{\N(a')^2H},$ 
 $z\mapsto a'{}^2 z+a'c'$.

\bigskip

Pour {\'e}tudier la VAHB universelle $\mathcal{A}^{\an}/M^{1,\an}$ 
au  voisinage de la pointe $\mathcal{C}$, 
trouvons quel r{\'e}seau est stable par 
$\gamma^{-1}\SL(\mathfrak{o}\oplus \mathfrak{c}^*)\gamma$.
 Par le th{\'e}or{\`e}me de Bezout on peut prendre 
$\gamma={\begin{pmatrix} a &b \\ c & d\end{pmatrix}}
\in \SL_2(F) \cap
\begin{pmatrix} \mathfrak{b}  & (\mathfrak{bc})^* \\ 
\mathfrak{bcd}& \mathfrak{b}^{-1}
\end{pmatrix}$, o{\`u} $\mathfrak{b}=a\mathfrak{o}+c\mathfrak{c}^*$.
Posons  $\mathfrak{a}=\mathfrak{bc}$.
Comme $\gamma^{-1}$ transforme le r{\'e}seau $L_0=\mathfrak{o}\oplus \mathfrak{c}^*$
en le r{\'e}seau $L=\mathfrak{b}\oplus\mathfrak{a}^*$
$$\gamma^{-1}\SL(\mathfrak{o}\oplus \mathfrak{c}^*)\gamma=
\SL(\mathfrak{b}\oplus\mathfrak{a}^*) =
\SL_2(F)\cap
\begin{pmatrix} \mathfrak{o}& \mathfrak{b}^{-2}\mathfrak{c}^* \\
\mathfrak{b}^2 \mathfrak{c}^{*-1} &\mathfrak{o}\end{pmatrix}.\text{ Donc}$$
$$\xymatrix@R=15pt@C=15pt{ 
\mathcal{A}^{\an }\ar[d] & B_{\Gamma^1,\mathcal{C}}\bs 
(\gamma W_{H}\times F\otimes\C )/\mathfrak{o}\oplus\mathfrak{c}^* 
\ar@{_{(}->}[l]\ar[d]\ar^{\hspace{-5mm}\sim}[r] &
\gamma^{-1}\Gamma^1\gamma\cap B_{\Q}\bs 
(W_{H}\times F\otimes\C)/\mathfrak{b}\oplus \mathfrak{a}^* \ar[d]\\
M^{1,\an } & B_{\Gamma^1,\mathcal{C}}\bs \gamma W_{H}\ar@{_{(}->}[l]
\ar^{\hspace{-5mm}\sim}[r]&
\gamma^{-1}\Gamma^1\gamma\cap B_{\Q}\bs W_{H}   }$$

\noindent car $\gamma:\mathcal{A}^{\an }\rightarrow \mathcal{A}^{\an }$\enspace
$(z,v)\mapsto (\gamma z, j(\gamma,z)^{-1}v)$ et $\gamma^{-1}$
envoie $\mathfrak{o}\oplus\mathfrak{c}^*$ sur 
$\mathfrak{b}\oplus \mathfrak{a}^* $.

\medskip
 A partir de l{\`a}, en posant 
${A}^{\an }_{\gamma,H}=\gamma^{-1}\Gamma^1\gamma\cap B_{\Q}\bs 
(W_{H}\times(F\otimes\C))/\mathfrak{b}\oplus \mathfrak{a}^*,$
 on a la description de la vari{\'e}t{\'e} universelle au
 voisinage de la pointe $\mathcal{C}$ :
$$\xymatrix@R=15pt@C=15pt{{A}^{\an }_{\gamma,H} 
 \ar[d] & X^* \bs 
(W_{H}\times(F\otimes\C))/\mathfrak{b}\oplus \mathfrak{a}^* 
 \ar[l]\ar[d]\ar@{^{(}->}[r] &
(\Gm\otimes X^*\times \Gm\otimes \mathfrak{a}^*)/
q_z^{\mathfrak{b}} \ar[d]\\
\gamma^{-1}\Gamma^1\gamma\cap B_{\Q}\bs W_{H} &
 X^* \bs  W_{H} 
\ar_{\hspace{4mm}\mathrm{mod}\enspace \mathfrak{o}_{\mathcal{C},1}^\times}[l]
\ar@{^{(}->}[r] &
\Gm\otimes X^*=:S_{\mathcal{C}},}$$

 Le groupe  $\mathfrak{b}$  agit
sur le tore $\Gm\otimes X^*\times \Gm\otimes\mathfrak{a}^*$ par 
 $(q_z,q_v)\cdot \beta=(q_z,q_v q_z^{\beta})$.
Le groupe  $\mathfrak{o}_{\mathcal{C},1}^\times$ agit sur 
$S_\mathcal{C}\times \Gm\otimes \mathfrak{a}^*$ 
par  $u \cdot (q_z,q_v)
=(q_z^{u^2}q_{\xi_{u}^*}^{u},q_v^{u}) $, o{\`u}
$\xi_{u}^*$ est un {\'e}l{\'e}ment de  $(\mathfrak{b}^2\mathfrak{c})^*$, 
bien d{\'e}fini modulo $X^*$, et tel que 
$\begin{pmatrix} u & \xi_{u}^* \\ 0 & u^{-1}\end{pmatrix}\in
\gamma^{-1}\Gamma^1\gamma$.

On rappelle que, par d{\'e}finition, pour tout 
$m\in F$, $z\in F\otimes\C$ on pose  
$$q_z^m=q(\phi(m)z)=q(\phi(m)z+\phi(n)) \text{ pour tout } n\in X^*\text{, o{\`u}}$$
$$0\rightarrow X^*\overset{\phi}{\longrightarrow}
 F\otimes\C \overset{q}{\longrightarrow}
\Gm\otimes\mathfrak{c}^*\rightarrow 1.$$

\begin{defin} \label{vahb-tate}
 La VAHB $\mathfrak{c}$-polaris{\'e}e au-dessus de $S_\mathcal{C}$ ainsi obtenue, 
s'appelle la VAHB analytique 
de Tate, not{\'e}e $\Tate_{\mathfrak{a},\mathfrak{b}}(q_z)$. 
Sa fibre au point $q_z\in S_\mathcal{C}$ est {\'e}gale {\`a} 
$\Gm\otimes \mathfrak{a}^*/q_z^{\mathfrak{b}}$.
\end{defin}

\vspace{-.4cm} 

\subsubsection*{Formes modulaires de Hilbert de niveau 
$\Gamma=\Gamma_1^D(\mathfrak{c},\mathfrak{n})$.}

Rappelons que $\Z[J_F]$ s'identifie au groupe 
des caract{\`e}res du tore $\Res^F_{\Q} \Gm$
par $\kappa=\sum_{\tau\in J_F}k_\tau\tau\mapsto 
(x\mapsto \prod \tau(x)^{k_\tau})$. On note ce caract{\`e}re $x\mapsto x^\kappa$
et on utilisera la notation additive pour la loi de groupe sur les caract{\`e}res.
Les {\'e}l{\'e}ments de $\Z[J_F]$ sont appel{\'e}s des {\it poids}.

On suppose d{\'e}sormais $F\neq \Q$.
Pour tout poids $\kappa=\sum_{\tau\in J_F}k_\tau\tau$, on peut d{\'e}finir
l'espace des formes automorphes de Hilbert holomorphes de poids $\kappa$ 
et niveau $\Gamma$ comme
l'espace des fonctions holomorphes $f:\mathfrak{H}_F\rightarrow \C$ 
telles que pour tout $\gamma\in \Gamma$
$$ f(\gamma(z))=\nu(\gamma)^{-\kappa/2}
j(\gamma,z)^{\kappa} f(z).$$

Ce sont les sections du fibr{\'e} inversible analytique $\underline{\omega}^\kappa$ 
sur $M^{\an}$  donn{\'e} par le cocycle 
$$\Gamma\rightarrow \mathcal{O}^\times_{\mathfrak{H}_F},\quad 
\gamma\mapsto \nu(\gamma)^{-\kappa/2}j(\gamma,z)^\kappa.$$

Cependant, on ne s'int{\'e}resse dans la suite de ce texte qu'aux formes qui peuvent
intervenir dans la cohomologie de la vari{\'e}t{\'e} de Hilbert {\`a} coefficients dans un syst{\`e}me local 
alg{\'e}brique (c'est-{\`a}-dire donn{\'e} par une repr{\'e}sentation
alg{\'e}brique de $G$). Ces repr{\'e}sentations sont de la forme
$$\bigotimes_{\tau\in J_F}\Sym^{n_\tau}\otimes \Det^{m_\tau}.$$ 
 Une telle repr{\'e}sentation ne d{\'e}finit un syst{\`e}me local sur $M^{\an}$ 
que si le centre de $\Gamma$ agit trivialement. Cette condition
{\'e}quivaut {\`a} la condition d'alg{\'e}bricit{\'e} de Clozel (\cite{clo} Sect.1.2.3):
 
\begin{defin}\label{alg} Un poids $\kappa \in \Z[J_F]$ est dit alg{\'e}brique si 
ses coordonn{\'e}es $k_\tau$ sont sup{\'e}rieures ou {\'e}gales {\`a} $2$ et sont de
m{\^e}me parit{\'e}.
On pose alors $k_0=\max\{k_\tau |\tau \in J_F\}$,  
$m_\tau=\frac{k_0-k_\tau}{2}\in\n$,
$t=\sum_{\tau\in J_F}\tau$ et $n_\tau=k_\tau-2\geq 0$ 
($n=\kappa-2t$ et $\kappa+2m=k_0t$).
\end{defin}

Pour toute fonction 
$f:\mathfrak{H}_F\rightarrow \C$ et  pour tout $\gamma\in G_{\Q}^+$, on pose :
$$ f|_\kappa \gamma(z)=\nu(\gamma)^{\kappa+m-t}
j(\gamma,z)^{-\kappa} f(\gamma z).$$

Consid{\'e}rons l'espace:
$$G_\kappa(\mathfrak{c},\mathfrak{n})^{\an }=\Big\{f:\mathfrak{H}_F\rightarrow \C
\enspace \Big| \enspace \forall \gamma\in\Gamma, 
f|_\kappa\gamma=f\enspace\mathrm{et} \enspace 
f\enspace\mathrm{ holomorphe} \enspace\mathrm{sur} \enspace 
\mathfrak{H}_F^{*}\Big\}.$$

On appelle cet espace l'espace des formes modulaires holomorphes de poids $\kappa$ et groupe de niveau
$\Gamma$. Il est isomorphe {\`a} l'espace des sections globales du fibr{\'e} analytique 
 $\underline{\omega}^\kappa\otimes \underline{\nu}^{-n_0t/2}$ 
sur $M^{\an}$ associ{\'e} au cocyle 
$\gamma\mapsto \nu(\gamma)^{-n_0t/2}j(\gamma,z)^\kappa$ avec $n_0=k_0-2$.

\begin{rque} La torsion par $\underline{\nu}^{-n_0t/2}$ induit un 
isomorphisme  d'espaces vectoriels complexes
$$\mathrm{H}^0(M^{\an},\underline{\omega}^\kappa)\cong 
\mathrm{H}^0(M^{\an},\underline{\omega}^\kappa\otimes 
\underline{\nu}^{-n_0t/2}).$$
\end{rque}

Pour chaque $f\in G_\kappa(\mathfrak{c},\mathfrak{n})^{\an }$, 
on se propose d'expliciter  la notion d'holomorphie en une pointe 
$\mathcal{C}=\gamma\infty \in \mathbb{P}^1(F)$.
La fonction $f_{\mathcal{C}}:=f|_\kappa\gamma$ est invariante 
par le groupe $\gamma^{-1}\Gamma\gamma$ et donc par 
son sous-groupe de translations 
 $\gamma^{-1}\Gamma\gamma\cap U_{\Q}\cong X^* $ (pour
le calcul de ce-dernier voir le paragraphe pr{\'e}c{\'e}dent). Par cons{\'e}quent, 
elle admet  un d{\'e}veloppement en s{\'e}rie de Fourier :
\begin{equation}\label{dev-fourier}
f_{\mathcal{C}}(z)=\sum_{\xi\in X} a_\xi 
e^{2i\pi \Tr_{F/\mathbb{Q}}(\xi z)}.
\end{equation}

 La condition d'holomorphie en la pointe $\mathcal{C}$ se lit alors :
\begin{equation}\label{holo}
a_\xi\neq 0\Rightarrow \xi\in X_+ \text{ ou }
 \xi=0.
\end{equation}

Pour tout  $(u,\epsilon)\in \mathfrak{o}_{\mathcal{C}}^\times$, il
existe $\xi_{u,\epsilon}^*\in (\mathfrak{b}^2\mathfrak{c})^*$, 
d{\'e}fini {\`a} $X^*$ pr{\`e}s, tel que
$\begin{pmatrix} u\epsilon & \xi_{u,\epsilon}^* \\ 0 & u^{-1}\end{pmatrix}\in
\gamma^{-1}\Gamma\gamma$. L'invariance de 
$f_{\mathcal{C}}$ par le groupe $\gamma^{-1}\Gamma\gamma\cap B_{\Q}$
nous donne pour tout $\xi\in X$ la relation :
\begin{equation}\label{units}
a_{u^2\epsilon\xi}=\epsilon^{\kappa+m-t}u^\kappa e^{2i\pi\Tr_{F/\mathbb{Q}}
(\xi u\xi_{u,\epsilon}^*)}a_{\xi}.
\end{equation}

\noindent{\bf Principe du $q$-d{\'e}veloppement : } Si
pour tout $\xi$ on a $a_\xi=0$, alors $f=0$.

\medskip
\noindent{\bf Principe de Koecher : } Si $F\neq \mathbb{Q}$, alors 
la condition (\ref{holo})
est toujours satisfaite.
Si $\kappa$ n'est pas parall{\`e}le, alors $a_0=0$ (pas de 
s{\'e}ries d'Eisenstein).

\smallskip
D'apr{\`e}s le (\ref{units}), pour tout
$u\in \gamma^{-1}\Gamma\gamma\cap T_{1,\Q}$ et 
$\xi\in X$, on a $a_{u^2\xi}=u^\kappa a_\xi$,
en particulier $a_0=u^\kappa  a_0$,
d'o{\`u} la deuxi{\`e}me propri{\'e}t{\'e}. 

V{\'e}rifions (\ref{holo}) par l'absurde : soient 
$\xi\in X$ et $\xi^* \in X_+^*$ 
tels que $a_\xi\neq 0$ et  $\langle\xi,\xi^*\rangle<0$. 
Alors, on peut choisir  $u\in \gamma^{-1}\Gamma\gamma\cap T_{1,\Q}$
de fa{\c c}on que la quantit{\'e}  $\langle u^2\xi,\xi^*\rangle$ 
soit arbitrairement proche de 
$-\infty$, ce qui contredit 
l'holomorphie de $f$ au point $\underline{i}\xi^*\in\mathfrak{H}_F$.
\hfill $\square$

%voir \cite{freitag} I.4.9.

\vspace{-.4cm}
\section[Compactifications toro{\"\i}dales analytiques.]{Compactifications toro{\"\i}dales analytiques}\label{toroidale}

R{\'e}f{\'e}rences :  \cite{AMRT}\cite{KKMS}.

\medskip
On a vu qu'en ajoutant {\`a} $M^{\an }$ un nombre fini de points (les $\Gamma$-pointes)
on obtient un espace analytique $M^{\an *}$, compact pour la topologie de Satake.
Il est aussi appel{\'e} compactification minimale et n'est pas lisse
si $d_F>1$, comme le montre 
un argument de topologie (voir \cite{freitag}).

Un voisinage typique de la pointe $\infty$ est de la forme 
$\mathfrak{o}_{\infty}^\times\bs q(D_{H})\subset 
\mathfrak{o}_{\infty}^\times\bs \C^{\times d}$. 
On aurait pu tenter de compactifier cette pointe en 
consid{\'e}rant l'adh{\'e}rence de $\mathfrak{o}_{\infty}^\times\bs q(D_{H})$
dans $\mathfrak{o}_{\infty}^\times\bs \C^d$. Le probl{\`e}me 
est que si $d_F>1$ le quotient de $\C^d$ par un groupe ab{\'e}lien,
ayant des points fixes isol{\'e}s, n'est jamais lisse (voir \cite{freitag} p.30).

Il est important de disposer de compactifications lisses de $M^{\an }$ 
avec diviseurs {\`a} croisements normaux {\`a} l'infini ({\it i.e.} au-dessus
des pointes). Par
exemple,  pour pouvoir donner une d{\'e}composition de Hodge de la cohomologie singuli{\`e}re
de $M^{\an }$, on doit introduire
des faisceaux coh{\'e}rents {\`a} singularit{\'e}s logarithmiques {\`a} l'infini. 
Pour obtenir une compactification lisse de $M^{\an }$, on 
utilise la th{\'e}orie des immersions toro{\"\i}dales, s'inspirant 
du fait qu'au voisinage d'une pointe, $M^{\an }$ ressemble au quotient
d'un tore par l'action d'un groupe.

\medskip
\vspace{-.4cm} \subsubsection*{Immersions toriques.}
Dans ce paragraphe on adopte les notations suivantes :

$k$  corps alg{\'e}briquement clos. 

$S\cong\mathbb{G}_m^d$ tore alg{\'e}brique sur $k$.

$X=\Hom(S,\Gm)\cong \Z^d$ groupe des caract{\`e}res de $S$. Pour $\xi\in X$ on notera
$q^\xi$ le caract{\`e}re correspondant.

$X^*=\Hom(\Gm,S)\cong \Z^d$ groupe des cocaract{\`e}res de $S$. Pour $\xi^{*}\in X^*$ on notera
$\lambda_{\xi^{*}}$ le cocaract{\`e}re correspondant.

On a un accouplement parfait  $\langle\enspace, \enspace\rangle
:X\times X^*\rightarrow \Z$.

Pour tout anneau commutatif $R$ et tout mono{\"\i}de $Q$, 
on notera $R[q^\xi;\xi\in Q]$  la $R$-alg{\`e}bre du mono{\"\i}de.

On a $S=\Spec\big(k[q^\xi;\xi\in X])$ 
et $S=\Gm\otimes X^*$.

\begin{rque}\label{additif}
Si $k=\C$ on peut identifier $\C/\Z$ et $\Gm$ par l'application 
$e^{2i\pi\cdot}$ et on a : 

\noindent$\mathrm{(i)}$ $X^*\cong \pi_1^{\mathrm{top}}(S)$.

\noindent$\mathrm{(ii)}$ $X^*_{\C}=X^*\otimes\C$  rev{\^e}tement universel de $S$.

\noindent$\mathrm{(iii)}$ $S\cong X^*_{\C}/X^*=S_c\times iX^*_{\R}$, 
o{\`u} $S_c\cong X^*_{\R}/X^*$ est le sous-groupe compact maximal de $S$.
On appelle $\ord : S \rightarrow X^*_{\R}$ 
l'application d{\'e}duite de la projection sur $iX^*_{\R}$.

\end{rque}

\begin{defin}
Une immersion torique normale (affine) de $S$, 
est une immersion ouverte de $S$ dans une 
vari{\'e}t{\'e} (=sch{\'e}ma int{\`e}gre de type fini, s{\'e}par{\'e} sur $k$) normale
(affine) munie d'une action de $S$ qui {\'e}tend l'action de $S$
sur lui-m{\^e}me.
\end{defin}

Dans la suite, on ne consid{\'e}rera que des c{\^o}nes poly{\'e}draux rationnels convexes de $X^*_{\R}$, ouverts dans l'espace vectoriel qu'ils engendrent et  stricts 
(i.e. qui
ne contiennent pas de droite); on abr{\'e}gera 
ces propri{\'e}t{\'e}s en parlant de c{\^o}nes p.r.c.o.s. 
Un tel c{\^o}ne $\sigma$ est dit {\it lisse}, si 
$\overline{\sigma}\cap X^*$ est engendr{\'e} par une partie d'une base de $X^*$.

\begin{theo}{\rm(\cite{KKMS}  Chap.I, Th{\'e}or{\`e}me 1') }
La correspondance :
$$\sigma\mapsto S_\sigma:=
\Spec\, k[q^\xi;\xi\in X\cap\check{\sigma}]$$
donne une bijection entre l'ensemble des {\it c{\^o}nes p.r.c.o.s.} de $X^*_{\R}$
 et l'ensemble des immersions toriques normales
affines de $S$. De plus $S_\sigma$ est lisse, si et seulement si, le c{\^o}ne $\sigma$ est lisse.
\end{theo}

\begin{ex}
Voici trois exemples d'immersions torique pour  $S=\Gm^2$ :

\hspace{-6mm}$\bullet$ $\overline{\sigma_1}=(1,0)\R_++(0,1)\R_+$, donne 
$\Gm^2\hookrightarrow \Spec(k[Z_1,Z_2])\cong \A^2$.

\hspace{-6mm}$\bullet$ $\overline{\sigma_2}=(1,0)\R_+$, donne $\Gm^2\hookrightarrow\Spec(k[Z_1,Z_2,Z_2^{-1}])\cong \A^1\times\Gm$.

\hspace{-6mm}$\bullet$ $\overline{\sigma_3}=(1,1)\R_++(1,-1)\R_+$, donne 

$\Gm^2\hookrightarrow 
\Spec(k[Z_1Z_2,Z_1,Z_1Z_2^{-1}])\cong \Spec(k[Z_1,Z_2,Z_3]/(Z_1Z_3-Z_2^2)).$
\end{ex}

\begin{prop}\label{KKMS2}{\rm(\cite{KKMS}  Chap.I, Th{\'e}or{\`e}me 3) }
Soient $S_{\sigma_1}$ et $S_{\sigma_2}$ deux immersions
toriques normales affines de $S$. Alors, il existe 
un morphisme $S$-{\'e}quivariant $S_{\sigma_1}\rightarrow S_{\sigma_2}$, 
si et seulement si $\sigma_1\subset \overline{\sigma_2}$.
\end{prop}

On veut maintenant d{\'e}crire le bord de $S_{\sigma}$ :
il est stratifi{\'e} en orbites sous $S$ de points {\`a} l'infini obtenus
comme des limites ``$\lim_{t\rightarrow 0}\lambda_{\xi^{*}}(t)$'',
pour $\xi^{*}\in X^*\cap \overline{\sigma}$. De mani{\`e}re rigoureuse, pour tout
$\xi^{*}\in \overline{\sigma}\cap X^*$, 
on d{\'e}finit le point $\lambda_{\xi^{*}}(0)\in S_\sigma$, par :
$$\forall \xi\in X\cap\check{\sigma}, \enspace 
q^\xi(\lambda_{\xi^{*}}(0))=
\begin{cases} 1\text{, si }   \langle\xi,\xi^{*}\rangle =0
\\ 0\text{, si }  \langle\xi,\xi^{*}\rangle >0 \end{cases}.$$

\begin{theo}{\rm(\cite{KKMS}  Chap.I, Th{\'e}or{\`e}me 2) }

\noindent{\rm (a)} Soient $\xi^{*}_1,\xi^{*}_2\in \overline{\sigma}\cap X^*$. Alors 
$\lambda_{\xi^{*}_1}(0)=\lambda_{\xi^{*}_2}(0)$, si et seulement si 
$\xi^{*}_1$ et $\xi^{*}_2$ appartiennent {\`a} l'int{\'e}rieur d'une m{\^e}me face de $\sigma$.

\noindent{\rm (b)} Chaque $S$-orbite de $S_\sigma$ contient un unique point du
type $\lambda_{\xi^{*}}(0)$, $\xi^{*}\in \overline{\sigma}\cap X^*$.

\noindent{\rm (c)} On a une bijection entre les faces de $\sigma$ et les
$S$-orbites de $S_\sigma$,  $\tau \mapsto o(\tau)$.

\noindent{\rm (d)} $\tau_1\subset\overline{\tau_2}$ si et seulement si 
$o(\tau_2)\subset \overline{o(\tau_1)}$.

\noindent{\rm (e)} $\dim(\tau)+\dim(o(\tau))=d$.
\end{theo}

Soit une face $\tau$ de $\sigma$. 
On a $o(\tau)=\Spec\Big( k[q^\xi; \xi\in X\cap \tau^\perp]\Big)$ 
et $\overline{o(\tau)}=\coprod_{\tau\subset \overline{\tau'}}o(\tau')$.
La strate $o(\tau)$ est ferm{\'e}e dans $S_\tau$ 
(donn{\'e}e par l'id{\'e}al engendr{\'e} par les $q^\xi$ tels que 
$\langle\xi,\xi^{*}\rangle >0$ pour
tout $\xi^{*}$
{\`a} l'int{\'e}rieur de $\tau$)
et $S_\tau$ est ouverte dans $S_\sigma$. De plus les strates de $S_\sigma$ contenues dans $S_\tau$
sont les strates de $S_\tau$.

\begin{defin}
Un {\it {\'e}ventail}  dans  $X^*_{\R}$
(= d{\'e}composition rationnelle partielle en c{\^o}nes poly{\'e}draux 
fortement convexes),
est la donn{\'e}e d'un ensemble $\Sigma$
de c{\^o}nes p.r.c.o.s. de
$X^*_{\R}$  deux {\`a} deux disjoints, 
tel que pour tout $\sigma\in\Sigma$ et pour toute face 
$\tau$ de $\sigma$, $\tau\in\Sigma$. L'{\'e}ventail est dit lisse si tout les c{\^o}nes qu'il contient sont lisses.
\end{defin}

Tout {\'e}ventail peut {\^e}tre raffin{\'e} en un {\'e}ventail lisse, par subdivision des c{\^o}nes.
D'apr{\`e}s la proposition \ref{KKMS2}, {\'e}tant donn{\'e} un
{\'e}ventail  $\Sigma$,
 on peut recoller les  $S_{\sigma_i}$, $i=1,2$,  
le long des $S_\tau$, $\tau$ d{\'e}signant l'int{\'e}rieur de
$\overline{\sigma}_1\cap\overline{\sigma}_2 $, et ainsi obtenir 
un sch{\'e}ma s{\'e}par{\'e}, normal, int{\`e}gre, {\it localement} de type fini
sur $k$, not{\'e} $S_\Sigma$ ou $S_{\{\sigma\}}$. Si $\Sigma$ est fini,
$S_\Sigma$ est une vari{\'e}t{\'e}.

\begin{theo}{\rm(\cite{KKMS}  Chap.I, Th{\'e}or{\`e}me 6) }

\noindent {\rm (a)} L'application $\Sigma\mapsto S_\Sigma$
donne une bijection entre les {\'e}ventails de $X^*_{\R}$
et les immersions toriques normales de $S$.

\noindent{\rm (b)} L'application $\sigma\mapsto S_{\sigma}$
donne une bijection entre les faces $\sigma$ et les ouverts affines
$S$-invariants de $S_\Sigma$.

\noindent{\rm (c)}  L'application $\sigma\mapsto O^{\sigma}$:=
l'unique orbite ferm{\'e}e de $S_\sigma$, est une bijection 
entre les faces $\sigma$ et les $S$-orbites de $S_\Sigma$.
De plus $\tau\subset\overline{\sigma}$, 
si et seulement si, $O^{\sigma}\subset \overline{O^\tau}$.
\end{theo}

\begin{prop}{\rm(\cite{KKMS}  Chap.I, Th{\'e}or{\`e}mes 7 et 8) }
Soient $S_\Sigma$ et $S_{\Sigma'}$ deux immersions
toriques normales de $S$. Alors, il existe 
un morphisme $S$-{\'e}quivariant $S_\Sigma
\rightarrow S_{\Sigma'}$, 
si et seulement si, $\Sigma\subset \Sigma'$.
De plus la fl{\`e}che $S_\Sigma
\rightarrow S_{\Sigma'}$ est propre, 
si et seulement si, $\bigcup_{\sigma\in\Sigma} 
\sigma=\bigcup_{\sigma\in\Sigma'} \sigma$.
\end{prop}

\begin{rque}
Si $k$ est un anneau (en particulier si $k=\Z$) la construction 
qui {\`a} $\Sigma$ associe $S_\Sigma$ 
reste inchang{\'e}e. En revanche, on n'obtient pas toutes les 
immersions toriques de cette mani{\`e}re-l{\`a}.
\end{rque}

\medskip
\vspace{-.4cm} 
\subsubsection*{Carte locale pour une  pointe de $M^{\an }$.}

Soit  une $\Gamma$-pointe $\mathcal{C}=\gamma\infty$.
On a vu dans la partie  \ref{vahb} qu'un 
syst{\`e}me de voisinages de $\mathcal{C}$ dans $M^{\an }$
est donn{\'e}, pour  $H>0$,  
par les $B_{\Gamma,\mathcal{C}}\bs \gamma W_{H}
\cong \mathfrak{o}_{\mathcal{C}}^\times
\bs D_{\gamma,H}$,
o{\`u} $D_{\gamma,H}=\gamma^{-1}\Gamma\gamma\cap U_{\R}\bs W_{H}=
X^*\bs W_{H}$ et o{\`u} 
 l'action de  $\mathfrak{o}_{\mathcal{C}}^\times$ sur 
$D_{\gamma,H}$ est donn{\'e}e par 
$(u,\epsilon)\cdot z=\phi(u^2\epsilon)z+\phi(u\xi_{u,\epsilon}^*)$ (voir 
le paragraphe qui pr{\'e}c{\`e}de la d{\'e}finition \ref{vahb-tate}).

Notons qu'avec les notations  de la remarque \ref{additif},
 $X^*_{\R}=F\otimes\R$, \\ $S_{\mathcal{C}}\overset{\sim}{\underset{q}{\leftarrow}}
\phi(X^*)\bs F\otimes\C$,
$S_{\mathcal{C},c}\overset{\sim}{\underset{q}{\leftarrow}}
 F\otimes\R/\phi(X^*)$
et $\ord:\phi(X^*)\bs F\otimes\C\rightarrow
F\otimes\R$ est  l'application ``partie imaginaire''. 
On a aussi $X^*\bs \mathfrak{H}_F=
\ord^{-1}((F\otimes\R)_+)$ et 
$X^*\bs W_{H}=\ord^{-1}\{y\in (F\otimes\R)_+ \enspace |
\enspace \prod_\tau y_\tau > H\}$.

L'exponentielle donne une injection $q : D_{\gamma,H}
\hookrightarrow S_{\mathcal{C}}$ et 
l'action de  $\mathfrak{o}_{\mathcal{C}}^\times$ 
s'{\'e}tend en une action 
sur le tore complexe $S_{\mathcal{C}}=\Gm\otimes X^*$ 
par $(u,\epsilon)\cdot q_z=q_z^{u^2\epsilon}q_{\xi_{u,\epsilon}^*}^{u}$.

 L'action de  $\mathfrak{o}_{\mathcal{C}}^\times$ sur 
$S_{\mathcal{C}}$ tout entier, 
n'est pas libre (l'{\'e}l{\'e}ment unit{\'e} de $S_{\mathcal{C}}$  est fixe
par cette action). Un autre probl{\`e}me est pos{\'e} par le   centre 
de $\mathfrak{o}_{\mathcal{C}}^\times$ : 
$$\mathfrak{o}_{\mathcal{C},Z}^\times=
\{ (u,\epsilon)\in \mathfrak{o}_{\mathcal{C}}^\times \quad | 
\quad \epsilon=u^{-2} \}.$$

\begin{lemme}
$\mathrm{(i)}$ Le groupe $\mathfrak{o}_{\mathcal{C},Z}^\times$ agit trivialement 
sur $S_{\mathcal{C}}$. 

$\mathrm{(ii)}$ Sous l'hypoth{\`e}se {\bf (NT)} 
$\mathfrak{o}_{\mathcal{C}}^\times/\mathfrak{o}_{\mathcal{C},Z}^\times$ 
agit librement et discontinuement sur  $q(D_{\gamma,H})$.
\end{lemme}
Un calcul direct montre que si $\epsilon=u^{-2}$, alors 
$\xi_{u,\epsilon}^*\in X^*$, d'o{\`u} le (i). Le (ii) d{\'e}coule 
du fait que sous l'hypoth{\`e}se {\bf (NT)} on a 
  $-1\notin \mathfrak{o}_{\mathcal{C}}^\times$.\hfill $\square$

\medskip
 On veut ajouter {\`a} $S_{\mathcal{C}}$ une 
fronti{\`e}re analytique $\mathcal{E}$ de fa{\c c}on que 
l'action de $\mathfrak{o}_{\mathcal{C}}^\times/
\mathfrak{o}_{\mathcal{C},Z}^\times$  sur $S_{\mathcal{C}}$ 
se prolonge en une action libre et discontinue sur $\mathcal{E}$.
Alors le quotient par $\mathfrak{o}_{\mathcal{C}}^\times$ de 
 l'adh{\'e}rence $\overline{q(D_{\gamma,H})}$  
de $q(D_{\gamma,H})$  dans $S_{\mathcal{C}}\cup \mathcal{E}$
sera notre carte locale pour la  compactification de
la pointe $\mathcal{C}$.

Pour ce faire  on consid{\`e}re  un {\'e}ventail  $\Sigma^\mathcal{C}$ de 
$X^*_{\R+}=\{0\}\cup (F\otimes\R)_+$ qui est complet 
(i.e. tel que $\bigcup_{\sigma\in\Sigma}
\sigma=X^*_{\R+}$), stable pour
l'action de $\mathfrak{o}^\times$ et qui 
contient un nombre fini d'{\'e}l{\'e}ments modulo cette action.
L'existence d'une telle d{\'e}composition d{\'e}coule
du th{\'e}or{\`e}me des unit{\'e}s de Dirichlet ($\mathfrak{o}^{\times 2}\cong\Z^{d-1}$).
En effet, il suffit de d{\'e}composer en cellules 
$\phi(X^*)\cap \{y\in X^*_{\R+}\enspace |\enspace 
\prod_{\tau}y_{\tau}=\min_{\xi^*\in
X^*\bs \{0\}}\N(\xi^*)\}
\overset{\sim}{\underset{\mathrm{exp}}{\longleftarrow}}
\Z^{d-1}$ et prendre chaque cellule comme base d'un c{\^o}ne.

Soit $S_{\mathcal{C}}\hookrightarrow S_{\Sigma^\mathcal{C}}$ 
l'immersion torique correspondante, avec 
action {\'e}quivariante de $\mathfrak{o}^\times$. 
Soit $\mathcal{E}=S_{\Sigma^\mathcal{C}}\bs S_{\mathcal{C}}$.

Quitte {\`a} raffiner notre d{\'e}composition (en subdivisant un
c{\^o}ne et subdivisant les autres c{\^o}nes de mani{\`e}re $\mathfrak{o}^{\times}$-
{\'e}quivariante) on peut toujours supposer $S_{\Sigma^\mathcal{C}}$ lisse.

Soit $\overline{q(D_{\gamma,H})}$ 
l'adh{\'e}rence de $q(D_{\gamma,H})$ dans $S_{\Sigma^\mathcal{C}}$.
On voit alors ais{\'e}ment que

\begin{prop}
$\mathrm{(i)}$ On a $\overline{q(D_{\gamma,H})}=q(D_{\gamma,H})\cup \mathcal{E}$.

$\mathrm{(ii)}$ Le groupe $\mathfrak{o}_{\mathcal{C},Z}^\times$ 
agit trivialement sur $\overline{q(D_{\gamma,H})}$. Le groupe  
$\mathfrak{o}_{\mathcal{C}}^\times/\mathfrak{o}_{\mathcal{C},Z}^\times$ 
agit librement et discontinument    sur $\overline{q(D_{\gamma,H})}$.

\end{prop}

L'espace analytique $\mathfrak{o}_{\mathcal{C}}^\times\bs 
\overline{q(D_{\gamma,H})}$ est la  carte de la pointe $\mathcal{C}$.
Pour compactifier la pointe $\mathcal{C}$ on recolle $M^{\an }$
et $\mathfrak{o}_{\mathcal{C}}^\times\bs \overline{q(D_{\gamma,H})}$
le long de $\mathfrak{o}_{\mathcal{C}}^\times\bs D_{\gamma,H}$.

\vspace{-.4cm} \subsubsection*{Recollement et compactification analytique.}

Soit $\overline{M^{\an }}=M^{\an }_\Sigma$ 
la vari{\'e}t{\'e} analytique complexe obtenue,
par la construction du paragraphe pr{\'e}c{\'e}dent, en recollant {\`a} 
$M^{\an }$ les cartes locales pour toutes les $\Gamma$-pointes.
Dans la suite nous {\'e}crirons juste $\overline{M^{\an }}$, bien que 
tout d{\'e}pend des {\'e}ventails $\Sigma^\mathcal{C}$.

\begin{prop}
$\overline{M^{\an }}$ est une vari{\'e}t{\'e} analytique complexe normale propre,
contenant  $M^{\an }$ comme sous-vari{\'e}t{\'e} ouverte dense. Elle  est lisse si tous les
{\'e}ventails $\Sigma_{\mathcal{C}}$ sont lisses.  
On a un  morphisme de vari{\'e}t{\'e}s analytiques 
$\pi: \overline{M^{\an }}\rightarrow M^{\an *}$  qui est un 
isomorphisme au-dessus de $M^{\an }$.
\end{prop}

\noindent{\bf D{\'e}monstration : } 
Pour d{\'e}montrer la propret{\'e} de $\overline{M^{\an }}$ nous 
utiliserons le crit{\`e}re de compacit{\'e} s{\'e}quentielle. Soit une 
suite de points $z_j\in \overline{M^{\an }}$. Comme $M^{\an }$
est ouvert dense dans $\overline{M^{\an }}$, il suffit de consid{\'e}rer le
cas o{\`u} $z_j\in M^{\an }$ (argument d'extraction diagonale). Puisque l'on sait d{\'e}j{\`a}
que $M^{\an *}$ est compact, on peut supposer que la
suite $\pi(z_j)$ converge vers une pointe $\mathcal{C}$ de 
$M^{\an *}$. Dans ce cas, pour $j$ assez grand, $z_j$
appartient {\`a} $D_{\gamma,H}$. Comme $\Sigma^{\mathcal{C}}$  poss{\`e}de  un 
nombre fini de c{\^o}nes modulo l'action de $\mathfrak{o}_{\mathcal{C}}^\times$,
on peut supposer qu'il existe un c{\^o}ne $\sigma\in \Sigma^{\mathcal{C}}$, tel que 
pour tout $j$, $q(z_j)$ appartient {\`a} $S_{\mathcal{C},\sigma}$.

Montrons alors qu'il existe une suite 
extraite de la suite $q(z_j)$ qui converge vers un point de 
$S_{\mathcal{C},\sigma}$. Consid{\'e}rons la suite $y_j=\ord(q(z_j))\in \sigma$.
On a $y_{j,\tau}> 0$ et $\lim_{j\rightarrow \infty}\prod_{\tau}y_{j,\tau}=\infty$.
Si l'on d{\'e}compose les $y_j$ dans une base de $\sigma$, on trouve ais{\'e}ment 
qu'au moins une coordonn{\'e} tend vers $+\infty$. 
D'apr{\`e}s la description de la topologie de $S_{\mathcal{C},\sigma}$, 
donn{\'e}e dans \cite{AMRT}, il est clair que l'on peut extraire de 
$q(z_j)$ une sous-suite convergente dans $S_{\mathcal{C},\sigma}$. 
\hfill $\square$

\vspace{-.4cm}

\section[Vari{\'e}t{\'e}s et formes de Hilbert arithm{\'e}tiques.]{Vari{\'e}t{\'e}s et formes de Hilbert arithm{\'e}tiques}\label{fmh}

\vspace{-.2cm} 
\subsubsection*{L'espace de modules  de Hilbert-Blumenthal.}

Soit $\mathfrak{c}$ un id{\'e}al de $F$, muni de sa
positivit{\'e} naturelle $\mathfrak{c}_+=\mathfrak{c}\cap(F\otimes\R)_+$. 
Posons  $\Delta=\N(\mathfrak{dn})=\Delta_F\N(\mathfrak{n})$.

On a un foncteur contravariant 
$\underline{\mathcal{M}}^1$  (resp. $\underline{\mathcal{M}}$)  de la
cat{\'e}gorie des   $\Z[{\frac{1}{\N(\mathfrak{n})}}]$-sch{\'e}mas  vers celle des ensembles,  
qui {\`a} un sch{\'e}ma $S$ associe 
l'ensemble des quadruplets $(A,\iota,\lambda,\alpha)/S$ 
(resp. $(A,\iota,\overline{\lambda},\alpha)/S$)
modulo isomorphisme, o{\`u} 
$(A,\iota)$ est une VAHB de dimension relative $d$,
$\lambda$ est une $\mathfrak{c}$-polarisation 
(resp. $\overline{\lambda}$ est un classe de $\mathfrak{c}$-polarisations;
voir D{\'e}f.\ref{classe}) sur $A$ et 
 $\alpha: (\mathfrak{o}/\mathfrak{n})(1)\hookrightarrow A[\mathfrak{n}]$ 
est une $\mu_{\mathfrak{n}}$-structure de niveau.

\vspace{-.4cm} 
\begin{theo}\cite{GIT}\cite{Za} \label{git}
Le foncteur $\underline{\mathcal{M}}^1$ est 
repr{\'e}sentable par un sch{\'e}ma quasi-projectif $M^1$ sur 
$\Z[{\frac{1}{\N(\mathfrak{n})}}]$ 
muni d'un quadruplet universel 
$(\mathcal{A},\iota,\lambda,\alpha)$. Le sch{\'e}ma $M^1$ est lisse au-dessus de 
 $\Z[{\frac{1}{\Delta}}]$. De plus $M^1(\C)\cong M^{1,\an }$ et donc $M^1$ est 
g{\'e}om{\'e}triquement connexe.
\end{theo}

\vspace{-.2cm} 
Soit  $f:\mathcal{A}\rightarrow M^1$ la projection canonique. On pose
$\underline{\omega}=\underline{\omega}_{\mathcal{A}/M^1}=
f_* \Omega^1_{\mathcal{A}/M^1}$ et $\mathcal{H}^1_{\dR}=
\mathcal{H}^1_{\dR}(\mathcal{A}/M^1)=R^1
f_* \Omega^\bullet_{\mathcal{A}/M^1}$. Au-dessus
de $\Z[{\frac{1}{\Delta}}]$ on a localement pour la topologie de
Zariski $\underline{\omega}\cong \mathfrak{o}\otimes\mathcal{O}_{M^1}$
et $\mathcal{H}^1_{\dR}\cong L_0 \otimes\mathcal{O}_{M^1}$, 
o{\`u} $L_0=\mathfrak{o}\oplus \mathfrak{c}^*$.

\begin{cor}
Le foncteur $\underline{\mathcal{M}}$ admet un sch{\'e}ma de modules 
grossier $M$  sur $\Z[{\frac{1}{\N(\mathfrak{n})}}]$  
quasi-projectif  et lisse au-dessus de 
 $\Z[{\frac{1}{\Delta}}]$. Le sch{\'e}ma $M$ est le quotient 
de $M^1$ par le groupe fini $\mathfrak{o}_{D+}^\times/
(\mathfrak{o}_{D+}^\times\cap\mathfrak{o}_{\mathfrak{n}}^{\times 2})$
qui agit proprement et librement par 
$$[\epsilon]:(\mathcal{A},\iota,\lambda,\alpha)/S \mapsto 
(\mathcal{A},\iota,\epsilon\lambda,\alpha)/S .$$
\end{cor}
Il est important de noter pour la suite que 
les automorphismes $[\epsilon]$ de $M^1$ d{\'e}finis dans le  
corollaire se prolongent en une action du groupe 
 $\mathfrak{o}_{D+}^\times/
(\mathfrak{o}_{D+}^\times\cap\mathfrak{o}_{\mathfrak{n}}^{\times 2})$
sur les  fibr{\'e}s  $\underline{\omega}$
et $\mathcal{H}^1_{\dR}$. L'action sur $\underline{\omega}$
 est donn{\'e}e  par la formule
$s\mapsto \epsilon^{-1/2}[\epsilon]^*s$, o{\`u}
$s$  est une section de $\underline{\omega}$.  
 L'action sur $\mathcal{H}^1_{\dR}$ vient de celle
sur le complexe $Rf_*\Omega_{\mathcal{A}/M^1}^\bullet$.
Ces  actions sont d{\'e}finis sur 
l'anneau des entiers  du corps de 
nombres $F''=F(\sqrt{\epsilon},\epsilon\in \mathfrak{o}_{D+}^\times)$.

Par quotient, on peut  d{\'e}finir des 
fibr{\'e}s encore not{\'e}s $\underline{\omega}$
et $\mathcal{H}^1_{\dR}$ sur $M$.  Au-dessus
de $\Z[{\frac{1}{\Delta}}]$ on a encore localement pour la topologie de
Zariski $\underline{\omega}\cong \mathfrak{o}\otimes\mathcal{O}_{M}$
et $\mathcal{H}^1_{\dR}\cong L_0 \otimes\mathcal{O}_{M}$, 
o{\`u} $L_0=\mathfrak{o}\oplus \mathfrak{c}^*$.

Pour chaque $\mu\in\mathfrak{c}_+$, on note $\mathcal{L}_\mu$ le 
faisceau inversible ample sur  $\mathcal{A}$ obtenu comme image inverse 
du fibr{\'e} de Poincar{\'e} sur $\mathcal{A}\times \mathcal{A}^t$
par le morphisme
$(\mathrm{id}_{\mathcal{A}},
\lambda\circ(\mathrm{id}_{\mathcal{A}}\otimes\mu))$.

\vspace{-.4cm} 
\subsubsection*{Formes modulaires de Hilbert arithm{\'e}tiques.}

Consid{\'e}rons le sch{\'e}ma en groupes  $\mathcal{T}_1=\Res^{\mathfrak{o}}_{\Z}\Gm$
qui est un mod{\`e}le entier du tore $F^\times=\Res^F_{\Q}\Gm$.
Ce n'est un tore que sur $\Z[\frac{1}{\Delta_F}]$ comme le 
montre l'exemple suivant :

\vspace{-.4cm} 
\begin{ex} Soit $F=\Q(\sqrt{D})$, avec $D\equiv 3\pmod 4$
sans facteurs carr{\'e}s. Alors 
$\mathcal{T}_1=\Spec(\Z[X,Y][\frac{1}{X^2-DY^2}])$  
et pour $p>2$   premier on a  :
$\mathcal{T}_1(\F_p)= \F_p^\times\times\F_p^\times$, si
$\left(\frac{D}{p}\right)=1$; $\mathcal{T}_1(\F_p)=\F_{p^2}^\times$, si
$\left(\frac{D}{p}\right)=-1$; $\mathcal{T}_1(\F_p)= \F_p^\times 
\times\F_p^+$, si $\left(\frac{D}{p}\right)=0$.
\end{ex}

On suppose d{\'e}sormais $F\neq \Q$ et on se place au-dessus de
$\Z[\frac{1}{\Delta}]$.
Consid{\'e}rons le faisceau $\underline{\mathfrak{M}}=
\underline{\Isom}_{\mathfrak{o}\otimes\mathcal{O}_M}
(\mathfrak{o}\otimes\mathcal{O}_M,\underline{\omega})$.
C'est un $\mathcal{T}_1$-torseur Zariski sur $M$.

Comme $\mathcal{T}_1$ est affine sur $M$, le faisceau 
$\underline{\mathfrak{M}}$ est repr{\'e}sentable par un 
sch{\'e}ma $f:\mathfrak{M}\rightarrow M$  (voir \cite{milne} III.4 Th{\'e}or{\`e}me 4.3). 
En particulier on a un isomorphisme $\mathcal{T}_1\underset{M}{\times}
\mathfrak{M}\cong\mathfrak{M}\underset{M}{\times}\mathfrak{M}$, \enspace
$(t,x)\mapsto (tx,x)$.

Sur le  sch{\'e}ma de modules fin $M^1$ le sch{\'e}ma correspondant 
$\mathfrak{M}^1$ repr{\'e}sente 
le foncteur :
$\underline{\mathfrak{M}}^1:\Z[{\frac{1}{\Delta}}]$-$\Sch
\rightarrow \Ens$,
qui {\`a} un $\Z[{\frac{1}{\Delta}}]$-sch{\'e}ma $S$ associe 
l'ensemble des quintuplets $(A,\iota,\lambda,\alpha,\omega)$ 
modulo isomorphisme, o{\`u} $(A,\iota,\lambda,\alpha)$ est une VAHB,
 comme plus haut, et o{\`u} $\omega$ est un isomorphisme 
de $\mathfrak{o}$-fibr{\'e}s inversibles $\omega:\mathfrak{o}\otimes\mathcal{O}_S
\cong \underline{\omega}$. La fl{\`e}che d'oubli fait de 
$\underline{\mathfrak{M}}^1$ un faisceau Zariski sur $M^1$.

\smallskip

Pour la  d{\'e}finition de l'espace des formes modulaires de Hilbert,
 nous suivons de pr{\`e}s le paragraphe 
6.8 dans \cite{rapoport}, r{\'e}dig{\'e} par P. Deligne.

Soit $\kappa\in \Z[J_F]=X(\mathcal{T}_1)$ un poids et soit 
$F'$ un corps de nombres, contenant $F''$ ainsi que les valeurs du 
caract{\`e}re $\kappa : F^\times \rightarrow \C^\times$.

\smallskip 
Si $D=\Gm$, on peut prendre, par exemple, $F'=\Q$ et 
 $\kappa=kt$ (poids parall{\`e}le), ou bien $F'=F^{\gal }$
 $\kappa\in\Z[J_F]$ poids quelconque.

\smallskip

Soit $\mathcal{O}'$ l'anneau des entiers de  $F'$. 
Le morphisme de groupes alg{\'e}briques $\kappa : \Res_{\Q}^{F}\Gm \rightarrow
\Res_{\Q}^{F'}\Gm$, se prolonge en un morphisme 
$\Res_{\Z}^{\mathfrak{o}}\Gm \rightarrow \Res_{\Z}^{\mathcal{O}'}\Gm$, qui 
{\'e}quivaut (par la formule d'adjonction) {\`a} un morphisme de groupes alg{\'e}briques
sur $\mathcal{O}'$,
 $\Res_{\Z}^{\mathfrak{o}}\Gm\times\Spec(\mathcal{O}')
\rightarrow \Gm\times\Spec(\mathcal{O}')$, not{\'e} encore $\kappa$.

\medskip

Pour tout $\Z[\frac{1}{\Delta}]$-sch{\'e}ma $Y$, on pose 
$Y'=Y\times \Spec(\mathcal{O}'[\frac{1}{\Delta}])$. On a ainsi un 
$\mathcal{T}'_1$-torseur  $f':\mathfrak{M}'\rightarrow M'$. 
Le tore d{\'e}ploy{\'e} $\mathcal{T}'_1$ agit sur $f'_* \mathcal{O}_{\mathfrak{M}'}$; la composante $-\kappa$-isotypique
 $(f'_* \mathcal{O}_{\mathfrak{M}'})[-\kappa]$
est un faisceau inversible sur $M'$, not{\'e} 
 $\underline{\omega}^\kappa$.

\begin{defin} 1) Soit $R$ une $\Z[\frac{1}{\Delta}]$-alg{\`e}bre. 
On d{\'e}finit  l'espace   $G(\mathfrak{c},\mathfrak{n};R)^{\geom}$
des formes modulaire de Hilbert de 
niveau $\Gamma$ et  {\`a} coefficients dans $R$, comme
$\mathrm{H}^0(\mathfrak{M}\times_{\Spec(\Z[\frac{1}{\Delta}])} \Spec(R)
,\mathcal{O}_\mathfrak{M})$.

2) Soit $R$ une $\mathcal{O}'[\frac{1}{\Delta}]$-alg{\`e}bre. 
Une forme modulaire de Hilbert arithm{\'e}tique 
de poids $\kappa$, de 
niveau $\Gamma$ et  {\`a} coefficients dans $R$, est 
une section  globale de $\underline{\omega}^\kappa$ sur
$M\times_{\Spec(\Z[\frac{1}{\Delta}])} \Spec(R)$.
On note 
$G_\kappa(\mathfrak{c},\mathfrak{n};R)^{\geom}:=
\mathrm{H}^0(M\times_{\Spec(\Z[\frac{1}{\Delta}])} \Spec(R)
,\underline{\omega}^{\kappa})$ l'espace  
de ces formes modulaire de Hilbert.
\end{defin}

\begin{rque} 
1) Le faisceau $\underline{\omega}^t$ ($t=\sum\tau$)
n'est autre que le faisceau 
$\wedge^d\underline{\omega}=\det(\underline{\omega})$ 
sur $M$, et $\underline{\omega}^{kt}$ - sa puissance $k$-i{\`e}me.
Les formes modulaires de Hilbert de poids parall{\`e}le $k\geq 1$,
s'{\'e}crivent donc $G_{kt}(\mathfrak{c},\mathfrak{n})^{\geom}=
\mathrm{H}^0(M,(\wedge^d \underline{\omega})^{\otimes k})$.

\noindent 2) Le  torseur $\mathfrak{M}$ n'est pas
trivial, car sinon pour tout $k\geq 1$ 
$(\wedge^d\underline{\omega})^{\otimes k}$ serait le fibr{\'e} trivial sur $M$,
et {\`a} fortiori sur $M^{\mathrm{\an }}$. Or, par le principe de
Koecher $\mathrm{H}^0(M^{\mathrm{\an }}, \mathcal{O}_{M^{\mathrm{\an }}})=
\mathrm{H}^0(\overline{M^{\mathrm{\an }}}, 
\mathcal{O}_{\overline{M^{\mathrm{\an }}}})=\C$, 
ce qui contredirait l'existence de formes modulaires de Hilbert cuspidales
non-nulles en poids $kt$.

\noindent 3) Si $F'\supset F^{\gal }$, on a
$\mathcal{T}'_1:=\mathcal{T}_1 \times \Spec(\mathcal{O}'[\frac{1}{\Delta}])\cong
\Gm^{J_F}\times \Spec(\mathcal{O}'[\frac{1}{\Delta}])$
et par le th{\'e}or{\`e}me de diagonalisabilit{\'e} des tores 
\cite{jantzen}, on a : 
$$f'_*\mathcal{O}_\mathfrak{M'}=
\bigoplus_{\kappa\in X(\mathcal{T}_1)}\underline{\omega}^\kappa \enspace,
\enspace \mathrm{H}^0(\mathfrak{M'},\mathcal{O}_\mathfrak{M'})=
\bigoplus_{\kappa\in X(\mathcal{T}_1)}\mathrm{H}^0(M',\underline{\omega}^\kappa).$$

Par ailleurs, l'action de $\mathfrak{o}$ permet de d{\'e}composer 
$\underline{\omega}=\Lie(\mathcal{A}'/M')^\vee\cong
\mathfrak{o}\otimes\mathcal{O}_{M'} $ en somme directe de fibr{\'e}s inversibles 
$\underline{\omega}^\tau$ sur $M'$, index{\'e}s par les diff{\'e}rents 
plongements $\tau$ de $\mathfrak{o}$ dans $\mathcal{O}'$.
On a $\underline{\omega}^\kappa=\otimes_\tau
(\underline{\omega}^\tau)^{\otimes k_\tau}$.

\noindent 4) Si $R$ est une $\mathcal{O}'[\frac{1}{\Delta}]$-alg{\`e}bre, avec 
$F'\supset F^{\gal }$, on a :
$$G(\mathfrak{c},\mathfrak{n};R)^{\geom}=
\bigoplus_{\kappa\in X(\mathcal{T}_1)}
G_\kappa(\mathfrak{c},\mathfrak{n};R)^{\geom}.$$
\end{rque}

\vspace{-.4cm} 
\subsubsection*{Constructions de fibr{\'e}s automorphes.}
Dans la partie \ref{vahb} on a introduit 
les formes modulaires de Hilbert classiques 
comme des sections globales de certains fibr{\'e}s de 
formes diff{\'e}rentielles holomorphes  sur 
$M^{\mathrm{\an }}$. Dans ce paragraphe nous donnons
 des constructions   de  fibr{\'e}s sur $M^{\mathrm{\an }}$ et   $M$, {\`a}
partir de repr{\'e}sentations de certains groupes.
Ces fibr{\'e}s  serviront {\`a} 
 red{\'e}finir et {\'e}tudier les formes modulaires de Hilbert arithm{\'e}tiques.

Soit   un poids alg{\'e}brique $\kappa$  et 
$n, m\in \Z[J_F]$ comme dans la  d{\'e}finition \ref{alg}. On notera 
$V_n$ la repr{\'e}sentation alg{\'e}brique de $G$ donn{\'e}e par
\begin{equation}\label{Vn}
V_n=\bigotimes_{\tau\in J_F}\Sym^{n_\tau}\otimes \Det^{m_\tau}.
\end{equation}

\medskip
\noindent$\bullet$ Consid{\'e}rons le rev{\^e}tement universel
 $u: \mathfrak{H}_F\rightarrow M^{\mathrm{\an }}$.
 Il est bien connu que l'on a une {\'e}quivalence
de cat{\'e}gories entre les repr{\'e}sentations de $\Gamma$ sur des
$K$-vectoriels de dimension finie, qui sont triviales sur le centre,
 et les syst{\`e}mes locaux en $K$-vectoriels $\bV^{\an}$ sur 
$M^{\mathrm{\an }}$, qui {\`a}
un $K$-vectoriel de dimension finie $V$ muni d'une telle  
action de $\Gamma$  associe le syst{\`e}me local $\bV^{\an}$
des sections continues de $\Gamma\bs(\mathfrak{H}_F\times V)
\rightarrow M^{\mathrm{\an }}$ ($V$ {\'e}tant muni de la topologie discr{\`e}te).

\begin{defin}
On note $\bV^{\an}_n$ le syst{\`e}me local associ{\'e} {\`a} 
repr{\'e}sentation $V_n$.
\end{defin}

\noindent$\bullet$ Une autre construction de fibr{\'e}s est donn{\'e}e par
la tour promodulaire $\widetilde{M_{\Q}}\rightarrow M_{\Q}$, o{\`u} 
$\widetilde{M_{\Q}}=\underset{r\geq 1}{\mathrm{limproj}}\enspace
M(\mathfrak{c},\mathfrak{n}r)_{\Q}$. On a une suite exacte 
$$1\rightarrow \pi_1(M_{\Q},x)^{\mathrm{mod,geom}}
\rightarrow \pi_1(M_{\Q},x)^{\mathrm{mod}}=\Gal(\widetilde{M_{\Q}}/M_{\Q}) \rightarrow 
\Gal(\Q^{\mathrm{ab}}/\Q) \rightarrow 1.$$
De plus, le groupe $\pi_1(M_{\Q},x)^{\mathrm{mod}}$
 est un sous-groupe ouvert de $\GL_2(\widehat{\mathfrak{o}})$; la
projection sur la $p$-composante fournit un morphisme continu canonique
$$\pi_1(M_{\Q},x)^{\mathrm{mod}}\rightarrow \GL_2(\mathfrak{o}\otimes \Z_p).$$ 
 
 On a  donc un foncteur de la cat{\'e}gorie des repr{\'e}sentations alg{\'e}briques de  $G$, 
 vers celle des faisceaux lisses sur $M_{\Q}$. Ce foncteur associe {\`a}
la repr{\'e}sentation $V$ le faisceau $\bV$ des sections continues de 
$\pi_1(M_{\Q},x)^{\mathrm{mod}}\bs(\widetilde{M_{\Q}}\times V)
\rightarrow M_{\Q}$.

\begin{defin}
On note $\bV_n$ le faisceau lisse sur $M_{\Q}$ associ{\'e} {\`a} $V_n$.
\end{defin}

\noindent$\bullet$ 
Dans le cadre arithm{\'e}tique le rev{\^e}tement universel $\mathfrak{H}_F$
de la premi{\`e}re construction est remplac{\'e} par
le torseur $\mathfrak{M'}\overset{\mathcal{T}'_1}{\rightarrow}M'$ du 
paragraphe pr{\'e}c{\'e}dent. On a un foncteur
de la cat{\'e}gorie des repr{\'e}sentations alg{\'e}briques du 
 $\mathcal{O'}[{\frac{1}{\Delta}}]$-sch{\'e}ma en groupes $\mathcal{T}'_1$, 
vers celle des fibr{\'e}s d{\'e}composables en fibr{\'e}s inversibles sur 
$M'$, qui {\`a} $W_1$ associe le produit contract{\'e}
 $\mathfrak{M'}\stackrel{\mathcal{T}'_1}{\times} W_1=:\mathcal{W}_1$,
d{\'e}fini comme le quotient par la relation d'{\'e}quivalence
$(mt,w)\sim (m,tw)$ pour $m\in \mathfrak{M'}$, $t\in\mathcal{T}'_1$ et $w\in W_1$.

\begin{rque} Pour chaque $\kappa\in\Z[J_F]=X(\mathcal{T}_1)$, notons 
$W_{1,\kappa}$ 
la $\mathcal{O'}[{\frac{1}{\Delta}}]$-repr{\'e}sentation de $\mathcal{T}'_1$ associ{\'e}e {\`a} 
$\kappa$. 
On a 
$\mathcal{W}_{1,\kappa}=\underline{\omega}^{-\kappa}$. On peut ainsi red{\'e}finir 
$G_\kappa(\mathfrak{c},\mathfrak{n})^{\geom}$ 
comme $\mathrm{H}^0(M',\mathcal{W}_{1,-\kappa})$.
\end{rque}

\noindent$\bullet$ 
On suppose que $D$ a un mod{\`e}le entier $\mathcal{D}$ sur $\Z[\frac{1}{\Delta_F}]$ (c'est le cas pour $D=\Gm$ ou $\Res^F_{\Q}\Gm$).
Rappelons que le $\mathfrak{o}$-fibr{\'e} projectif de rang deux $\mathcal{H}^1_{\mathrm{\dR}}=R^1f_*\Omega^\bullet_{\mathcal{A}'/M'}$ 
est muni d'un accouplement parfait symplectique $\mathfrak{o}$-lin{\'e}aire associ{\'e} au choix d'un repr{\'e}sentant $\lambda$
de la classe de $\mathfrak{c}$-polarisations universelle $\overline{\lambda}=
\mathfrak{o}_{D+}^\times\cdot\lambda$.
On d{\'e}finit alors le $\mathcal{D}$-torseur 
$$\mathfrak{M}_{\mathcal{D}}=\mathrm{Isom}^{\mathcal{D}}_{\mathfrak{o}\otimes\mathcal{O}_{M}}
(\mathfrak{o}\otimes\mathcal{O}_{M},\wedge^2_{\mathfrak{o}\otimes\mathcal{O}_{M}}\mathcal{H}^1_{\dR})$$
 au-dessus de $M$, dont les $S$-points sont ceux de
$\mathrm{Isom}_{\mathfrak{o}\otimes\mathcal{O}_{M}}
(\mathfrak{o}\otimes\mathcal{O}_{M},\wedge^2_{\mathfrak{o}\otimes\mathcal{O}_{M}}\mathcal{H}^1_{\dR})$
induisant via $\lambda$ un {\'e}l{\'e}ment de $\mathcal{D}(\mathcal{O}_S)$ dans $(\mathfrak{o}\otimes \mathcal{O}_S)^\times$.

\noindent$\bullet$ 
On choisit pour mod{\`e}le entier du tore maximal standard $T$ de $B$ 
le sch{\'e}ma en groupes $\mathcal{T}=\mathcal{T}_1\times \mathcal{D}$.
On en d{\'e}duit un mod{\`e}le entier de $B$ dont $\mathcal{T}$ 
est tore maximal standard via
$(u,\epsilon)\in \mathcal{T}\mapsto \left(\begin{array}{cc}u\cdot\epsilon&0\\0&u^{-1}\end{array}\right)$. 
On va d{\'e}finir un ${\mathcal{B}}'$-torseur 
$\mathfrak{M}_{\mathcal{B}}'\overset{\mathcal{B'}}{\rightarrow}M'$ {\`a} l'aide du
$\mathfrak{o}$-fibr{\'e} 
$\mathcal{H}^1_{\dR}$ muni de la filtration de Hodge 
$$0 \rightarrow \underline{\omega}\rightarrow \mathcal{H}^1_{\dR}
\rightarrow \underline{\omega}^{\vee}\otimes \mathfrak{cd}^{-1}\rightarrow 0.$$

Soit $L_0=\mathfrak{o}\oplus \mathfrak{c}^*$.
On munit $L_0 \otimes\mathcal{O}_{M'}$ de la filtration canonique {\`a} deux 
crans associ{\'e}e {\`a} $\mathcal{B}'$ : $0\subset \mathfrak{c}^* \otimes\mathcal{O}_{M'}
\subset L_0 \otimes\mathcal{O}_{M'}$.
On d{\'e}finit alors
$\mathfrak{M}_{\mathcal{B}}'$ comme le produit fibr{\'e} de
$\mathrm{Isom}^{\mathrm{fil}}_{\mathfrak{o}\otimes\mathcal{O}_{M'}}
(L_0\otimes\mathcal{O}_{M'},\mathcal{H}^1_{\dR})$ et $\mathrm{Isom}^{\mathcal{D}}_{\mathfrak{o}\otimes\mathcal{O}_{M'}}
(\mathfrak{o}\otimes\mathcal{O}_{M'},\wedge^2_{\mathfrak{o}\otimes\mathcal{O}_{M'}}\mathcal{H}^1_{\dR})$
au-dessus de $\mathrm{Isom}_{\mathfrak{o}\otimes\mathcal{O}_{M'}}
(\mathfrak{o}\otimes\mathcal{O}_{M'},\wedge^2_{\mathfrak{o}\otimes\mathcal{O}_{M'}}\mathcal{H}^1_{\dR})$.
C'est un $\mathcal{B}'$-torseur sur $M'$.

Il d{\'e}finit un foncteur $\mathcal{F}_{\mathcal{B}'}$
de la cat{\'e}gorie des repr{\'e}sentations alg{\'e}briques du 
 $\mathcal{O}'[{\frac{1}{\Delta}}]$-sch{\'e}ma en groupes $\mathcal{B}'$ vers celle 
des fibr{\'e}s sur $M'$ qui sont des extensions successives de
fibr{\'e}s inversibles. 
Il est donn{\'e} par le produit contract{\'e} :
 $V \quad \mapsto \quad \mathcal{V}:=
\mathfrak{M}_{\mathcal{B}}'\stackrel{\mathcal{B}'}{\times} V$,
(c'est-{\`a}-dire le quotient par la relation 
$(\widetilde{m}b,v)\sim(\widetilde{m},bv)$,
pour $\widetilde{m}\in \mathfrak{M}_{\mathcal{B}}'$, 
$b\in\mathcal{B}'$ et $v \in V$).

\begin{defin}\label{Vfil}
On note $\mathcal{V}_n$ le fibr{\'e} filtr{\'e} sur $M'$ image de $V_n$
par  $\mathcal{F}_{\mathcal{B}'}$.
\end{defin}

\noindent$\bullet$ Si $W$ est une $\mathcal{O}'[{\frac{1}{\Delta}}]$-repr{\'e}sentation
du tore $\mathcal{T}'$ (sur 
un $\mathcal{O}'[{\frac{1}{\Delta}}]$-module libre de type fini) 
on peut la voir comme une repr{\'e}sentation de 
$\mathcal{B}'$, en faisant agir le radical unipotent
$\mathcal{U}'$ trivialement. Le foncteur $\mathcal{F}_{\mathcal{B}'}$
associe {\`a} $W$ un fibr{\'e} $\mathcal{W}$ d{\'e}composable en somme directe 
de fibr{\'e}s inversibles. 

On pourrait {\'e}galement construire $\mathcal{W}$ {\`a} l'aide du 
$\mathcal{T}'$-torseur  $\mathfrak{M}'\times_{M'}\mathfrak{M}_{\mathcal{D}}'$.

\begin{defin}\label{W}
Soient $n,m\in\Z[J_F]$ et $c\in \Z$ tels que $n+2m=ct$. 
Soit $W_{n,c}$ 
la repr{\'e}sentation irr{\'e}ductible de $\mathcal{T}'$, donn{\'e}e par
le caract{\`e}re 
$$(u,\epsilon)\in \mathcal{T}'_1\times\mathcal{D}'\mapsto u^n\epsilon^m.$$

On note $\mathcal{W}_{n,c}$ le fibr{\'e} inversible  sur $M'$ image de $W_{n,c}$
par le foncteur $\mathcal{F}_{\mathcal{B}'}$.
\end{defin}

\noindent$\bullet$ Consid{\'e}rons le mod{\`e}le entier $\mathcal{G}$ de $G$ sur 
$\Z[\frac{1}{\Delta_F}]$ par 
$$\mathcal{G}=\Res^{\mathfrak{o}}_{\Z}
\GL_2\times_{\Res^{\mathfrak{o}}_{\Z}\Gm}\mathcal{D}. $$ 
 
On introduit pour finir un $\mathcal{G'}$-torseur
$\mathfrak{M}_{\mathcal{G}'}\overset{\mathcal{G'}}{\rightarrow}M'$ {\`a} l'aide de   
$\mathcal{H}^1_{\dR}=R^1f_*\Omega^\bullet_{\mathcal{A}'/M'}$
muni de sa connexion de Gauss-Manin qui est int{\'e}grable. Plus pr{\'e}cis{\'e}ment, on munit 
$L_0 \otimes\mathcal{O}_{M'}$ de la connexion plate ${\rm Id}\otimes d$
et on pose
$$\mathfrak{M}_{\mathcal{G}}'=
\mathrm{Isom}^{\mathcal{D}}_{\mathfrak{o}\otimes\mathcal{O}_{M'}}
(L_0\otimes\mathcal{O}_{M'},\mathcal{H}^1_{\dR}).$$

Il d{\'e}finit un foncteur $\mathcal{F}_{\mathcal{G}'}$
de la cat{\'e}gorie des repr{\'e}sentations alg{\'e}briques du 
 $\mathcal{O}'[{\frac{1}{\Delta}}]$-sch{\'e}ma en groupes $\mathcal{G}'$ vers celle 
des fibr{\'e}s sur $M'$ munis d'une connexion int{\'e}grable.
Il est donn{\'e} par le produit contract{\'e} :
 $V \quad \mapsto \quad \mathcal{V}^{\nabla}:=
\mathfrak{M}_{\mathcal{G}}'\stackrel{\mathcal{G}'}{\times} V$,
(c'est-{\`a}-dire le quotient par la relation 
$(\widetilde{m}g,v)\sim(\widetilde{m},gv)$, pour $\widetilde{m}\in 
\mathfrak{M}_{\mathcal{G}}'$, $g\in\mathcal{G}'$ et $v \in V$).

\begin{defin}\label{Vnabla}
On note $\mathcal{V}_n^\nabla$ le fibr{\'e} {\`a} connexion sur $M'$ image de $V_n$
par le foncteur $\mathcal{F}_{\mathcal{G}'}$.
\end{defin}

\vspace{.2cm} 
Pour  une $\mathcal{O}'[{\frac{1}{\Delta}}]$-repr{\'e}sentation alg{\'e}brique 
$V$  de $G$, on peut comparer 
$\bV^{\an}$, $\bV$, $\mathcal{V}$ et $\mathcal{V}^\nabla$   comme suit
 
\begin{prop}
1) Sur $M'$, on a   $\mathcal{V}=\mathcal{V}^\nabla$ et

2) Sur $M^{\mathrm{\an }}$, on a
$\mathcal{V}\otimes_{\mathcal{O}_M}{\mathcal{O}_{M^{\mathrm{\an }}}}
\cong\bV^{\an}\otimes {\mathcal{O}_{M^{\mathrm{\an }}}}$ et 
$\bV^{\an}\otimes \Z_p\cong (\bV)^{\an}$.
\end{prop}

Pour la d{\'e}monstration de ce r{\'e}sultat, voir  \cite{MoTi} Sect.5.2.2, lemme 9.

\section[Compactifications  arithm{\'e}tiques de la vari{\'e}t{\'e} de Hilbert.]{Compactifications   arithm{\'e}tiques de la vari{\'e}t{\'e} de Hilbert.}

Dans cette partie nous {\'e}non{\c c}ons le r{\'e}sultat principal de \cite{dimdg}
en conservant les notations de cette r{\'e}f{\'e}rence. En particulier, nous utiliserons
la notion de  $(R,\mathfrak{n})$-composante $\mathcal{C}$ (D{\'e}f.3.2 de \cite{dimdg})
{\`a} qui sont  associ{\'e}s les  objets suivants : 
des id{\'e}aux  $\mathfrak{b}$, $\mathfrak{b'}$,
$\mathfrak{a}=\mathfrak{bc}$, $X=\mathfrak{cbb}'$ ;
une racine de l'unit{\'e}  $\zeta_{\mathcal{C}}$ d'ordre l'exposant $n$ du 
groupe $\mathfrak{b'}/\mathfrak{b}$; 
des groupes d'unit{\'e}s $\mathfrak{o}_{\mathcal{C}}^\times$,
$\mathfrak{o}_{\overline{\mathcal{C}}}^\times$, 
$\mathfrak{o}_{\mathcal{C},1}^\times$,
$\mathfrak{o}_{\overline{\mathcal{C}},1}^\times$; des sous-groupes 
$H_{\mathcal{C}}=\mathfrak{o}_{\overline{\mathcal{C}}}^\times/
\mathfrak{o}_{\mathcal{C}}^\times $,
$H_{\mathcal{C},1}=\mathfrak{o}_{\overline{\mathcal{C}},1}^\times/
\mathfrak{o}_{\mathcal{C},1}^\times $ du groupe   $(\Z/n\Z)^\times$.

\medskip
De fa{\c c}on impr{\'e}cise mais suggestive, on peut penser {\`a}  une $(R,\mathfrak{n})$-composante 
comme {\`a} une 
orbite sous le groupe de Galois d'une  $(R,\mathfrak{n})$-pointe ({\it loc. cit.} pour plus de d{\'e}tails).

\begin{defin} Un {\it {\'e}ventail $\Gamma$-admissible} 
$\Sigma=(\Sigma^{\mathcal{C}})_{\mathcal{C}}$ 
est la donn{\'e}e pour chaque $(R,\mathfrak{n})$-composante  
$\mathcal{C}$  d'un {\'e}ventail complet
$\Sigma^{\mathcal{C}}$ de $X_+^*$, 
stable par $\mathfrak{o}_{\mathcal{C}}^{\times}$ et 
contenant un nombre fini d'{\'e}l{\'e}ments modulo
cette action, de sorte que les donn{\'e}es soient compatibles aux isomorphismes 
de $(R,\mathfrak{n})$-composantes $\mathcal{C}\cong \mathcal{C}'$.
\end{defin}

On se fixe un {\'e}ventail lisse 
$\Gamma$-admissible $\Sigma=(\Sigma^{\mathcal{C}})_{\mathcal{C}}$.

\smallskip

Soit $R_{\mathcal{C}}=\Z[q^\xi;\xi\in X]$.
Soit  $S_{\mathcal{C}}=\Spec(R_{\mathcal{C}})=\Gm\otimes X^*$ le tore sur $\Z$ 
de groupe des caract{\`e}res  $X=\mathfrak{cbb}'$. 

Soit $S_{\mathcal{C}}\hookrightarrow
S_{\Sigma^{\mathcal{C}}}$, l'immersion torique associ{\'e}e. On rappelle
qu'elle est obtenue en recollant, pour $\sigma\in \Sigma^{\mathcal{C}}$, 
les immersions toriques
affines $S_{\mathcal{C}}\hookrightarrow
S_{\mathcal{C},\sigma}=\Spec({R_{\mathcal{C},\sigma}})$, o{\`u} 
$R_{\mathcal{C},\sigma}=\Z[q^\xi;\xi\in X\cap\check{\sigma}]$.
Soit $S_{\mathcal{C}, \sigma}^{\wedge}$ le compl{\'e}t{\'e} de $S_{\mathcal{C},\sigma}$
le long de $S_{\mathcal{C}, \sigma}^{\infty}:=S_{\mathcal{C}, \sigma}\bs 
S_{\mathcal{C}}$ et 
$S_{\Sigma^{\mathcal{C}}}^{\wedge}$ le compl{\'e}t{\'e} de $S_{\Sigma^{\mathcal{C}}}$
le long de $S_{\Sigma^{\mathcal{C}}}^{\infty}:=S_{\Sigma^{\mathcal{C}}}\bs
S_{\mathcal{C}}$.

\smallskip

Posons  $\overline{S}_{\mathcal{C},\sigma}=\Spec(R_{\mathcal{C},\sigma}^{\wedge})$ 
et   $\overline{S}{}_{\mathcal{C},\sigma}^{0}=
S_{\mathcal{C}}\underset{S_{\mathcal{C},\sigma}}{\times}\overline{S}_{\mathcal{C},\sigma}=
\Spec(R_{\mathcal{C},\sigma}^{\wedge}
\otimes_{R_{\mathcal{C},\sigma}} R)$.  
Si $\sigma'\subset\sigma$, on a une fl{\`e}che $\overline{S}_{\mathcal{C},\sigma'}
\rightarrow \overline{S}_{\mathcal{C},\sigma}$.

\medskip
Le th{\'e}or{\`e}me suivant est une variante pour le groupe $\Gamma$ des 
th{\'e}or{\`e}mes de Rapoport \cite{rapoport},
et Chai \cite{chai}; les modifications n{\'e}cessaires pour son {\'e}nonc{\'e} 
et sa d{\'e}monstration sont donn{\'e}es dans \cite{dimdg}.

\bigskip

\begin{theo} 
Soit $\Sigma=\{\Sigma^{\mathcal{C}}\}_{\mathcal{C}}$ 
un {\'e}ventail  $\Gamma$-admissible lisse.

$\mathrm{(i)}$ Il existe un $\Z[\frac{1}{\N(\mathfrak{n})}]$-sch{\'e}ma propre
$\overline{M}=M_{\Sigma}$  
 lisse au-dessus de $\Z[\frac{1}{\Delta}]$,   une  
 immersion ouverte $j:M\hookrightarrow  \overline{M}$
et un isomorphisme de sch{\'e}mas formels 
$$\varphi: \coprod_{(R,\mathfrak{n})\mathrm{-composantes}/\sim} 
\left(S_{\Sigma^{\mathcal{C}}}^{\wedge}/\mathfrak{o}_{\mathcal{C}}^{\times} \right)\times \textstyle
\Spec(\Z[\frac{1}{\N(\mathfrak{n})},\zeta_{\mathcal{C}}]^{H_{\mathcal{C}}})
\overset{\sim}{\longrightarrow} \overline{M}^{\wedge},$$
o{\`u} $\overline{M}^{\wedge}$ est le compl{\'e}t{\'e} formel
de $\overline{M}$ le long du diviseur {\`a} croisements normaux 
{\`a} l'infini $\overline{M}\bs M$.

$\mathrm{(ii)}$ Il  existe un unique sch{\'e}ma en groupes semi-ab{\'e}lien 
$\overline{f}:\mathfrak{G}\rightarrow \overline{M^1}$ qui prolonge 
la VAHB  universelle $f:\mathcal{A} \rightarrow M^1$.
 Ce sch{\'e}ma en groupes est muni d'une action de 
$\mathfrak{o}$ au-dessus de $\overline{M^1}$ prolongeant celle sur $\mathcal{A}$.
C'est un tore au-dessus de $\overline{M^1}\bs M^1$.

$\mathrm{(iii)}$
On a un isomorphisme de Kodaira-Spencer logarithmique :
$$\Omega_{\overline{M^1}}(\dlog\,\infty)\cong
\underline{\omega}_{\mathfrak{G}/\overline{M^1}}\otimes_{\mathfrak{o}
\otimes\mathcal{O}_{\overline{M^1}}}
(\underline{\omega}_{\mathfrak{G}/\overline{M^1}}
\otimes_{\mathfrak{o}} \mathfrak{dc}^{-1}), $$
o{\`u} $\underline{\omega}_{\mathfrak{G}/\overline{M^1}}=
\overline{e}^*\Omega_{\mathfrak{G}/\overline{M^1}}$, en 
notant $\overline{e}:\overline{M^1}\rightarrow \mathfrak{G}$
la section  unit{\'e} de $\overline{f}$. En outre
$\underline{\omega}_{\mathfrak{G}/\overline{M}}$ co{\"\i}ncide avec
le faisceau $(\overline{f}_*\Omega_{\mathfrak{G}/\overline{M^1}})^{\mathfrak{G}}$
des $\mathfrak{G}$-invariants de 
$\overline{f}_*\Omega_{\mathfrak{G}/\overline{M^1}}$.

$\mathrm{(iv)}$ Le $\Z[\frac{1}{\N(\mathfrak{n})}]$-sch{\'e}ma 
$$M^{1*}=\Proj_{\Z[\frac{1}{\N(\mathfrak{n})}]}\left(\oplus_{k\geq 0}
\Gamma(\overline{M^1},
\underline{\omega}_{\mathfrak{G}/\overline{M^1}}^{kt})\right),$$
est ind{\'e}pendant du choix de  $\Sigma$. 
Le morphisme canonique  $ \pi :\overline{M^1}\rightarrow M^{1*} $
est surjectif et {\'e}quivariant pour l'action  du groupe fini 
$\mathfrak{o}_{D+}^\times/
(\mathfrak{o}_{D+}^\times\cap\mathfrak{o}_{\mathfrak{n}}^{\times 2})$.
Le quotient de $M^{1*}$ pour cette action est un sch{\'e}ma projectif,
normal, de type fini, not{\'e} $M^*$. La restriction {\`a} $M$
de la surjection canonique $ \pi :\overline{M}\rightarrow M^*$
 induit   un isomorphisme sur un ouvert dense de $M^*$, 
 not{\'e} encore $M$.

 $\mathrm{(v)}$
Le sch{\'e}ma $M^* \backslash M$ est fini et {\'e}tale sur 
$\Z[\frac{1}{\N(\mathfrak{n})}]$ et il est 
isomorphe {\`a} : $$ \coprod_{(R,\mathfrak{n})\mathrm{-composantes}/\sim}\textstyle
\Spec(\Z[\frac{1}{\N(\mathfrak{n})},\zeta_{\mathcal{C}}]^{H_{\mathcal{C}}}).$$
Les composantes connexes de $M^* \backslash M$ sont appel{\'e}es les
pointes de $M$. Cependant celles-ci ne sont des points ferm{\'e}s que pour les
 $(R,\mathfrak{n})$-composantes non-ramifi{\'e}es.

 $\mathrm{(vi)}$ La compl{\'e}tion formelle de 
$\overline{M}$ le long de l'image r{\'e}ciproque $\pi^{-1}(\mathcal{C})$ d'une   
$(R,\mathfrak{n})$-composante non-ramifi{\'e}e $\mathcal{C}$, est canoniquement
isomorphe {\`a} 
$$\left(S_{\Sigma^{\mathcal{C}}}^{\wedge}
/{}_{\mathfrak{o}_{\mathfrak{n}}^{\times}\times\mathfrak{o}_{D+}^\times}\right)
\times \textstyle \Spec(\Z[\frac{1}{\N(\mathfrak{n})}]).$$
\end{theo}

Pour le (iii), on voit localement  en passant aux vari{\'e}t{\'e}s de Tate au voisinage 
de chaque pointe que la fl{\`e}che de Kodaira-Spencer induit un isomorphisme.

\vspace{-.4cm}
\section[Compactifications toro{\"\i}dales des vari{\'e}t{\'e}s de Kuga-Sato.]{Compactifications toro{\"\i}dales des vari{\'e}t{\'e}s de Kuga-Sato}\label{kugasato}
 
Dans toute cette partie, on se limite au cas du groupe 
de niveau $\Gamma^1$, ($D=\Gm$) et donc $M^{\an}=M^{1,\an}$ de sorte qu'il
y a une VAHB  analytique universelle $\mathcal{A}^{\an}\rightarrow M^{\an}$.

\vspace{-.2cm}
\subsubsection*{Compactifications toro{\"\i}dales  des  vari{\'e}t{\'e}s analytiques
de Kuga-Sato.}

\medskip
La vari{\'e}t{\'e} analytique de Kuga-Sato $\mathcal{A}^{\an, s}$ 
est d{\'e}finie comme le produit fibr{\'e} $s$-fois de $\mathcal{A}^{\an }$
au-dessus de $M^{\an }$. Soit  $\overline{M^{\an }}$ une compactification toro{\"\i}dale
 de $M^{\an }$, comme  dans la partie pr{\'e}c{\'e}dente. On note de m{\^e}me $\mathfrak{G}^s$ le 
produit fibr{\'e} $s$ fois de $\mathfrak{G}$ par lui-m{\^e}me au-dessus de $\overline{M^{\an }}$.
On veut  compactifier  $\mathcal{A}^{\an, s}$  en une 
vari{\'e}t{\'e} analytique $\overline{\mathcal{A}^{\an,s }}$ projective lisse au-dessus de 
$\overline{M^{\an }}$ de mani{\`e}re {\`a} ce que le sch{\'e}ma semi-ab{\'e}lien $\mathfrak{G}^s$
op{\`e}re naturellement dessus, en prolongeant l'action par translation  de 
$\mathcal{A}^{\an, s}$ sur elle-m{\^e}me. Le caract{\`e}re naturel de ce prolongement sera
d{\'e}taill{\'e} dans l'{\'e}nonc{\'e} du th{\'e}or{\`e}me \ref{thmkugasato} plus bas.

Nous proc{\'e}derons en effectuant  des compactifications partielles de chaque  pointe {\`a} l'aide 
d'immersions toro{\"\i}dales, puis  en recollant ces derni{\`e}res 
on obtiendra la compactification cherch{\'e}e.

\medskip
Si on pose $(\mathcal{A}^{\an }_{\gamma,H})^s=
\gamma^{-1}\Gamma \gamma\cap B_{\R}\bs 
(W_{H}\times(F\otimes\C)^{s})/(\mathfrak{b}\oplus \mathfrak{a}^*)^{s} $, 
la description de la VAHB au voisinage de la pointe $\mathcal{C}=\gamma\infty$, 
faite dans  la partie  \ref{vahb}, donne :

\hspace{-1.35cm}\xymatrix@R=15pt@C=15pt{
 (\mathcal{A}^{\an }_{\gamma,H})^s
 \ar[d] &
X^*\bs
(W_{H}\times(F\otimes\C)^{s})/(\mathfrak{b}\oplus \mathfrak{a}^*)^{s}
 \ar[l]\ar[d]\ar@{^{(}->}[r] &
(\Gm\otimes X^*\times (\Gm\otimes \mathfrak{a}^*)^{s})/
\mathfrak{b}^{s} \ar[d],\\
  \gamma^{-1}\Gamma \gamma\cap B_{\R}\bs  W_{H}&
 X^*\bs W_{H}=D_{\gamma,H} 
\ar@{^{(}->}[r]\ar_{\mod{\mathfrak{o}_{\mathcal{C}}^{\times}}}[l] & 
\Gm\otimes X^*=:S_\mathcal{C} }

\noindent o{\`u} on rappelle que $\mathfrak{a}=\mathfrak{bc}$ et  $X=\mathfrak{cbb'}$.

\medskip

Le groupe $ \mathfrak{b}^s\rtimes \mathfrak{o}_{\mathcal{C}}^{\times} $ 
(produit semi-direct donn{\'e} par $(\beta_1,..,\beta_s;(u,\epsilon))
(\beta'_1,..,\beta'_s;(u',\epsilon'))\\
=(\beta_1+\beta'_1u^{-1}\epsilon^{-1},..,\beta_s+\beta'_su^{-1}\epsilon^{-1};
(uu',\epsilon\epsilon'))$)
agit {\`a} gauche sur $X^*_{+}\times (\mathfrak{a}^*)^s$, ainsi que sur 
 $(F\otimes\R)_+\times (F\otimes\R)^s$, par :
$$ (\beta_1,..,\beta_s;(u,\epsilon))\cdot(q;l_1,...,l_s)=
(u^2\epsilon q,ul_1+u^2\epsilon q\beta_1,...,ul_s+u^2\epsilon q\beta_s)$$
Notons que cette action est bien d{\'e}finie, 
car $X^*\mathfrak{b}\subset \mathfrak{a}^*$.

On aimerait ajouter {\`a} $\Gm\otimes X^*\times 
(\Gm\otimes \mathfrak{a}^*)^{s}$ une fronti{\`e}re 
analytique $\mathcal{F}$ au  dessus de la fronti{\`e}re 
analytique $\mathcal{E}$ de $S_\mathcal{C}$ et  sur laquelle
$\mathfrak{b}^s\rtimes\mathfrak{o}_{\mathcal{C}}^{\times}$ 
agit discontinument  et de mani{\`e}re compatible avec l'action 
de $\mathfrak{o}_{\mathcal{C}}^{\times}$ sur $\mathcal{E}$. 
Le quotient par $\mathfrak{b}^s\rtimes \mathfrak{o}_{\mathcal{C}}^{\times}$
 de l'adh{\'e}rence de $q(D_{\gamma,H})\times(\Gm\otimes \mathfrak{a}^*)^{s}$ 
dans $\Gm\otimes X^*\times (\Gm\otimes \mathfrak{a}^*)^{s}\cup \mathcal{F}$
serait alors la compactification partielle de la pointe $\mathcal{C}$ (voir
la partie \ref{toroidale}).

Le probl{\`e}me se traduit en le probl{\`e}me combinatoire suivant :

Soit un  {\'e}ventail complet $\Sigma$  de 
$X^*_{\R+}=\{0\}\cup \R_+^{*J_F}$, stable pour
l'action de $\mathfrak{o}_+^\times$ et qui 
contient un nombre fini d'{\'e}l{\'e}ments modulo cette action. Trouver 
un {\'e}ventail complet $\widetilde{\Sigma}$ de
$X^*_{\R+} \times (\mathfrak{a}^*_{\R})^s$ stable pour
l'action de $\mathfrak{b}^s\rtimes \mathfrak{o}_+^\times$,
 contenant un nombre fini d'{\'e}l{\'e}ments modulo cette action
et tel que la projection de chaque $\tau\in\widetilde{\Sigma}$ sur
$X^*_{\R+}$ soit un des $\sigma\in\Sigma$.

\medskip

Soit $\xi_0^*$ un {\'e}l{\'e}ment de $X^*_{+}$ de norme minimale.
Soit $\varepsilon_1,..,\varepsilon_{d-1}$
une base de $\mathfrak{o}^{\times2}$ et posons
$\Pi=\{\prod_{i\in I}\varepsilon_i \enspace |\enspace I\subset\{1,..,d-1\}\}$.
Alors l'ensemble $\Sigma$ des int{\'e}rieurs des 
$ \bigcap_{u\in U}\R_+\Conv_{\varepsilon\in\Pi}(u\varepsilon \xi_0^*)$,
avec $U$ d{\'e}crivant les sous-ensembles finis de $\mathfrak{o}^{\times2}$, 
est  un  {\'e}ventail complet de $X^*_{\R+}$, stable pour
l'action de $\mathfrak{o}^\times$ et  
contenant  un nombre fini d'{\'e}l{\'e}ments modulo cette action.

Soit $e^{(1)},..,e^{(d)}$ une base de $\mathfrak{b}$ et posons
$\Pi'=\{\sum_{i\in I}e^{(i)}\enspace |\enspace I\subset\{1,..,d\}\}$.
Consid{\'e}rons l'ensemble $\widetilde{\Sigma}''$ des c{\^o}nes ferm{\'e}s suivants :
$$\R_+\Conv_{\varepsilon\in\Pi;e_1,..,e_s\in\Pi'}
\Big((u^2\varepsilon^2\xi_0^*;(\beta_1+e_1)u\varepsilon^2,..,
(\beta_s+e_s)u\varepsilon^2)\Big),$$
 avec $\beta=(u;\beta_1,..,\beta_s)$ 
d{\'e}crivant $(\mathfrak{o}^{\times}/\{{\pm}1\})\times \mathfrak{b}^s$.

Il ne suffit pas de prendre les int{\'e}rieurs des intersections 
finies de tels c{\^o}nes pour obtenir l'{\'e}ventail cherch{\'e} $\widetilde{\Sigma}$, car
{\it l'intersection de deux c{\^o}nes de $\widetilde{\Sigma}''$
n'est pas forc{\'e}ment une face de chacun}.
Cependant, on n'en est pas loin, 
car une face donn{\'e}e d'un c{\^o}ne de $\widetilde{\Sigma}''$ ne rencontre qu'un 
nombre fini parmi les autres c{\^o}nes. De ce fait :

1) si on  d{\'e}coupe un des c{\^o}nes de $\widetilde{\Sigma}''$ (en prenant par exemple un 
point rationnel {\`a} l'int{\'e}rieur, bien que cela ne soit pas 
forc{\'e}ment n{\'e}cessaire) de fa{\c c}on que les autres c{\^o}nes intersectent 
les nouveaux c{\^o}nes ainsi obtenus en des faces de ces-derniers, et

2) si on d{\'e}coupe les autres c{\^o}nes en cons{\'e}quence, en translatant 
le d{\'e}coupage du 1) par le groupe
 $\mathfrak{b}^s\rtimes\mathfrak{o}_+^{\times}$,

alors on obtiendra  un nouvel ensemble de c{\^o}nes ferm{\'e}s $\widetilde{\Sigma}'$
qui sera plus fin que  $\widetilde{\Sigma}''$ et dans lequel 
l'intersection de deux c{\^o}nes sera une face de chacun.

  L'ensemble $\widetilde{\Sigma}$ des int{\'e}rieurs des intersections 
finies de  c{\^o}nes de $\widetilde{\Sigma}'$ sera alors
un {\'e}ventail complet  de $X^*_{\R+} \times (\mathfrak{a}^*_{\R})^s$
 stable pour l'action de $ \mathfrak{b}^s\rtimes\mathfrak{o}_+^\times$,
 contenant un nombre fini d'{\'e}l{\'e}ments modulo cette action
et tel que la projection de chaque $\tau\in \widetilde{\Sigma}$ sur
$(F\otimes\R)_+$ soit un des $\sigma\in\Sigma$. \hfill $\square$

Quitte {\`a} raffiner $\widetilde{\Sigma}$, on pourra le supposer lisse.

Par  la m{\^e}me m{\'e}thode que dans la partie \ref{toroidale} on obtient 
alors une  compactification lisse de la forme voulue  de  
$\mathcal{A}^{\an,s}$. L'{\'e}nonc{\'e} pr{\'e}cis
sera donn{\'e} plus tard dans cette partie, dans le cas arithm{\'e}tique
(le cas analytique en d{\'e}coule par les arguments habituels).

\vspace{-.4cm} 
\subsubsection*{Compactifications arithm{\'e}tiques des vari{\'e}t{\'e}s de Kuga-Sato.} 

On rappelle que dans cette partie, on se limite au cas 
du sch{\'e}ma de modules fin $M=M^1$, de sorte qu'il
existe une VAHB universelle $\mathcal{A}\rightarrow M$ par le Th{\'e}or{\`e}me 4.1.
Soit $\Sigma$ un {\'e}ventail  $\Gamma$-admissible
et $M\hookrightarrow \overline{M}_{\Sigma}= \overline{M}\hookleftarrow D$ 
la compactification toro{\"\i}dale 
associ{\'e}e, avec $D$ diviseur {\`a} croisements normaux.

Pour chaque entier $s\geq 1$ on d{\'e}finit la vari{\'e}t{\'e} de Kuga-Sato
$\mathcal{A}^s=\mathcal{A}\underset{M}{\times}...\underset{M}{\times} \mathcal{A}$ ($s$-fois), qui 
est munie d'un morphisme projectif lisse $f_s:\mathcal{A}^s\rightarrow M$. 

Le but est de construire (en s'inspirant de \cite{FaCh}) des 
 compactifications toro{\"\i}dales 
$\mathcal{A}^s\hookrightarrow \overline{\mathcal{A}^s}=
\mathcal{A}^s_{\widetilde{\Sigma}}\hookleftarrow E$,
avec $E$ diviseur {\`a} croisements normaux relatif, au-dessus des 
compactifications toro{\"\i}dales  de $M$. En d'autres termes on veut obtenir 
un diagramme : 
$$\xymatrix@R=15pt{ \mathcal{A}^s \ar@{^{(}->}[r]\ar^{f_s}[d] & 
\overline{\mathcal{A}^s}\ar^{\overline{f_s}}[d] & E \ar@{_{(}->}[l]\ar[d] \ar[d]\\
M \ar@{^{(}->}[r]\ar[dr] &  \overline{M}_{\Sigma}\ar[d] &  D \ar@{_{(}->}[l]\ar[dl] \\
& \Spec(\Z[\frac{1}{\N(\mathfrak{n})}])&}$$
 avec $\overline{f_s}$ semi-stable et projectif.

L'importance de l'existence d'une vari{\'e}t{\'e} $\overline{\mathcal{A}^s}\rightarrow \overline{M}$
pour chaque valeur de l'entier $s$ appara{\^\i}tra clairement dans la section sur la th{\'e}orie
de Hodge. Pour d{\'e}montrer le th{\'e}or{\`e}me de d{\'e}g{\'e}n{\'e}rescence de la suite spectrale BGG duale
vers de Rham, on doit en effet recourir, suivant \cite{FaCh}, VI.5.5,  au th{\'e}or{\`e}me de Deligne de d{\'e}g{\'e}n{\'e}rescence de la
suite spectrale de Hodge vers de Rham pour $\overline{f_s}:\overline{\mathcal{A}^s}\rightarrow \overline{M}$.

L'{\'e}ventail consid{\'e}r{\'e} dans la partie pr{\'e}c{\'e}dente 
pour la compactification analytique ne peut pas {\^e}tre r{\'e}utilis{\'e} ici,
car la m{\'e}thode de \cite{FaCh} utilise des {\'e}ventails
munis d'une fonction de polarisation. On utilise 
la d{\'e}composition de Voronoi-Delaunay, qui est naturellement
munie d'une fonction de polarisation (voir aussi \cite{kunn}).

\subsubsection*{Donn{\'e}es de d{\'e}g{\'e}n{\'e}rescence.}

\begin{defin} Soient $\mathfrak{a}$ et $\mathfrak{b}$ deux id{\'e}aux 
de $\mathfrak{o}$ tels que $\mathfrak{ab}^{-1}=\mathfrak{c}$.
Des   donn{\'e}es de d{\'e}g{\'e}n{\'e}rescence, pour la vari{\'e}t{\'e} de Kuga-Sato consistent en :

\noindent$\relbar$ (polarisation) des morphismes $\mathfrak{o}$-lin{\'e}aires 
$\phi_i:\mathfrak{b}\rightarrow \mathfrak{a}$, $1\leq i \leq s$.

\smallskip
\noindent$\relbar$ (fl{\`e}che tautologique) une forme bilin{\'e}aire  
$b:\mathfrak{b}\times \mathfrak{a}\rightarrow \Z$ telle
que pour tout $m\in \mathfrak{o}$, $\alpha \in \mathfrak{a}$ et 
$\beta\in\mathfrak{b}$ on ait $b(m\beta,\alpha )=b(\beta,m\alpha )$ 
et  telle pour tout $1\leq i \leq s$ l'application 
 $b(\cdot,\phi_i(\cdot))$ soit une forme bilin{\'e}aire d{\'e}finie positive sur 
$\mathfrak{b}$.
\end{defin}

Soit $\mathcal{C}$ une $(R,\mathfrak{n})$-pointe, donn{\'e}e par un r{\'e}seau $L$ de $F^2$, 
une  suite exacte de $\mathfrak{o}$-modules 
$0\rightarrow \mathfrak{a}^* \rightarrow L \rightarrow \mathfrak{b} \rightarrow 0$ 
et une application $\mathfrak{o}$-lin{\'e}aire injective 
$\alpha:\mathfrak{n}^{-1}\mathfrak{d}^{-1}/\mathfrak{d}^{-1}\hookrightarrow
\mathfrak{n}^{-1}L/L$. On associe {\`a} $\mathcal{C}$ l'id{\'e}al 
$X=\mathfrak{ab}'\supset \mathfrak{ab}$ (voir \cite{dimdg} section 3).
A chaque $(\xi^*,(\mu_i)_{1\leq i\leq s})\in X_{+}^{*}\times  \mathfrak{c}_{+}^s$,  
 on peut  associer des donn{\'e}es de d{\'e}g{\'e}n{\'e}rescence
$\phi_i=\phi_{\mu_i}$ 
et $b=b_{\xi^*}$, d{\'e}finies par : pour tout $\alpha\in \mathfrak{a}$,
$\beta\in \mathfrak{b}$  et $1\leq i \leq s$ \enspace  
$\phi_i(\beta)=\mu_i\beta$ et  $b(\beta,\alpha)=\Tr_{F\!/\!\Q}(\xi^*\alpha\beta)$.

\smallskip

On pose $C_+=X_{+}^{*}$ et $\widetilde{C}_+=C_+\times (\mathfrak{a}^*)^s$.
Le groupe $ \mathfrak{b}^s\rtimes \mathfrak{o}_{\mathcal{C}}^{\times} $ 
agit {\`a} gauche sur $\widetilde{C}_+$  (de m{\^e}me que dans (\ref{vahb-univ})
le groupe $(\mathfrak{o}\oplus\mathfrak{c}^*)\rtimes \Gamma$ agit sur 
$\mathfrak{H}_F\times (F\otimes \C)$) par
$$\begin{cases} 
(u,\epsilon)(q;l_1,..,l_s)=(u^2\epsilon q;ul_1,..,ul_s) \\
(\beta_1,..,\beta_s)(q;l_1,..,l_s)=(q;l_1+\beta_1q,..,l_s+\beta_sq)
\end{cases}.$$

\vspace{-4mm}
\subsubsection*{Fonctions de polarisation.} Le but est de construire :

\noindent$\bullet$ Un {\'e}ventail $\Sigma^{\mathcal{C}}$ de $C_{\R+}$  qui est 
$\mathfrak{o}^{\times}$-admissible. 

\noindent$\bullet$ Un {\'e}ventail $\widetilde{\Sigma}^{\mathcal{C}}$ 
de $\widetilde{C}_{\R+}$  qui est 
$ \mathfrak{b}^s\rtimes\mathfrak{o}^{\times}$-admissible 
et tel que pour tout $\tau\in \widetilde{\Sigma}^{\mathcal{C}}$, 
il existe $\sigma \in \Sigma^{\mathcal{C}}$ 
tel que $\pr_1(\tau)\subset\sigma$. Si de plus  cette 
inclusion est une {\'e}galit{\'e} l'{\'e}ventail sera dit  {\it {\'e}quidimensionnel}. 

\noindent$\bullet$ Une {\it fonction de support} sur $\widetilde{\Sigma}^{\mathcal{C}}$,
i.e. une fonction $\varphi: \widetilde{C}_{\Q+}\rightarrow \Q$ qui est 
continue, convexe,  enti{\`e}re sur $X^*\times (\mathfrak{a}^*)^s$
lin{\'e}aire sur chaque $\tau\in \widetilde{\Sigma}^{\mathcal{C}}$ 
(donc lin{\'e}aire par morceaux),
$ \mathfrak{b}^s\rtimes\mathfrak{o}^{\times}$-invariante 
et telle que pour tout $\lambda\geq 0$ 
$\varphi(\lambda\cdot)=\lambda\varphi(\cdot)$.

Si de plus $\varphi$ est strictement convexe au-dessus de chaque 
$\sigma\in \Sigma^{\mathcal{C}}$ (i.e. pour tout 
$\tau\in \widetilde{\Sigma}^{\mathcal{C}}$, il existe 
$\sigma\in \Sigma^{\mathcal{C}}$, $n\in \n$ et 
$\widetilde{l}^*\in X\times \mathfrak{a}^s$ tels que
$$\tau=\{\widetilde{l}=(q,l)\in  \widetilde{C}_+| q\in \sigma, \enspace
n\varphi(\widetilde{l})=\langle\widetilde{l}^*,\widetilde{l}\rangle \},$$ 
alors $ \varphi$ est appel{\'e}e une {\it fonction de polarisation}. 

\vspace{-2mm}
\subsubsection*{D{\'e}composition de Voronoi-Delaunay.} 

Fixons  une $(R,\mathfrak{n})$-pointe  $\mathcal{C}$,
ainsi que des  $\mu_i\in \mathfrak{c}_{+}$, $1\leq i \leq s$.
On a ainsi des polarisations 
$\phi_i=\phi_{\mu_i}$, $1\leq i \leq s$.

Pour tout choix de  $\beta=(\beta_i)_{1\leq i \leq s}\in \mathfrak{b}^s$,
 on d{\'e}finit une  fonction : 
$$\chi_{\beta}:\widetilde{C}_{\Q+} \rightarrow   \Q\text{ , }
(q,l_1,..,l_s)\mapsto\sum_{1\leq i \leq s}
b_q(\beta_i,\phi_i(\beta_i))+2l_i(\phi_i(\beta_i)).$$

L'application $\chi_{\beta}$ est la compos{\'e}e de 
$$(q,l_1,..,l_s)\mapsto \sum_{1\leq i \leq s}
q\mu_i\beta_i^2+2l_i\mu_i\beta_i$$ avec l'application 
trace  $\Tr_{F\!/\!\Q}:F\rightarrow \Q$.

On pose pour $\widetilde{l}=(q,l=(l_1,..,l_s))\in  \widetilde{C}_{\Q+}$
$$\varphi(\widetilde{l})=\underset{\beta\in \mathfrak{b}^s}{\min}\chi_\beta(\widetilde{l}),$$

L'application  $\varphi$ est $1$-{\it tordue}, au sens de \cite{FaCh}, 
car pour tout  $\beta\in \mathfrak{b}^s$ on a 
$$\hspace{-5mm}
\varphi(\beta\cdot(q,l))=
\underset{\beta'\in \mathfrak{b}^s}{\min}\chi_{\beta'}(q,l+q\beta)=
\underset{\beta'\in \mathfrak{b}^s}{\min}\chi_{\beta+\beta'}(q,l)-\chi_\beta(q,l)=
\varphi(\widetilde{l})-\chi_\beta(\widetilde{l})$$

Pour $\sigma\in \Sigma$ fix{\'e} et $\mathfrak{B}\subset \mathfrak{b}^s$ un sous-ensemble fini on  d{\'e}finit :
$$\tau_{\sigma,\mathfrak{B}}=\{  \widetilde{l}=(q,l)\in  \widetilde{C}_+ | 
 q\in \sigma , \enspace \forall \beta\in \mathfrak{B} \enspace 
\chi_\beta(\widetilde{l})=\varphi(\widetilde{l}) \}.$$

\begin{prop} L'{\'e}ventail $\widetilde{\Sigma}=\{\tau_{\sigma,\mathfrak{B}}\}$ est 
un  {\'e}ventail complet 
$\mathfrak{b}^s\rtimes \mathfrak{o}^{\times}-$ admissible, 
{\'e}quidimensionnel de $\widetilde{C}_+$ et $\varphi$ est une 
fonction de polarisation $1$-tordue. 
Il existe une subdivision lisse de $\widetilde{\Sigma}$, 
muni d'une polarisation $k$-tordue, pour un certain $k\geq 1$. 
\end{prop}

On se propose de calculer l'action de 
$\mathfrak{b}^s\rtimes \mathfrak{o}^{\times}$ sur
l'{\'e}ventail $\widetilde{\Sigma}$. 

Pour $\beta\in \mathfrak{b}^s$ on a $\chi_\beta(q,l)=\varphi(q,l)$, si et
seulement si pour tout
$,e\in \mathfrak{b}^s$ on a $\chi_{\beta+e}(q,l)-\chi_\beta(q,l) 
=\Tr_{F\!/\!\Q}(e(2l+q(2\beta+e))\mu )\geq 0$. On en d{\'e}duit 
$$\tau_{\sigma,\mathfrak{B}}=\{  (q,l)\in  \widetilde{C}_+ |
q\in \sigma , \enspace \forall \beta\in \mathfrak{B}, \enspace 
\forall e \in \mathfrak{b}^s\enspace \Tr_{F\!/\!\Q}(e(2l+q(2\beta+e))\mu)\geq 0\}.$$

Pour tout  $u\in \mathfrak{o}^{\times}$ on a
$u\cdot \tau_{\sigma,\mathfrak{B}}=$

$=\{  (u^2q,ul)\in  \widetilde{C}_+ |
q\in \sigma, \forall \beta\in \mathfrak{B}, 
\forall e \in \mathfrak{b}^s\enspace \Tr_{F\!/\!\Q}(e(2l+q(2\beta+e))\mu )\geq 0\}=$

$\hspace{-3mm}\{  (u^2q,ul)\in  \widetilde{C}_+ |
u^2q\in u^2\sigma, \forall \beta,  \forall e 
\Tr_{F\!/\!\Q}(u^{-1}e(2ul+u^2q(2u^{-1}\beta+u^{-1}e))\mu )\geq 0\}$

$ =\tau_{u^2\sigma,u^{-1}\mathfrak{B}}$. 

\smallskip
Si $y\in \mathfrak{b}^s$ on a $y\cdot \tau_{\sigma,\mathfrak{B}}=$

$=\{  (q,l+qy)\in  \widetilde{C}_+ |
q\in \sigma ,  \forall \beta\in \mathfrak{B},  
\forall e \in \mathfrak{b}^s\enspace \Tr_{F\!/\!\Q}(e(2l+q(2\beta+e))\mu )\geq 0\}$

$=\{  (q,l-2qy)\in  \widetilde{C}_+ | q\in \sigma , \forall \beta, \forall e 
\enspace \Tr_{F\!/\!\Q}(e(2(l+qy)+q(2(\beta-y)+e))\mu )\geq 0\}$

$ =\tau_{\sigma,\mathfrak{B}-y}$.

\medskip
Le diagramme suivant d{\'e}crit l'action de 
$\mathfrak{b}^s\rtimes \mathfrak{o}^{\times}$ sur l'{\'e}ventail $\widetilde{\Sigma}$.
$$\xymatrix@R=15pt{ \tau_{\sigma,\mathfrak{B}}\ar^{\cdot u}[rr]\ar^{\cdot y}[d] &&
\tau_{u^2\sigma,u^{-1}\mathfrak{B}} \ar^{\cdot y u^{-1}=u\cdot y}[d]\\
\tau_{\sigma,\mathfrak{B}-y}\ar^{\cdot u}[rr] &&
\tau_{u^2\sigma,u^{-1}\mathfrak{B}-u^{-1}y}}$$

\vspace{-2mm}
\subsubsection*{Mod{\`e}les relativement complets faibles.}

On introduit la notion de {\it mod{\`e}les relativement 
complets faibles polaris{\'e}s} dans le  cas totalement d{\'e}g{\'e}n{\'e}r{\'e} 
qui nous int{\'e}resse (voir \cite{FaCh} VI.1.7, ainsi que la partie 2 de \cite{mum}).

Soit $R$ un anneau excellent, int{\'e}gralement clos, noeth{\'e}rien, complet 
pour la  topologie $I$-adique, pour un   id{\'e}al radiciel $I=\sqrt{I}$.
Soit $K$ le corps des fractions de $R$. 
Soit $S=\Spec(R)$, $\eta$ son point g{\'e}n{\'e}rique  et 
 $S_0=\Spec(R/I)$ le sous-sch{\'e}ma ferm{\'e} d{\'e}fini par $I$.

Consid{\'e}rons le tore d{\'e}ploy{\'e} $\widetilde{G}=(\mathbb{G}_m\otimes \mathfrak{a}^*)^s
\times S=\Spec(R[\mathfrak{X}^\alpha; \alpha\in \mathfrak{a}^s])$  sur $S$. 
Un ensemble de p{\'e}riodes $\mathfrak{b}^s
\subset  \widetilde{G}(K)$ {\'e}quivaut {\`a} la donn{\'e}e d'une application bilin{\'e}aire
non-d{\'e}g{\'e}n{\'e}r{\'e}e $\mathfrak{b}^s\times \mathfrak{a}^s\rightarrow K^\times, \enspace
(\beta,\alpha)\mapsto \mathfrak{X}^{\alpha}(\beta)$.  Une 
polarisation $\phi$ sur l'ensemble des p{\'e}riodes $\mathfrak{b}^s$
est la donn{\'e}e d'un homomorphisme $\mathfrak{o}$-lin{\'e}aire  $\phi:\mathfrak{b}^s\rightarrow \mathfrak{a}^s$,
tel que :

$\mathrm{(i)}$  $\mathfrak{X}^{\phi(\beta)}(\beta')=\mathfrak{X}^{\phi(\beta')}(\beta)$, pour tout  
$\beta,\beta'\in \mathfrak{b}^s$,

$\mathrm{(ii)}$ $\mathfrak{X}^{\phi(\beta)}(\beta)\in I$, pour tout  $\beta\in \mathfrak{b}^s\bs \{0\}$.

\begin{defin}
Un {\it mod{\`e}le relativement complet faible polaris{\'e}} de $\widetilde{G}$, 
par rapport {\`a} $(\mathfrak{b}^s,\phi)$,  est la donn{\'e}e des {\'e}l{\'e}ments suivants :

\noindent$\mathrm{(a)}$  Un sch{\'e}ma int{\`e}gre $\widetilde{P}$, 
localement de type fini sur $R$, dont la fibre g{\'e}n{\'e}rique est isomorphe 
{\`a} $\widetilde{G}_\eta$.

\noindent$\mathrm{(b)}$  Un faisceau inversible $\widetilde{\mathcal{L}}$ sur $\widetilde{P}$.

\noindent$\mathrm{(c)}$  Une action du tore $\widetilde{G}$ sur 
$(\widetilde{P},\widetilde{\mathcal{L}})$, {\'e}tendant l'action par translation sur la fibre 
g{\'e}n{\'e}rique et son faisceau structural.
On note  cette action $S_g:\widetilde{P}\rightarrow\widetilde{P}$, $S_g^*:
\widetilde{\mathcal{L}}\rightarrow\widetilde{\mathcal{L}}$,
pour tout point fonctoriel $g$ de $\widetilde{G}$.

\noindent$\mathrm{(d)}$  Une action de $\mathfrak{b}^s$ 
sur ($\widetilde{P},\widetilde{\mathcal{L}}$),
not{\'e}e $T_{\beta}:\widetilde{P}\rightarrow\widetilde{P}$ et 
$T_{\beta}^*:\widetilde{\mathcal{L}}\rightarrow\widetilde{\mathcal{L}}$, {\'e}tendant l'action 
de $\mathfrak{b}^s$ sur $\widetilde{G}_\eta$ par translation (via $\mathfrak{b}^s\subset \widetilde{G}(K)$).

\smallskip

 satisfaisant aux conditions suivantes :
\smallskip

\noindent$\mathrm{(i)}$  Il existe un ouvert $\widetilde{G}$-invariant 
$U\subset \widetilde{P}$ de type fini sur $S$ et tel que 
$\widetilde{P}=\cup_{\beta\in\mathfrak{b}^s}T_{\beta}(U)$.

\noindent$\mathrm{(ii)}$  $\widetilde{\mathcal{L}}$ est ample sur  $\widetilde{P}$, au sens 
que les compl{\'e}ments des lieux des z{\'e}ros des sections globales 
 $\Gamma(\widetilde{P},\widetilde{\mathcal{L}}^{\otimes n})$, $n\geq 1$, forment
une base de la topologie de Zariski de $\widetilde{P}$.

\noindent$\mathrm{(iii)}$  Pour toute valuation $v$ sur $R(\widetilde{G})$ (le corps des
fonctions rationnelles sur $\widetilde{G}$), qui est positive sur $R$, on a :

\noindent $v$ a du centre sur $\widetilde{P}$  $\iff$ pour tout 
$\alpha\in \mathfrak{a}^s$, il existe  $\beta\in\mathfrak{b}^s$ avec 
$v(\mathfrak{X}^{\alpha}(\beta)\mathfrak{X}^{\alpha})\geq 0$.

\end{defin}

L'int{\'e}r{\^e}t des  mod{\`e}les relativement complets faibles polaris{\'e}s
$(\widetilde{P},\widetilde{\mathcal{L}})$ est qu'en suivant les fl{\`e}ches du diagramme :
\xymatrix@R=15pt{  \widetilde{P} & & 
\widetilde{\mathfrak{P}}
\ar_{\text{compl{\'e}tion}}[ll]
\ar_{\text{quotient formel par } \mathfrak{b}^s}[d]\\
P  && \mathfrak{P}\ar^{\text{alg{\'e}brisation}}[ll]  }

\noindent l'on peut construire 
``le quotient'' $(P,\mathcal{L})$ de $(\widetilde{P},\widetilde{\mathcal{L}})$ par 
le groupe des p{\'e}riodes $\mathfrak{b}^s$. Nous allons utiliser cette construction dans
le th{\'e}or{\`e}me suivant.

\subsubsection*{{\'E}nonc{\'e} du th{\'e}or{\`e}me.} 

Soient $\mu_i\in\mathfrak{c}_+$, $1\leq i\leq s$.
Soit $\mathcal{A}^s=\mathcal{A}\times_{M}\ldots\times_{M} \mathcal{A}$.
Soit $\Sigma=(\Sigma^{\mathcal{C}})_{\mathcal{C}}$
 un {\'e}ventail complet et lisse de  $C_{\R+}$  
qui est  $\mathfrak{o}^{\times}$-admissible et  
soit  $\widetilde{\Sigma}=(\widetilde{\Sigma}^{\mathcal{C}})_{\mathcal{C}}$ {\'e}ventail 
complet et lisse de  $\widetilde{C}_{\R+}$ qui est 
$\mathfrak{b}^s\rtimes \mathfrak{o}^{\times}$-admissible,
{\'e}quidimensionnel au-dessus de $\Sigma$ et muni d'une  
fonction de polarisation $k$-tordue $\varphi$.

\begin{theo}\label{thmkugasato}
Il existe un $\Z[\frac{1}{\N(\mathfrak{n})}]$-sch{\'e}ma $\overline{\mathcal{A}^s}=\mathcal{A}^s_{\widetilde{\Sigma}}$
propre (et m{\^e}me projectif) sur $\overline{M}= M_{\Sigma}$,
muni d'un faisceau inversible ample $\mathcal{L}$ tel que:

$(\mathrm{i})$ 
$\overline{\mathcal{A}^s}|_{M}=\mathcal{A}^s$ est la vari{\'e}t{\'e} de Kuga-Sato
universelle au-dessus de $M$ et  $\mathcal{L}|_{\mathcal{A}^s}$ 
s'identifie avec la puissance tensorielle $k$-i{\`e}me du faisceau 
inversible ample $\otimes_i \mathrm{pr}_i^* \mathcal{L}_{\mu_i}$, 
o{\`u} pour $1\leq i\leq s$,  
$\mathrm{pr}_i:\mathcal{A}^s\rightarrow \mathcal{A}$ 
d{\'e}signe la $i$-i{\`e}me projection et $\mathcal{L}_{\mu_i}$ d{\'e}signe  le 
faisceau ample  inversible canonique sur $\mathcal{A}$, obtenu par 
pull-back du faisceau  de Poincar{\'e} par le morphisme  
$(\mathrm{id}_{\mathcal{A}},
\lambda\circ(\mathrm{id}_{\mathcal{A}}\otimes\mu_i))$.

$(\mathrm{ii})$ $\overline{\mathcal{A}^s}$ poss{\`e}de une stratification 
naturelle param{\'e}tr{\'e}e par  $\widetilde{\Sigma}/(\mathfrak{b}^s\rtimes
 \mathfrak{o}^{\times})$.

$(\mathrm{iii})$ Le sch{\'e}ma $\overline{\mathcal{A}^s}$ est lisse sur 
$\Z[\frac{1}{\Delta}]$ et 
$\overline{\mathcal{A}^s}\bs \mathcal{A}^s$ est un diviseur {\`a} 
croisements normaux relatif sur $\overline{M}$. Le morphisme 
$\overline{f_s}:\overline{\mathcal{A}^{s}}\rightarrow\overline{M}$ 
est semi-stable.

\medskip
Supposons que  pour tout $\sigma\in \Sigma$, il existe 
$\tau\in \widetilde{\Sigma}$ tel que $\sigma\times\{0\}=\tau$. Alors :

$(\mathrm{iv})$ Le sch{\'e}ma semi-ab{\'e}lien $\mathfrak{G}^s$ est contenu 
comme ouvert dense dans $\overline{\mathcal{A}^s}$  et 
la restriction de $\mathcal{L}$ {\`a} $\mathfrak{G}^s$ co{\"\i}ncide, comme dans
le $(\mathrm{i})$   avec  la puissance tensorielle $k$-i{\`e}me du faisceau 
inversible ample canonique $\otimes_i \mathrm{pr}_i^* \mathcal{L}_{\mu_i}$.
 De plus $\mathfrak{G}^{s}\rightarrow
 \overline{M}$  agit sur $\overline{\mathcal{A}^{s}}$ en prolongeant 
l'action de $\mathcal{A}^s$ sur lui-m{\^e}me par translation. 

$(\mathrm{v})$ Le faisceau $\Omega^1_{\overline{\mathcal{A}^s}/\overline{M}}
(\dlog\infty)$ est isomorphe {\`a} 
$\overline{f_s}^*({\underline{\omega}_{\mathfrak{G}/\overline{M}}}^{\oplus s})$.

$(\mathrm{vi})$ Pour tout couple d'entiers $a,b\geq 0$, on a des 
isomorphismes canoniques
 $$R^a \overline{f_s}_*\left(\bigwedge^b\Omega^1_{\overline{\mathcal{A}^s}/
\overline{M}}(\dlog\infty)\right)
\cong \bigwedge^a (\mathfrak{cd}^{-1}\otimes
\underline{\omega}_{\mathfrak{G}/\overline{M}}^\vee)^{\oplus s}
\otimes \bigwedge^b
{\underline{\omega}_{\mathfrak{G}/\overline{M}}}^{\oplus s}.$$
\end{theo}

Dans le reste de l'article on abr{\'e}gera 
$\Omega^\bullet_{\overline{\mathcal{A}^s}/\overline{M}} 
(\dlog\infty_{\overline{\mathcal{A}^s}/\overline{M}})$ 
en  $\Omega^\bullet_{\overline{\mathcal{A}^s}/\overline{M}}(\dlog\infty)$.

\begin{rque}
La canonicit{\'e} des isomorphismes de $(\mathrm{vi})$ montre 
en particulier que les  faisceaux
$R^a \overline{f_s}_*\bigwedge^b
\Omega^1_{\overline{\mathcal{A}^s}/\overline{M}}(\dlog\infty)$ 
sont :

{\rm 1)}
ind{\'e}pendants du choix de la compactification toro{\"\i}dale 
de $\mathcal{A}^s$ choisie,

{\rm 2)}
munis d'une action naturelle de $\mathfrak{G}^s$ et de 
$\mathfrak{o}$,

{\rm 3)} localement
 libres sur $\mathcal{O}_{\overline{M}}\otimes \mathfrak{o}$.
\end{rque}

\subsubsection*{D{\'e}monstration  du th{\'e}or{\`e}me.} 

On construit $P=\overline{\mathcal{A}^s}$ en suivant les {\'e}tapes 
de \cite{FaCh} VI.1 :
Pour chaque $(R,\mathfrak{n})$-pointe $\mathcal{C}$ la projection de
$\widetilde{\Sigma}^{\mathcal{C}}$ sur $\Sigma^{\mathcal{C}}$ nous donne
un morphisme d'immersions toro{\"\i}dales (voir \cite{KKMS}) :
\xymatrix@R=15pt{ \widetilde{S}_{\mathcal{C}}\ar@{^{(}->}[r]\ar[d] & 
\widetilde{S}_{\widetilde{\Sigma}^{\mathcal{C}}}\ar[d]\\
S_{\mathcal{C}}\ar@{^{(}->}[r] & S_{\Sigma^{\mathcal{C}}}}

\noindent o{\`u} $S_{\mathcal{C}}$ (resp. $\widetilde{S}_{\mathcal{C}}$)
d{\'e}signe le tore d{\'e}ploy{\'e}  de groupe des caract{\`e}res  $X$ 
(resp. $X\times \mathfrak{a}^s$). Le morphisme 
$\widetilde{S}_{\widetilde{\Sigma}^{\mathcal{C}}}\rightarrow 
S_{\Sigma^{\mathcal{C}}}$  est {\'e}quivariant pour l'action des
tores $\widetilde{S}_{\mathcal{C}}\rightarrow S_{\mathcal{C}}$
et pour l'action des groupes $\mathfrak{b}^s\rtimes \mathfrak{o}^\times
\rightarrow \mathfrak{o}^\times$.

\medskip
Il est  crucial de noter alors que :

$\relbar$ la fonction de polarisation 
$\varphi:\widetilde{C}_+\rightarrow \Z$ induit un faisceau inversible 
relativement ample $\mathcal{L}$ sur $S_{\widetilde{\Sigma}^{\mathcal{C}}}$
(voir \cite{KKMS}).

$\relbar$ le fait que $\varphi$ est $k$-tordue nous donne, pour 
tout $\beta\in \mathfrak{b}^s$ et tout point $g$ de $\widetilde{G}^s=(\Gm\otimes \mathfrak{a}^*)^s$
 la relation 
$$S^*_g T^*_\beta=\mathfrak{X}^{\phi(\beta)}(g)^{2k}T^*_\beta S^*_g, $$
 qui est similaire {\`a} celle impos{\'e}e en plus dans la d{\'e}finition des
mod{\`e}les relativement complets (voir \cite{mum}2.1(iv)).

$\relbar$ pour tout $\sigma\in \Sigma^{\mathcal{C}}$ le pull-back
de $(\widetilde{S}_{\widetilde{\Sigma}^{\mathcal{C}}},\mathcal{L})$ par 
le morphisme $\overline{S}_{\mathcal{C}, \sigma}\rightarrow 
S_{\Sigma^{\mathcal{C}}}$ est un mod{\`e}le relativement complet faible 
polaris{\'e} du tore $\widetilde{G}^s\times \overline{S}_{\mathcal{C}, \sigma}$,
relativement {\`a} $(\mathfrak{b}^s, (\phi_{\mu_i})_{1\leq i\leq s})$.
\medskip

Ainsi, par le r{\'e}sultat principal sur ces mod{\`e}les \cite{FaCh} VI.1.10, on obtient un sch{\'e}ma 
propre $P_\sigma$ sur 
$\overline{S}_{\mathcal{C}, \sigma}[\frac{1}{\N(\mathfrak{n})},\zeta_{\mathcal{C}}]$,
prolongeant le pull-back $\mathcal{A}^s_\sigma$ de la vari{\'e}t{\'e} de 
Kuga-Sato universelle  {\`a}
$\overline{S}{}_{\mathcal{C}, \sigma}^0[\frac{1}{\N(\mathfrak{n})},\zeta_{\mathcal{C}}]$, 
et un faisceau inversible ample  $\mathcal{L}_\sigma$ sur $P_\sigma$,
prolongeant le faisceau inversible ample canonique 
$\otimes_i \mathrm{pr}_i^* \mathcal{L}_{\mu_i}$
de $\mathcal{A}^s_\sigma$.

\medskip

Par compatibilit{\'e} des immersions toriques, comme dans 
\cite{FaCh} IV.3 p.104, on obtient  un sch{\'e}ma  propre
$$g_{\mathcal{C}}: P_{\Sigma^{\mathcal{C}}}\rightarrow \textstyle
S_{\Sigma^{\mathcal{C}}}^\wedge[\frac{1}{\N(\mathfrak{n})},\zeta_{\mathcal{C}}],$$

\noindent appel{\'e} le ``{\it bon mod{\`e}le formel compact}'' 
en la pointe   $\mathcal{C}$,   et  un faisceau inversible ample 
$\mathcal{L}_{\Sigma^{\mathcal{C}}}$ sur $P_{\Sigma^{\mathcal{C}}}$.
 Le couple  $(P_{\Sigma^{\mathcal{C}}},\mathcal{L}_{\Sigma^{\mathcal{C}}})$ 
  descend {\`a}  $S_{\Sigma^{\mathcal{C}}}^\wedge/\mathfrak{o}^\times_{\mathcal{C}}
\times\Spec(\Z[\frac{1}{\N(\mathfrak{n})},\zeta_{\mathcal{C}}]^{H_{\mathcal{C}}})$.

On peut alors  alg{\'e}briser et  recoller ces sch{\'e}mas, ainsi que les faisceaux 
inversibles amples, pour obtenir un morphisme $\overline{f}:P\rightarrow 
\overline{M}$ et un  faisceau inversible ample $\mathcal{L}$ sur $P$, de sorte que :

1)  $\overline{f}$ est projectif sur $\overline{M}$,

2) on a canoniquement $P\vert_M\cong \mathcal{A}^s$,

3)  le sch{\'e}ma semi-ab{\'e}lien $\mathfrak{G}^s$  
 agit sur $P$ en  prolongeant l'action  par 
translation de $\mathcal{A}^s$ sur lui-m{\^e}me au-dessus de $M$,

4)  si pour tout $\mathcal{C}$, $\widetilde{\Sigma}^{\mathcal{C}}$ 
contient les $\sigma\times \{0\}$ 
pour tout $\sigma\in \Sigma^{\mathcal{C}}$, le sch{\'e}ma $\mathfrak{G}^s$ est muni d'une
immersion ouverte $j:\mathfrak{G}^s\hookrightarrow P$ d'image dense dans $P$,
et le faisceau  
$j^*\mathcal{L}$ co{\"\i}ncide avec le faisceau ample canonique sur $\mathfrak{G}^s$
  associ{\'e} {\`a}  $(\mu_1,\ldots,\mu_s)$,
via la $\mathfrak{c}$-polarisation canonique 
sur le sch{\'e}ma semi-ab{\'e}lien $\mathfrak{G}^s$ prolongeant celle de $\mathcal{A}^s$.

5) Pour tout c{\^o}ne  $\tau\in\widetilde{\Sigma}^{\mathcal{C}}$ 
et  $\sigma\in\Sigma^{\mathcal{C}}$, avec $\pr_1(\tau)=\sigma$, la compl{\'e}tion formelle de 
$P\rightarrow \overline{M}$ le long de la $\tau$-strate (au-dessus de la $\sigma$-strate 
de $\overline{M}$)
s'identifie, localement pour la topologie {\'e}tale, au morphisme d'immersion toriques 
$$\widetilde{S}_{\mathcal{C},\tau}\rightarrow S_{\mathcal{C},\sigma}.$$

\begin{rque} 
1) Le qualificatif ``{faible}'' fait r{\'e}f{\'e}rence au fait que la construction de $P$
{\`a} partir de $\widetilde{G}$, bien que du type de celle de Mumford (compl{\'e}tion, quotient par les 
p{\'e}riodes, puis alg{\'e}brisation),
ne suppose pas que le sch{\'e}ma $\widetilde{P}$ associ{\'e} au tore d{\'e}ploy{\'e} $\widetilde{G}$ contienne ce tore,
(on a encore cependant $\widetilde{G}_\eta=\widetilde{P}_\eta$).   

2) $j:\mathfrak{G}^s\hookrightarrow P$ n'est pas une immersion toro{\"\i}dale au-dessus de $\overline{M}$.

\end{rque}

On pose $\overline{\mathcal{A}^s}=P$. 
Il reste {\`a} v{\'e}rifier les {\'e}nonc{\'e}s $(\mathrm{v})$ et $(\mathrm{vi})$ du th{\'e}or{\`e}me.
A partir de $j:\mathfrak{G}^s\hookrightarrow \overline{\mathcal{A}^s}$, on 
obtient 
$\quad j^*:\Omega^1_{\overline{\mathcal{A}^s}/\overline{M}}
(\dlog \infty)\rightarrow 
\Omega^1_{\mathfrak{G}^s/\overline{M}}$ qui induit un 
isomorphisme sur 
$\overline{f_s}^*e^* \Omega^1_{\mathfrak{G}^s/\overline{M}}=
\overline{f_s}^*({\underline{\omega}_{\mathfrak{G}/\overline{M}}}^{\oplus s})$,
d'o{\`u} le $(\mathrm{v})$.

 Le  $(\mathrm{vi})$ se  d{\'e}duit {\`a} partir du  $(\mathrm{v})$ et du cup-produit 
\begin{equation}
\quad \bigwedge^a R^1\overline{f_s}_*
\left(\mathcal{O}_{\overline{\mathcal{A}^s}}\right)\rightarrow 
R^a \overline{f_s}_*(\mathcal{O}_{\overline{\mathcal{A}^s}})
\end{equation}

Pour  montrer que cette fl{\`e}che est un  isomorphisme on se 
 ram{\`e}ne d'abord par  compl{\'e}tion aux bons mod{\`e}les formels compacts 
(voir  VI.1.11 de \cite{FaCh}), qui permettent de remplacer le morphisme 
$\overline{f_s}: \overline{\mathcal{A}^s}\rightarrow \overline{M}$ par 
les morphismes d'immersions toriques 
$$g:\widetilde{S}_{\mathcal{C},\tau} \rightarrow S_{\mathcal{C},\sigma}.$$

On exploite alors  l'action de  $\widetilde{G}$ sur 
$R^a g_*(\mathcal{O}_{\widetilde{S}_{\mathcal{C},\tau}})$ qui
 permet de calculer la cohomologie des immersions toriques comme au 
bas de la page 208 de  \cite{FaCh}. \hfill $\square$

\medskip
Les points $(\mathrm{v})$ et $(\mathrm{vi})$ du th{\'e}or{\`e}me pr{\'e}c{\'e}dent 
sont en partie cons{\'e}quence du fait plus g{\'e}n{\'e}ral suivant que  le complexe 
$R\overline{f}_*\Omega^\bullet_{\overline{\mathcal{A}}/\overline{M}}
(\dlog\,\infty)$ ne d{\'e}pend pas
du choix de la compactification toro{\"\i}dale $\overline{\mathcal{A}}$
(voir le lemme VI.3.4 de \cite{FaCh} qui se transpose sans changement
{\`a} notre cas). On en d{\'e}duit en particulier que le fibr{\'e} 
$\overline{\mathcal{H}}^1_{\dR}$ ne d{\'e}pend pas du choix de 
la compactification toro{\"\i}dale $\overline{\mathcal{A}}$ 
au-dessus de $\overline{M}$ et  qu'il est muni d'une action de $\mathfrak{o}$.
En fait, si on pose  $j_M:M\hookrightarrow \overline{M}$, alors 
$\overline{\mathcal{H}}^1_{\dR}$ s'identifie au sous-faisceau de 
$j_{M*}\mathcal{H}^1_{\dR}(\mathcal{A}/M)$ 
des sections ${\mathfrak{G}}$-invariantes 
de $\mathcal{H}^1_{\dR}(\mathfrak{G}/\overline{M})$.

La suite spectrale de Hodge vers de Rham 
logarithmique  fournit une suite exacte courte

$$0\rightarrow \overline{f}_*\Omega_{\overline{\mathcal{A}}/\overline{M}}(\dlog \,\infty)\rightarrow 
\overline{\mathcal{H}}^1_{\dR}\rightarrow 
R^1\overline{f}_*\mathcal{O}_{\overline{\mathcal{A}}} \rightarrow 0$$

qui est, elle-aussi, ind{\'e}pendante de $\overline{\mathcal{A}}$.
La filtration de Hodge sur $\overline{\mathcal{H}}^1_{\dR}$ est 
donc ind{\'e}pendante de $\overline{\mathcal{A}}$. On a la 
premi{\`e}re partie de la 

\bigskip
\begin{prop}\label{dRind} 
{\'E}tant donn{\'e} un $\Gamma$-{\'e}ventail complet $\Sigma$ de $C_+$, le 
fibr{\'e} $\overline{\mathcal{H}}^1_{\dR}$ ne d{\'e}pend pas du choix de 
la compactification toro{\"\i}dale $\overline{\mathcal{A}}$ au-dessus de 
$\overline{M}$; il en est de m{\^e}me pour sa filtration de Hodge et 
pour sa connexion logarithmique prolongeant la connexion de Gauss-Manin.
Il est muni d'une action naturelle de ${\mathfrak{G}}$ et de $\mathfrak{o}$.
La connexion de Gauss-Manin logarithmique  est compatible avec 
la fl{\`e}che de Kodaira-Spencer  
$$\underline{\omega}_{\mathfrak{G}/\overline{M}}\rightarrow 
(\underline{\omega}_{\mathfrak{G}/\overline{M}}^\vee\otimes_{\mathfrak{o}}
 \mathfrak{cd}^{-1}) \otimes \Omega_{\overline{M}}(\dlog\,\infty).$$
\end{prop}

\noindent{\bf D{\'e}monstration :}
On d{\'e}montre que la connexion de Gauss-Manin poss{\`e}de un prolongement 
ind{\'e}pendant de $\overline{\mathcal{A}}$ et que ce prolongement est unique.
Pour une compactification donn{\'e}e $\overline{\mathcal{A}}$, 
on peut d{\'e}finir la connexion de Gauss-Manin logarithmique 
(voir la section 2 de \cite{KaOd} dans le cas  non-logarithmique) comme suit.
Posons
$ \Fil^i \Omega_{\overline{\mathcal{A}}}^{\bullet }(\dlog\,\infty) =
\im(\overline{f}^*\Omega_{\overline{M}}^i(\dlog\,\infty)\otimes 
\Omega_{\overline{\mathcal{A}}}^{\bullet -i}(\dlog\,\infty)\rightarrow  
\Omega_{\overline{\mathcal{A}}}^\bullet(\dlog\,\infty))$
 et consid{\'e}rons  la suite exacte de complexes
\begin{equation}\label{GM}
0\rightarrow \Fil^1/\Fil^2\rightarrow \Fil^0/\Fil^2\rightarrow 
\Fil^0/\Fil^1\rightarrow 0
\end{equation}

La connexion de Gauss-Manin s'identifie alors au morphisme connectant
$$R^1\overline{f}_* \gr^0\rightarrow R^2\overline{f}_*\, \gr^1.$$

 Si l'on pose 
$\Fil^i_{\mathfrak{G}}=\Fil^i \Omega_{\mathfrak{G}}^\bullet(\dlog\,\infty)=
\im( \overline{f}^*\Omega_{\overline{M}}^i(\dlog\,\infty)\otimes
\Omega_{\mathfrak{G}}^{\bullet-i}(\dlog\,\infty)
\rightarrow \Omega_{\mathfrak{G}}^\bullet(\dlog\,\infty))$, 
o{\`u} le $\dlog\,\infty$  n'est relatif qu'aux p{\^o}les le long du diviseur 
vertical $\overline{f}^*\infty$ de $\mathfrak{G}$, 
on peut identifier (\ref{GM}) {\`a} la sous-suite (exacte) des 
${\mathfrak{G}}$-invariants de
$$0\rightarrow \Fil^1_{\mathfrak{G}}/\Fil^2_{\mathfrak{G}}\rightarrow 
\Fil^0_{\mathfrak{G}}/\Fil^2_{\mathfrak{G}}\rightarrow 
\Fil^0_{\mathfrak{G}}/\Fil^1_{\mathfrak{G}} \rightarrow 0,$$

\noindent  qui ne d{\'e}pend pas de $\overline{\mathcal{A}}$. 
Encore une fois, ceci r{\'e}sulte de ce que
la connexion de Gauss-Manin
sur $\mathcal{A}$ est $\mathcal{A}$-invariante; 
elle se prolonge donc localement de fa{\c c}on unique via l'identification
 $\overline{\mathcal{H}}^1_{\dR}=\mathcal{H}^1_{\dR}(\mathfrak{G}/\overline{M})^{\mathfrak{G}}$.
\hfill $\square$

\vspace{-4mm}
\section[Applications des compactifications toro{\"\i}dales arithm{\'e}tiques.]{Applications de la compactification toro{\"\i}dale arithm{\'e}tique}

\subsubsection*{Irr{\'e}ductibilit{\'e} du sch{\'e}ma $M\otimes \F_p$ ($p\nmid \Delta$).}

Le sch{\'e}ma $M$ est g{\'e}om{\'e}triquement irr{\'e}ductible sur $\Z[\frac{1}{\Delta}]$.
Il en est de m{\^e}me pour le $\mathcal{T}_1$-torseur
 $\mathfrak{M}$ sur $M$.

La d{\'e}monstration est la m{\^e}me que dans \cite{FaCh} IV.5.10 :
la  fibre g{\'e}n{\'e}rique de $M$ est g{\'e}om{\'e}triquement connexe 
 par la description transcendante
 de $M^{\an }$ et le principe GAGA; il en est de m{\^e}me pour 
la fibre g{\'e}n{\'e}rique d'une compactification toro{\"\i}dale  
$\overline{M}$. 
Soit $a:\overline{M}\rightarrow S=\Spec(\Z[\frac{1}{\Delta}])$ 
le morphisme structural. Ce morphisme est lisse donc plat. 
Il est propre donc $a_*\mathcal{O}_{\overline{M}}$ est un
$\Z[\frac{1}{\Delta}]$-module de type fini. Par platitude, 
ce module est libre de rang $r$. 
En passant {\`a} la fibre g{\'e}n{\'e}rique, on voit que $r=1$ parce que cette fibre est connexe (et propre).
Le Th{\'e}or{\`e}me de Connexit{\'e} de Zariski montre que la condition
$a_*\mathcal{O}_{\overline{M}}=  \Z[\frac{1}{\Delta}]$ entra{\^\i}ne
la connexit{\'e} des  fibres $\overline{M}\otimes \F_p$  ($p\nmid \Delta$). 
La lissit{\'e} de $\overline{M}\otimes \F_p$
entra{\^\i}ne alors l'irr{\'e}ductibilit{\'e} g{\'e}om{\'e}trique  de $M\otimes \F_p$.

\vspace{-4mm}

\subsubsection*{Prolongement des fibr{\'e}s automorphes.}

\label{prolfaisceaux}
On peut reprendre la construction de  fibr{\'e}s automorphes de 
la partie \ref{fmh}  {\`a} l'aide de torseurs sur $\overline{M}$.
Dans ce paragraphe on se place au-dessus de $\Z[\frac{1}{\Delta}]$.

On commence par le cas de $\overline{M^1}$.
Fixons un {\'e}ventail $\Gamma$-admissible lisse $(\Sigma^{\mathcal{C}})_{\mathcal{C}}$
et une compactification toro{\"\i}dale $\overline{M^1}=M_{\Sigma}^1$. Soit
$\overline{f}:\mathfrak{G}\rightarrow \overline{M^1}$ 
le sch{\'e}ma semi-ab{\'e}lien construit sur $\overline{M^1}$. Rappelons
que l'on pose $\underline{\omega}_{\mathfrak{G}/\overline{M^1}}=
\overline{e}^*\Omega_{\mathfrak{G}/\overline{M^1}}$, o{\`u}
$\overline{e}:\overline{M^1}\rightarrow \mathfrak{G}$ d{\'e}signe la
section unit{\'e}. On a 
$\underline{\omega}_{\mathfrak{G}/\overline{M^1}}\cong
\mathcal{O}_{\overline{M^1}}\otimes \mathfrak{o}$ localement pour la
topologie de Zariski (voir \cite{dimdg}). 

Alors $\overline{\mathfrak{M}^1}:=
\Isom_{\overline{M^1}}\Big(\mathcal{O}_{\overline{M^1}}\otimes \mathfrak{o},
\underline{\omega}_{\mathfrak{G}/\overline{M^1}}\Big)$ est 
 un $\mathcal{T}_1$-torseur de Zariski sur $\overline{M^1}$.

Un point de $\overline{\mathfrak{M}^1}$ est 
un couple $(x,\omega)$ constitu{\'e} d'un point
$x\in \overline{M^1}$, et d'une 
$(\mathfrak{o}\otimes\mathcal{O}_{\overline{M^1}})$-base 
$\omega$ de $\underline{\omega}_{\mathfrak{G}/\overline{M^1}}$.

Comme dans la partie  \ref{fmh}, le $\mathfrak{o}$-fibr{\'e} 
inversible $\underline{\omega}_{\mathfrak{G}/\overline{M^1}}$
descend en un  $\mathfrak{o}$-fibr{\'e}  inversible sur $\overline{M}$,
not{\'e} encore $\underline{\omega}$. Alors 
le $\mathcal{T}_1$-torseur de Zariski 
$$\overline{\mathfrak{M}}=
\Isom_{\overline{M}}\Big(\mathcal{O}_{\overline{M}}\otimes \mathfrak{o},
\underline{\omega}\Big)$$ prolonge le $\mathcal{T}_1$-torseur 
$\mathfrak{M}$ sur $M$, d{\'e}fini dans la partie \ref{fmh}.

\medskip
Soit  $\mathcal{O}'$ l'anneau des entiers de 
$F^{\gal }(\sqrt{\epsilon},\epsilon\in \mathfrak{o}_{D+}^\times)$.

Pour tout $\Z[\frac{1}{\Delta}]$-sch{\'e}ma $Y$, on pose
$Y'=Y\times \Spec(\mathcal{O}'[\frac{1}{\Delta}])$.

\medskip

 On a  un foncteur 
$\overline{\mathcal{F}}_{\mathcal{T}_1}$ des 
repr{\'e}sentations alg{\'e}briques du 
 $\mathcal{O'}[{\frac{1}{\Delta}}]$-sch{\'e}ma en groupes $\mathcal{T}'_1$, 
vers les fibr{\'e}s d{\'e}composables en fibr{\'e}s inversibles sur 
$\overline{M}'$, qui {\`a} $W$ associe le produit contract{\'e}
 $\overline{\mathfrak{M}}'\stackrel{\mathcal{T}'_1}{\times} 
W=:\overline{\mathcal{W}}$

Si $W_{\mathrm{st}}$ est la repr{\'e}sentation standard de $\mathcal{T}'_1$ ({i.e.}  
$\mathfrak{o}\otimes\mathcal{O}'[{\frac{1}{\Delta}}]$
avec action de $\mathfrak{o}^\times$),  on a 
$\overline{\mathcal{W}}_{\mathrm{st}}=\underline{\omega}^\vee$.
Pour tout caract{\`e}re $\kappa$,  vu comme $\mathcal{O}'$-repr{\'e}sentation
de $\mathcal{T}'_1$, on obtient le prolongement canonique de 
$\underline{\omega}^\kappa=\mathcal{W}_{1,-\kappa}$ {\`a} $\overline{M}$,
comme
$\overline{\mathfrak{M}}'\stackrel{\mathcal{T}'_1}{\times}
\mathcal{O}'[{\frac{1}{\Delta}}](-\kappa)$.

Pour all{\'e}ger les notations, on note encore $\underline{\omega}^\kappa$ 
le prolongement 
canonique  de $\underline{\omega}^\kappa$ {\`a} $\overline{M}$.

\subsubsection*{Principe de Koecher.} Dans toute cette partie on suppose $F\neq \Q$. 
Pour tout poids  $\kappa\in\Z[J_F]$ et 
pour toute $\Z[\frac{1}{\Delta}]$-alg{\`e}bre $R$ contenant les valeurs de
$\kappa$, on a :

\begin{theo}(Principe de Koecher \cite{rapoport} 4.9) \label{koecher}
$$\Gamma\big(M \times \Spec(R),
\underline{\omega}^\kappa
\big)=\Gamma\big(\overline{M}\times \Spec(R),
\underline{\omega}^\kappa \big)$$
\end{theo}

\noindent {\bf D{\'e}monstration : } Il suffit de v{\'e}rifier l'holomorphie 
d'une section globale de $\underline{\omega}^\kappa$ 
le long du diviseur {\`a} l'infini de 
$\overline{M}$. Il suffit donc de montrer que
pour toute $(R,\mathfrak{n})$-pointe $\mathcal{C}$, 
les sections globales m{\'e}romorphes du pull-back de $\underline{\omega}^\kappa$ 
sur $S_{\Sigma^\mathcal{C}}^\wedge\times 
\Spec(R)$ qui sont $\mathfrak{o}_{\mathcal{C}}^\times$-invariantes 
sont holomorphes.
Le  pull-back $\underline{\omega}_{\mathcal{C}}$ 
de $\underline{\omega}$ sur $S_{\Sigma^\mathcal{C}}^\wedge \times  
\Spec(R)$ est canoniquement isomorphe {\`a} 
$\mathfrak{a}\otimes \mathcal{O}_{S_{\Sigma^\mathcal{C}}^\wedge}
\otimes R$.
On peut donc identifier une section m{\'e}romorphe de $\underline{\omega}_{\mathcal{C}}$
{\`a} une s{\'e}rie $f_{\mathcal{C}}=\sum_{\xi\in X} a_\xi q^\xi$, telle que 
pour tout $(u,\epsilon)\in \mathfrak{o}_{\mathcal{C}}^\times$ d'apr{\`e}s
(\ref{units})
on a  $a_{u^2\epsilon\xi}=u^\kappa e^{2i\pi\Tr_{F/\mathbb{Q}}
(\xi u\xi_{u,\epsilon}^*)}a_{\xi}$.
Supposons que $f_{\mathcal{C}}$ ne soit pas holomorphe. Il 
existe donc $\xi_0$ non totalement positif tel que $a_{\xi_0}\neq 0$. 
C'est donc qu'il existe  $\xi_0^*\in X^*_{\R+}$ avec 
$\Tr_{F/\Q}(\xi_0 \xi_0^*)$ strictement n{\'e}gatif. Comme $F\neq\Q$,
on peut choisir des unit{\'e}s  $(u,\epsilon)\in \mathfrak{o}_{\mathcal{C}}^\times$ 
de mani{\`e}re {\`a} rendre la quantit{\'e} $\Tr_{F/\Q}(u^2\epsilon \xi_0 \xi_0^*)$ 
arbitrairement  proche de $-\infty$. 
 Soit $\sigma$ un c{\^o}ne poly{\'e}dral  de $\Sigma^{\mathcal{C}}$
contenant $\xi_0^*$. Par d{\'e}finition de $S_{\sigma}^\wedge$, on 
voit que $f_{\mathcal{C}}$ n'est pas
m{\'e}romorphe sur $S_{\sigma}^\wedge$, ce qui est absurde. 

\subsubsection*{$q$-d{\'e}veloppement.}
Pour all{\'e}ger les notations on se cantonne au cas des pointes non-ramifi{\'e}es.
Voir la partie 8 de \cite{dimdg} pour le cas g{\'e}n{\'e}ral.
Soit $\kappa\in\Z[J_F]$ et soit $R$ une  $\mathcal{O}'[\frac{1}{\Delta}]$-alg{\`e}bre,
o{\`u} $\mathcal{O}'$ d{\'e}signe l'anneau des entiers d'une cl{\^o}ture galoisienne de $F$.
Soit $\mathcal{C}$ une $(R,\mathfrak{n})$-pointe uniformis{\'e}e non ramifi{\'e}e; posons
$$R_\mathcal{C}^{(\kappa)}(R):=
\left\{ \sum_{\xi\in X_+\cup\{0\}} a_\xi q^\xi\Big|
a_\xi\in R, \enspace a_{u^2\epsilon\xi}=\epsilon^{\kappa/2}u^\kappa a_\xi
, \enspace \forall (u,\epsilon) \in\mathfrak{o}_{\mathcal{C}}^\times\right\}$$

Comme le  pull-back $\underline{\omega}_{\mathcal{C}}$ de 
$\underline{\omega}$ sur $S_{\Sigma^\mathcal{C}}^\wedge \times  
\Spec(R)$ est canoniquement isomorphe {\`a} 
$\mathfrak{a}\otimes \mathcal{O}_{S_{\Sigma^\mathcal{C}}^\wedge}
\otimes R$, on a 
$$\underline{\omega}_{\mathcal{C}}^{\kappa}= 
\left(\mathfrak{a}\otimes \mathcal{O}_{S_{\Sigma^\mathcal{C}}^\wedge}
\otimes  R \right)^{-\kappa}$$

Or 
$\left(\mathfrak{a}\otimes \mathcal{O}_{S_{\Sigma^\mathcal{C}}^\wedge}
\otimes  R\right)^{-\kappa}=
\left(\mathfrak{a}\otimes \mathcal{O}'\right)^{-\kappa}
\otimes_{\mathcal{O}'} \left(\mathfrak{o}\otimes \mathcal{O}_{S_{\Sigma^\mathcal{C}}^\wedge}
\otimes  R\right)^{-\kappa} $.

Notons que 
$R_\mathcal{C}^{(\kappa)}(R)=\left(\mathfrak{o}\otimes 
\mathcal{O}_{S_{\Sigma^\mathcal{C}}^\wedge}
\otimes  R\right)^{-\kappa}$ est un 
$\mathcal{O}_{S_{\Sigma^\mathcal{C}}^\wedge} \otimes  R$-module 
inversible et $\mathfrak{a}^{(\kappa)}:=
(\mathfrak{a}\otimes \mathcal{O}')^{-\kappa}$ est un 
$\mathcal{O}'$-module inversible. 

On  peut donc associer {\`a} toute 
forme modulaire de Hilbert $f$ de poids $\kappa$,  
niveau $\Gamma$ d{\'e}finie sur $R$, un {\'e}l{\'e}ment $f_\mathcal{C}\in 
\mathfrak{a}^{(\kappa)}  \otimes_{\mathcal{O}'} R_\mathcal{C}^{(\kappa)}(R)$.

\begin{defin}
La s{\'e}rie $f_\mathcal{C}$  est appel{\'e}e le $q$-d{\'e}veloppement de la forme $f$ 
en la pointe $\mathcal{C}$.
On note $\ev_{\mathcal{C},\kappa}$ l'application $f\mapsto f_\mathcal{C}$.
\end{defin}

 Le principe du $q$-d{\'e}veloppement en une pointe $\mathcal{C}$ non-ramifi{\'e}e s'{\'e}nonce alors :

\begin{prop} Pour toute  $\mathcal{O}'[\frac{1}{\Delta}]$-alg{\`e}bre
 $R$,

1) l'application
 $$\ev_{\mathcal{C},\kappa}:G_\kappa(\mathfrak{c},\mathfrak{n};R)\rightarrow 
 \mathfrak{a}^{(\kappa)}  \otimes_{\mathcal{O}'} R_\mathcal{C}^{(\kappa)}(R)$$
 est
injective,

2)  pour toute inclusion $R\subset R'$ d'alg{\`e}bres, 
si $f\in G_\kappa(\mathfrak{c},\mathfrak{n};R')$ et $f_{\mathcal{C}}\in 
\mathfrak{a}^{(\kappa)}  \otimes_{\mathcal{O}'} R_\mathcal{C}^{(\kappa)}(R)$,
alors $f\in G_\kappa(\mathfrak{c},\mathfrak{n};R)$.
\end{prop}

L'{\'e}nonc{\'e} 2) dans le cas de l'anneau nul $R=0$ redonne 1).

\noindent{\bf D{\'e}monstration : } Les deux {\'e}nonc{\'e}s r{\'e}sultent du suivant:
soit $R$ un groupe ab{\'e}lien; l'application $f\mapsto f_{\mathcal{C}}$:
$$\mathrm{H}^0(\overline{M},\underline{\omega}^\kappa\otimes R)\rightarrow 
\mathfrak{a}^{(\kappa)}  \otimes
R[[q^\xi;\xi\in X_+\cup\{0\}]]$$
est injective
(on  utilise le principe de Koecher pour passer {\`a} 
 $\overline{M}$). 

Par commutation des deux membres aux limites inductives, on se ram{\`e}ne ais{\'e}ment au cas
$R=\Z$ ou $R=\Z/r\Z$.
Par l'irr{\'e}ductibilit{\'e} g{\'e}om{\'e}trique 
de $\mathfrak{M}$ sur
$\Z[\frac{1}{\Delta}]$, une section globale $s$ de 
$\underline{\omega}^\kappa$ sur $\overline{M}\otimes \Z/r\Z$
est nulle, si et seulement si son pull-back {\`a} la compl{\'e}tion de 
$\overline{M}\otimes \Z/r\Z$ le long d'un
diviseur est nul. Soit $\mathcal{C}$ une pointe; 
$S_{\Sigma^{\mathcal{C}}}^\wedge/
\mathfrak{o}_{\mathcal{C}}^\times$ s'identifie {\`a} 
la compl{\'e}tion de $M_\Sigma$ le long de 
$\pi^{-1}(\mathcal{C})$.
 Il suffit donc que le pull-back de $s$ {\`a} 
$S_{\Sigma^{\mathcal{C}}}^\wedge/
\mathfrak{o}_{\mathcal{C}}^\times$ soit nul. C'est-{\`a}-dire que $\ev_{\mathcal{C},\kappa}(s)$ soit nul.

\begin{rque} 1) L'application $\ev_{\mathcal{C}}$,
somme des $\ev_{\mathcal{C},\kappa}$,
$$G(\mathfrak{c},\mathfrak{n};R)=
\bigoplus_{\kappa\in \Z[J_F]}G_\kappa(\mathfrak{c},\mathfrak{n};R)
\rightarrow R\otimes \Z[[q^\xi;\xi\in X_+\cup\{0\}]]$$
n'est pas injective en g{\'e}n{\'e}ral comme le montre l'exemple de 
$F=\Q$, $R=\F_p$, $\Gamma=\SL_2(\Z)$:
le noyau de $\ev_{\mathcal{C}}\otimes \id_{\F_p}$ 
$$ \bigoplus_{\kappa\in \Z}G_\kappa(\mathfrak{c},\mathfrak{n};\F_p)
\rightarrow \F_p[[q]]$$
est l'id{\'e}al engendr{\'e} par $E_{p-1}-1$.

2) Pour $F$ totalement r{\'e}el quelconque, le noyau de la 
fl{\`e}che $\ev_{\mathcal{C}}\otimes \F_p$  {\`a} {\'e}t{\'e} calcul{\'e}e par Goren
(voir \cite{goren} Chap.5, Corollaire 4.5).

3) En fait, si $R$ est une $\Z$-alg{\`e}bre  sans torsion, 
$\ev_{\mathcal{C}}$ est injective
gr{\^a}ce au th{\'e}or{\`e}me de Dedekind d'ind{\'e}pendance lin{\'e}aire des caract{\`e}res distincts.
\end{rque}

\subsubsection*{Prolongement de fibr{\'e}s filtr{\'e}s et {\`a} connexions.}

Fixons un {\'e}ventail admissible principalement polaris{\'e} 
$(\widetilde{\Sigma},\widetilde{\phi})$  pour $\mathcal{A}$ au-dessus de 
l'{\'e}ventail admissible $ \Sigma=(\Sigma^{\mathcal{C}})$ fix{\'e} pour $M^1$; 
on a ainsi un morphisme de compactifications toro{\"\i}dales $\overline{f}:\overline{\mathcal{A}}\rightarrow \overline{M^1}$ prolongeant $f:\mathcal{A}\rightarrow M^1$. 

Dans ce qui suit, on posera pour abr{\'e}ger 
$\Omega^\bullet_{\overline{\mathcal{A}}/\overline{M^1}}(\dlog \infty)=
\Omega^\bullet_{\overline{\mathcal{A}}/\overline{M^1}}(\dlog \infty_{\overline{\mathcal{A}}/\overline{M^1}})$ et 
 $\overline{\mathcal{H}}^1_{\dR}=
R^1\overline{f}_*\Omega_{\overline{\mathcal{A}}/\overline{M^1}}^\bullet
(\dlog\,\infty)$.
Ce dernier faisceau est localement libre de rang $2$ sur 
$\mathfrak{o}\otimes \mathcal{O}_{\overline{M^1}}$.
En outre, il est muni d'une filtration {\`a} deux crans donn{\'e}e par 
la suite spectrale de Hodge vers de Rham :
$$ 0\rightarrow \overline{f}_*\Omega_{\overline{\mathcal{A}}/\overline{M^1}}
(\dlog \infty)\rightarrow \overline{\mathcal{H}}_{\dR}^1\rightarrow 
R^1\overline{f}_*\mathcal{O}_{\overline{\mathcal{A}}} \rightarrow 0$$

Par le  th{\'e}or{\`e}me  \ref{thmkugasato}(vi), on a des 
isomorphismes canoniques de faisceaux
$$ \overline{f}_*\Omega_{\overline{\mathcal{A}}/\overline{M^1}}(\dlog \infty)\cong \underline{\omega}_{\mathfrak{G}/\overline{M^1}}
\text{ et }
R^1\overline{f}_*\mathcal{O}_{\overline{\mathcal{A}}}\cong \underline{\omega}_{\mathfrak{G}/\overline{M^1}}^{\vee}\otimes\mathfrak{cd}^{-1}.$$

La filtration de $\overline{\mathcal{H}}_{\dR}^1$ se r{\'e}{\'e}crit donc  
$\Fil^0\,\overline{\mathcal{H}}_{\dR}^1=\overline{\mathcal{H}}_{\dR}^1$,
 $\Fil^1\,\overline{\mathcal{H}}_{\dR}^1= 
\underline{\omega}_{\mathfrak{G}/\overline{M^1}}$ et 
$\gr^0\,\overline{\mathcal{H}}_{\dR}^1= 
\underline{\omega}_{\mathfrak{G}/\overline{M^1}}^\vee\otimes\mathfrak{cd}^{-1}$.

Comme dans la partie  \ref{fmh} le fibr{\'e} $\overline{\mathcal{H}}_{\dR}^1$
descend en un fibr{\'e} sur $\overline{M}$ jouissant aux m{\^e}mes propri{\'e}t{\'e}s.

\medskip
On d{\'e}finit un $\mathcal{D}$-torseur $\overline{\mathfrak{M}}_{\mathcal{D}}=
\mathrm{Isom}^{\mathcal{D}}_{\mathfrak{o}\otimes\mathcal{O}_{\overline{M}}}
(\mathfrak{o}\otimes\mathcal{O}_{\overline{M}},
\wedge^2_{\mathfrak{o}\otimes\mathcal{O}_{\overline{M}}}
\overline{\mathcal{H}}^1_{\dR})$
 au-dessus de $\overline{M}$, dont les $S$-points sont ceux de
$\mathrm{Isom}_{\mathfrak{o}\otimes\mathcal{O}_{\overline{M}}}
(\mathfrak{o}\otimes\mathcal{O}_{\overline{M}},
\wedge^2_{\mathfrak{o}\otimes\mathcal{O}_{\overline{M}}}
\overline{\mathcal{H}}^1_{\dR})$
induisant via $\lambda$ un {\'e}l{\'e}ment de 
$\mathcal{D}(\mathcal{O}_S)$ dans $(\mathfrak{o}\otimes \mathcal{O}_S)^\times$.

\medskip
On d{\'e}finit un $\mathcal{B}$-torseur $\overline{\mathfrak{M}}_{\mathcal{B}}
 \overset{\mathcal{B}}{\rightarrow}\overline{M}$ comme le produit fibr{\'e} de $\overline{\mathfrak{M}}_{\mathcal{D}}$ et de \\
$\mathrm{Isom}^{\mathrm{fil}}_{\mathfrak{o}\otimes\mathcal{O}_M}
((\mathfrak{o}\otimes\mathcal{O}_M)^2,\overline{\mathcal{H}}^1_{\dR})$
  au-dessus de 
$\mathrm{Isom}_{\mathfrak{o}\otimes\mathcal{O}_{\overline{M}}}
(\mathfrak{o}\otimes\mathcal{O}_{\overline{M}},
\wedge^2_{\mathfrak{o}\otimes\mathcal{O}_{\overline{M}}}
\overline{\mathcal{H}}^1_{\dR})$.

Il d{\'e}finit un foncteur $\overline{\mathcal{F}}_{\mathcal{B}'}$
de la cat{\'e}gorie des repr{\'e}sentations alg{\'e}briques du 
 $\mathcal{O}'[{\frac{1}{\Delta}}]$-sch{\'e}ma en groupes $\mathcal{B}'$ vers celle 
des fibr{\'e}s sur $\overline{M}'$ qui sont des extensions successives de
fibr{\'e}s inversibles.  Le foncteur est donn{\'e} par
$V \mapsto \overline{\mathcal{V}}:= 
\overline{\mathfrak{M}}_{\mathcal{B'}} \stackrel{\mathcal{B'}}{\times} V$.

Si $V_{\mathrm{st}}=\Big(\mathfrak{o}[{\frac{1}{\Delta}}]\Big)^2$ est la repr{\'e}sentation standard 
de $\mathcal{B}$, on a
$\overline{\mathcal{V}}_{\mathrm{st}}=\overline{\mathcal{H}}^1_{\dR}$.

\begin{defin}
Pour tout  poids alg{\'e}brique $\kappa$  et 
$n, m\in \Z[J_F]$ comme dans la  d{\'e}finition \ref{alg}, on 
note $\overline{\mathcal{V}}_n$ et $\overline{\mathcal{W}}_{n,c}$
les prolongements {\`a} $\overline{M}'$ construits {\`a} l'aide de
  $\overline{\mathcal{F}}_{\mathcal{B}'}$  des fibr{\'e}s $\mathcal{V}_n$ et 
$\mathcal{W}_{n,c}$ des d{\'e}finitions  \ref{Vfil} et \ref{W}.
\end{defin}

\begin{rque}
1) Pour toute $\mathcal{O}'[{\frac{1}{\Delta}}]$-repr{\'e}sentation alg{\'e}brique $V$ 
de $\mathcal{B}'$, le fibr{\'e} $\overline{\mathcal{V}}$ est 
le prolongement de Mumford de $\mathcal{V}$, 
c'est-{\`a}-dire que son pull-back {\`a} toute carte 
locale donn{\'e}e par une immersion torique 
$S_\mathcal{C}\stackrel{j}{\hookrightarrow} 
S_{\mathcal{C},\sigma_\alpha^{\mathcal{C}}}$
est engendr{\'e} par les sections $S_\mathcal{C}$-invariantes de $j_*\mathcal{V}$.
En effet, c'est vrai pour $\overline{\mathcal{H}}_{\dR}^1$ par la 
proposition \ref{dRind} et 
donc pour tout les fibr{\'e}s  obtenus par pl{\'e}thysme {\`a} partir de
 $\overline{\mathcal{H}}_{\dR}^1$.
Ceci implique par exemple, que pour tout un poids alg{\'e}brique 
$\kappa$ et $m,n\in\Z[J_F]$, comme dans la d{\'e}finitions \ref{alg},  
le foncteur $\overline{\mathcal{F}}_{\mathcal{B}'}$ 
fournit sur $\C$ ( sur $\Q$ ou  sur $\Q_p$) le prolongement de Mumford de 
$\Sym^n\mathcal{H}^1_{\dR} \otimes (\wedge^2 \mathcal{H}^1_{\dR})^{\otimes m}$ 
et de $\underline{\omega}^\kappa$.

2) Rappelons que sur une cl{\^o}ture galoisienne $F^{\gal}$ de $F$, on a en posant
$$\mathcal{H}^1_{dR,\tau,F^{\gal}}=\mathcal{H}^1_{\dR}\otimes_{F,\tau} 
F^{\gal},\quad 
\underline{\omega}^\tau=\underline{\omega}\otimes_{F,\tau} F^{\gal},$$
$$\Sym^n\mathcal{H}^1_{dR,F^{\gal}}=\bigotimes_\tau 
\Sym^{n_\tau}\mathcal{H}^1_{dR,\tau,F^{\gal}} \text{, et }
\underline{\omega}^\kappa=\bigotimes_\tau  \Sym^{k_\tau}_{F^{\gal}}
\underline{\omega}^\tau.$$

3) Soit $p$ premier ne divisant pas $\Delta$.
Le foncteur $\overline{\mathcal{F}}_{\mathcal{B}}$ ne donne le prolongement de Mumford des
faisceaux $R^s f_*\Omega^\bullet_{\mathcal{A}/M}$ sur $\overline{M}\otimes \Z_p$ que lorsque $s<p$. 
En effet, pour tout $s<p$, Illusie \cite{illusie2} a montr{\'e} que
$R^s\overline{f}\Omega^\bullet_{\overline{\mathcal{A}}/\overline{M}}(\dlog\,\infty)$ est libre sur $(\mathfrak{o}\otimes 
\mathcal{O}_{\overline{M}})$. Il en r{\'e}sulte que le 
foncteur $\overline{\mathcal{F}}_{\mathcal{B}'}$
fournit le prolongement de Mumford de $M\otimes \Z_p$ {\`a} 
$\overline{M}\otimes \Z_p$ des faisceaux
$\Sym^n \mathcal{H}^1_{\dR}\otimes (\wedge^2 \mathcal{H}^1_{\dR})^{\otimes m} $ 
lorsque $p-1> \sum_\tau (n_\tau+1)$.
\end{rque}

On d{\'e}finit de plus un $G$-torseur de Zariski :
$$\mathfrak{M}_{\nabla}=\Isom_{\mathfrak{o}\otimes \mathcal{O}_{\overline{M}}}
(L_0\otimes \mathcal{O}_{\overline{M}},\overline{\mathcal{H}}_{\dR}^1)$$

On d{\'e}finit ainsi un foncteur de la cat{\'e}gorie des 
repr{\'e}sentations alg{\'e}briques de $G$ sur
$\mathcal{O}'[\frac{1}{\Delta}]$ vers celle des fibr{\'e}s sur $\overline{M}'$ 
munis d'une connexion int{\'e}grable {\`a} singularit{\'e}s logarithmique
 et dont la r{\'e}duction modulo $p$ est 
quasi-nilpotente en chaque $p\nmid \Delta$  (voir \cite{MoTi} Sect.5.2). 
On observera  d'ailleurs que l'utilisation du torseur $\mathfrak{M}_{\nabla}$ 
d{\'e}fini ici aurait simplifi{\'e} substantiellement la pr{\'e}sentation de 
la partie 5.2 cit{\'e}e en rendant la partie 5.2.3 inutile. 

\begin{defin}
Pour tout  poids alg{\'e}brique $\kappa$  et 
$n, m\in \Z[J_F]$ comme dans la  d{\'e}finition \ref{alg}, on 
note $\overline{\mathcal{V}}{}_n^\nabla$ 
le prolongement {\`a} $\overline{M}'$ du fibr{\'e} 
$\mathcal{V}_n^\nabla$ construit dans la d{\'e}finition \ref{Vnabla}.
\end{defin}

\vspace{-4mm}

\subsubsection*{D{\'e}composition de Hodge-Tate de $\mathrm{H}^\bullet(M\otimes \overline{\Q}_p, \mathbb{V})$.}
\label{hodge}

Dans cette section, nous ne consid{\'e}rons que la filtration de Hodge dite aussi 
$F$-filtration (et son gradu{\'e} associ{\'e}).
C'est-{\`a}-dire
que la filtration par le poids dont les gradu{\'e}s sont purs est ici ignor{\'e}e: les gradu{\'e}s 
que nous faisons appara{\^\i}tre sont encore munis d'une filtration par le poids.

\noindent {\bf Sur $\C$ : }  Nous remercions 
H. Hida pour avoir attir{\'e} notre attention sur le point de vue transcendant suivant.
Soit $V$ une $\Q$-repr{\'e}sentation de $G$, de syst{\`e}me local sur 
$M^{\an }$ associ{\'e} $\mathbb{V}$. On a
 $\GL_2(F\otimes \R)/(F\otimes \R)^\times \mathrm{O}_2(F\otimes \R)=
\SL_2(F\otimes \R)/\mathrm{SO}_2(F\otimes \R)$.
Par cette identification on voit  que 
$\left(\begin{array}{cc}1&0\\0&-1\end{array}\right)\cdot z=-\overline{z}$; 

Le groupe de Weyl de $G$, $
W=\mathrm{O}_2(F\otimes \R)/\mathrm{SO}_2(F\otimes \R)\cong 
\{\pm 1\}^{J_F}$ agit donc sur $(M^{\an },\mathbb{V})$ :
si $\epsilon_J=(-1_J,1_{\overline{J}})\in W$ et $z=(z_J,z_{\overline{J}})\in \mathfrak{H}_F$, 
$\epsilon_J \cdot (z,v)=((-\overline{z}_J,z_{\overline{J}}),v)$.

\medskip
Sur $\mathrm{H}^\bullet (M^{\an }, \mathbb{V}_{\C})=\mathrm{H}^\bullet(M^{\an }, 
\mathbb{V})\otimes \C$,  on a donc
l'action de $W$ d'une part et celle de la conjugaison complexe 
$c$ sur les coefficients de l'autre. 
Soit $\epsilon_\tau=(-1_\tau,1^\tau)\in W$. On d{\'e}compose en espaces propres pour
l'action des $\epsilon_\tau\otimes c$:
$$\mathrm{H}^\bullet(M^{\an }, \mathbb{V}_{\C})=\bigoplus_{J\subset J_F} 
\mathrm{H}^{J,\overline{J}}(M^{\an }, \mathbb{V}_{\C})$$

o{\`u}, en notant $\chi_J$ la fonction caract{\'e}ristique d'une partie $J$ de $J_F$, on a 
$$\mathrm{H}^{J,\overline{J}}(M^{\an }, \mathbb{V}_{\C})=
\{x;(\epsilon_\tau\otimes c) (x)=(-1)^{\chi_{\overline{J}}(\tau)}x \}$$

Cette d{\'e}composition est plus fine que la d{\'e}composition donn{\'e}e par la filtration b{\^e}te: pour tout
entier $a$ tel que $0\leq a\leq d$, on a
$$ \gr_{\mathrm{b\hat{e}te}}^a\mathrm{H}^d(M^{\an }, \mathbb{V}_{\C})=
\bigoplus_{J\subset J_F, \mathrm{card}(J)=a} \mathrm{H}^{J,\overline{J}}(M^{\an }, \mathbb{V}_{\C})
$$

\medskip

 Si $F=\Q$ , $J_{\Q}=\{\id_{\Q}\}$ et $V=\Q$, 
la d{\'e}composition de Hodge de $\mathrm{H}^1(M^{\an },{\C})$ en 
$\mathrm{H}^{J_{\Q}, \emptyset}(M^{\an }, {\C}) 
=\mathrm{H}^{1,0}\cong \mathrm{H}^0(M^{\an  *},\Omega_{M^{\an }} 
(\dlog\,\infty))$ et 
$\mathrm{H}^{\emptyset,J_{\Q}}(M^{\an  *},{\C})=\mathrm{H}^{0,1} \cong \mathrm{H}^1(M^{\an },\mathcal{O}_M^{\an })$,
o{\`u} $M^{\an  *}$ d{\'e}signe la compactification toro{\"\i}dale, qui co{\"\i}ncide ici avec la compactification de Satake.

On voit {\`a} l'isomorphisme d'Eichler-Shimura, on voit que la partie Eisenstein
du $\mathrm{H}^1$ est concentr{\'e}e dans le $\mathrm{H}^{1,0}$.

\medskip

Cette d{\'e}composition de nature transcendante a un parall{\`e}le
 alg{\'e}brique semblable au Th.5.5 Chap.VI de
 \cite{FaCh}. La simplicit{\'e} de l'{\'e}criture ci-dessus 
 vient de ce que le groupe d{\'e}riv{\'e} de $G(\R)$ est un produit de copies de $\SL_2$.

\bigskip

\noindent {\bf Sur $\Q_p$ : }
Soit $p$ un nombre premier et $V_n$ la $\Q_p$-repr{\'e}sentation alg{\'e}brique de $G$ 
d{\'e}finie dans (\ref{Vn}). On peut lui associer 
un faisceau lisse $\mathbb{V}_n$ sur $M^1\otimes \Q$ 
et des  fibr{\'e}s $\overline{\mathcal{V}}_n$ (resp.  
$\overline{\mathcal{V}}_n^\nabla$)  sur
$\overline{M^1}\otimes \Q_p$ qui sont 
filtr{\'e}s (resp. {\`a} connexion {\`a} singularit{\'e}s logarithmiques int{\'e}grable et 
quasi-nilpotente). 

 La d{\'e}monstration de la d{\'e}g{\'e}n{\'e}rescence des suites spectrales des Th.5.5 et 6.2
du Chapitre VI de  \cite{FaCh} dont les termes $\mathrm{E}_1$ 
sont donn{\'e}s en termes du
complexe de Bernstein-Gelfand-Gelfand 
(dualis{\'e} et faisceautis{\'e}) s'adapte au cas de la vari{\'e}t{\'e}
de Hilbert. Il est important de noter que c'est la d{\'e}monstration de ce th{\'e}or{\`e}me qui requiert l'utilisation des compactifications toro{\"\i}dales
de toutes les puissances de la vari{\'e}t{\'e} de Kuga-Sato et pas seulement de la puissance premi{\`e}re. En effet, par un  th{\'e}or{\`e}me de Deligne \cite{De-hodge2} Sect. 3.2, 
la suite spectrale de Hodge vers de Rham 
$$\mathrm{E}_1^{i,j}=\mathrm{H}^j(\overline{\mathcal{A}^s}, 
\Omega^i_{\overline{\mathcal{A}^s}}(\dlog\,\infty))
\Rightarrow  \mathrm{H}^{i+j}
(\Omega^\bullet_{\overline{\mathcal{A}^s}}(\dlog\,\infty))$$

\noindent d{\'e}g{\'e}n{\`e}re en $\mathrm{E}_1$. On en d{\'e}duit comme dans
\cite{FaCh} p.234 la d{\'e}g{\'e}n{\'e}rescence 
en $\mathrm{E}_1$ de la suite  spectrale 
$$\mathrm{E}_1^{i,j}= \mathrm{H}^{i+j}(\overline{M^1}, \gr^i_F (R^sf_{s*} 
\Omega^\bullet_{\overline{\mathcal{A}^s}/\overline{M^1}}(\dlog\,\infty)\otimes 
\Omega^\bullet_{\overline{M^1}}(\dlog\,\infty)))
\Rightarrow $$ 
$$\Rightarrow \mathrm{H}^{i+j}(R^sf_{s*} 
\Omega^\bullet_{\overline{\mathcal{A}^s}/\overline{M^1}}(\dlog\,\infty)\otimes 
\Omega^\bullet_{\overline{M^1}}(\dlog\,\infty)),$$

\noindent o{\`u} la $F$-filtration est obtenue en faisant le 
produit tensoriel des deux filtrations de Hodge. En prenant $s=n_0d$, le fibr{\'e} 
$\overline{\mathcal{V}}_n$ est par pl{\'e}thysme (voir \cite{MoTi} Appendice II)
un facteur direct de 
$(R^1f_* \Omega^\bullet_{\overline{\mathcal{A}}/\overline{M^1}}(\dlog\,\infty)
)^{\otimes s}$ qui, par la formule  de  K{\"u}nneth, 
est lui-m{\^e}me  un facteur direct de 
$R^sf_{s*}\Omega^\bullet_{\overline{\mathcal{A}^s}/\overline{M^1}}(\dlog\,\infty)$.
On en d{\'e}duit la d{\'e}g{\'e}n{\'e}rescence 
en $\mathrm{E}_1$ de la suite  spectrale de Hodge vers de Rham
$$\mathrm{E}_1^{i,j}=\mathrm{H}^{i+j}(\overline{M^1},\gr^i_F 
(\overline{\mathcal{V}}_n\otimes 
\Omega^\bullet_{\overline{M^1}}(\dlog\,\infty)))
\Rightarrow \mathrm{H}^{i+j}(\overline{M^1},\overline{\mathcal{V}}_n\otimes 
\Omega^\bullet_{\overline{M^1}}(\dlog\,\infty)). $$

Par le Th{\'e}or{\`e}me de comparaison de Faltings \cite{Fa-jami}, la 
$\Gal(\overline{\Q}_p/\Q_p)$-repr{\'e}sentation $\mathrm{H}^\bullet 
(M^1\otimes \overline{\Q}_p,\mathbb{V}_n)$
est de de Rham et pour toute compactification toro{\"\i}dale 
$\overline{f}:\overline{\mathcal{A}}\rightarrow \overline{M^1}$ de 
$f:\mathcal{A}\rightarrow M^1$, on a un isomorphisme canonique
$$\mathrm{H}^\bullet(M^1\otimes \overline{\Q}_p,\mathbb{V}_n)\otimes B_{\dR}\cong
\mathrm{H}^\bullet(\overline{M^1},\overline{\mathcal{V}}_n\otimes 
\Omega^\bullet_{\overline{M^1}}(\dlog\,\infty))\otimes B_{\dR}.$$

Les poids de Hodge-Tate de  $\mathrm{H}^\bullet 
(M^1\otimes \overline{\Q}_p,\mathbb{V}_n)$ sont donc donn{\'e}s par 
les sauts de la filtration de Hodge sur 
$\mathrm{H}^\bullet(\overline{M^1},\overline{\mathcal{V}}_n\otimes 
\Omega^\bullet_{\overline{M^1}}(\dlog\,\infty))$ venant de la
suite spectrale ci-dessus. Nous allons calculer ces derniers
  comme dans \cite{FaCh} Th.5.5 ou  \cite{MoTi}
{\`a} l'aide d'un sous-complexe  facteur direct  de 
$\overline{\mathcal{V}}_n\otimes 
\Omega^\bullet_{\overline{M^1}}(\dlog\,\infty)$, appel{\'e} le  complexe BGG.
Avant d'{\'e}noncer le th{\'e}or{\`e}me nous allons introduire quelques notations. 
\medskip

On identifie l'ensemble des parties de $J_F$ avec le groupe de Weyl $W$  de $G$, en 
associant {\`a}  $J\subset J_F$ l'{\'e}l{\'e}ment $\epsilon_J=(-1_J,1_{\overline{J}})\in W$.
Pour tout $J\subset J_F$ on pose 
$p(J)=\sum_{\tau\in J}(k_0\!-\!m_\tau\!-\!1)\tau+
\sum_{\tau\in J_F\bs J}m_\tau\tau$; de m{\^e}me, pour  
$a=\sum_{\tau\in J_F}a_\tau\tau\in\Z[J_F]$, on pose
$|a|=\sum_{\tau\in J_F}a_\tau\in \Z$. 
Le complexe BGG est d{\'e}fini comme :
$$\overline{\mathcal{K}}_n^i=\bigoplus_{J \subset J_F, |J|=i} 
\overline{\mathcal{W}}_{\epsilon_J(n+t)-t,n_0}.$$

Soit $H=\Big(\Big(\begin{array}{cc}0&0\\0&-1\end{array}\Big)_{\tau}\Big)_{\tau\in J_F}\in(\mathfrak{gl}_2)^{J_F}$. 
On a $-(\epsilon_J(n+t)-t,n_0)(H)=|p(J)|$ ; en effet 
le caract{\`e}re de $T$ correspondant {\`a} $(n;n_0)$ est donn{\'e} par la formule 
$\Big(\begin{array}{cc}a & 0 \\ 0 & d\end{array}\Big)
\mapsto a^{(n_0t+n)/2}d^{(n_0t-n)/2}$. Ainsi pour tout $\tau\in J_F$ on a 
$$-(\epsilon_\tau(n_\tau+1)-1;n_0)(H)=\begin{cases}
\frac{n_0-n_\tau}{2}=m_\tau\text{, si }\epsilon_\tau=1 
(\iff \tau\in \overline{J}),\\
\frac{n_0+n_\tau+2}{2}=k_0-m_\tau-1 \text{, si }\epsilon_\tau=-1
(\iff \tau\in J).\end{cases}$$ 

La $F$-filtration sur $\overline{\mathcal{K}}_n^\bullet$ est donn{\'e}e par 
$\Fil^i\overline{\mathcal{K}}_n^\bullet=\bigoplus_{J \subset J_F,
 |p(J)|\geq i} 
\overline{\mathcal{W}}_{\epsilon_J(n+t)-t,n_0}.$

\begin{theo}\label{hodge-pad} 

{\rm (i)} On a un quasi-isomorphisme de complexes filtr{\'e}s
$$\overline{\mathcal{K}}_n^\bullet\hookrightarrow \overline{\mathcal{V}}_n\otimes 
\Omega^\bullet_{\overline{M^1}}(\dlog\,\infty).$$
 
{\rm (ii)} La suite spectrale donn{\'e}e par la  $F$-filtration 
$$\mathrm{E}_1^{i,j}=\bigoplus_{J \subset J_F, |p(J)|=i} 
\mathrm{H}^{i+j-|J|} (\overline{M^1}, 
\overline{\mathcal{W}}_{\epsilon_J(n+t)-t,n_0})
\Rightarrow
\mathrm{H}^{i+j}(\overline{M^1},\overline{\mathcal{V}}_n\otimes 
\Omega^\bullet_{\overline{M^1}}(\dlog\,\infty))$$
d{\'e}g{\'e}n{\`e}re en $\mathrm{E}_1$.

{\rm (iii)} Pour tout entier $j$, $0\leq j\leq d$, les poids 
de Hodge-Tate de la repr{\'e}sentation $p$-adique 
$\mathrm{H}^j(M^1\otimes \overline{\Q}_p,\mathbb{V}_n)$ 
appartiennent {\`a} l'ensemble $\{|p(J)|\enspace , \enspace |J|\leq j\}$.
\end{theo}

Ce th{\'e}or{\`e}me admet un corollaire, donnant des propri{\'e}t{\'e}s $p$-adiques
des repr{\'e}sentations galoisiennes associ{\'e}es aux formes modulaires de 
Hilbert. Prenons  $G=\Res^F_{\Q}\GL_2$ de sorte qu'on connaisse l'existence 
de ces repr{\'e}sentations galoisiennes.  
Soit $f\in G_{\kappa}(\mathfrak{c},\mathfrak{n})$ 
une  forme de Hilbert cuspidale pour  $\Res^F_{\Q} \GL_2$  propre  pour tous 
les op{\'e}rateurs  de Hecke, primitive de poids alg{\'e}brique $\kappa$ 
(voir \ref{alg}).
Soit $\rho_f$ la repr{\'e}sentation de $\Gal(\overline{\Q}/F)$ dans 
$\GL_2(\overline{\Q}_p)$, associ{\'e}e {\`a} $f$ et {\`a} un plongement de 
$\overline{\Q}$ dans $\overline{\Q}_p$. Soit $W_f$ la $f$-partie  de 
$\mathrm{H}^d_!(M_{\overline{\Q}},\mathbb{V}_{n}(\overline{\Q}_p))$.
Soit $\widetilde{F}$ la cl{\^o}ture galoisienne de $F$ dans $\overline{\Q}$.
Pour tout $\tau\in J_F$ on note $f_\tau$ le conjugu{\'e} interne de $f$ 
par $\tau$.  D'apr{\`e}s un r{\'e}sultat 
de Brylinski et Labesse \cite{BL} les semi-simplifications des
restrictions {\`a} $\Gal(\overline{\Q}/\widetilde{F})$ des repr{\'e}sentations 
$W_f$ et  
$\underset{\tau\in J_F}{\otimes}\rho_{f_\tau}$ sont isomorphes.

En prenant les invariants de la $F$-filtration de 
$\overline{\mathcal{V}}_n \otimes 
\Omega^\bullet_{\overline{M^1}}(\dlog\,\infty)$ par le groupe
 de Galois du rev{\^e}tement {\'e}tale 
$\overline{M^1}\rightarrow \overline{M}$,
on obtient une filtration sur le complexe $\overline{\mathcal{V}}_n \otimes 
\Omega^\bullet_{\overline{M}}(\dlog\,\infty)$ sur $\overline{M}$, 
appel{\'e}e encore $F$-filtration.
De m{\^e}me, on d{\'e}finit le complexe BGG sur $\overline{M}$ en prenant 
les invariants du complexe BGG
sur $ \overline{M^1}$. La suite spectrale associ{\'e}e est donn{\'e}e par 
invariants de la suite spectrale 
du Th{\'e}or{\`e}me \ref{hodge-pad} (ii). D'o{\`u} la premi{\`e}re assertion du 
 
\begin{cor}
{\rm (i)} La suite spectrale donn{\'e}e par la  $F$-filtration 
$$\mathrm{E}_1^{i,j}=\bigoplus_{J \subset J_F, |p(J)|=i} 
\mathrm{H}^{i+j-|J|} (\overline{M}, 
\overline{\mathcal{W}}_{\epsilon_J(n+t)-t,n_0})
\Rightarrow
\mathrm{H}^{i+j}(\overline{M},\overline{\mathcal{V}}_n\otimes 
\Omega^\bullet_{\overline{M}}(\dlog\,\infty))$$
d{\'e}g{\'e}n{\`e}re en $\mathrm{E}_1$.

{\rm (ii)}
 Les poids de Hodge-Tate de $W_f$ sont les entiers 
$\vert p(J)\vert$, $J\subset J_F$, compt{\'e}s avec multiplicit{\'e}.
\end{cor}

\noindent{\bf D{\'e}monstration :}  (ii)
On a   $\overline{\mathcal{W}}_{\epsilon_J(n+t)-t,n_0}=
\underline{\omega}^{-\epsilon_J(n+t)+t}\otimes \underline{\Det}^{p(J)}$. 

 Il r{\'e}sulte du th{\'e}or{\`e}me \ref{hodge-pad} (comme dans \cite{FaCh} Th.5.5
 et \cite{MoTi} Sect.2.3) que les sauts de la filtration de 
Hodge figurent parmi les $|p(J)|$, $J\subset J_F$. De plus, 
$$ \gr^{|p(J)|}\mathrm{H}^{d}(\overline{M},\overline{\mathcal{V}}_n\otimes 
\Omega^\bullet_{\overline{M}}(\dlog\,\infty))= \mathrm{H}^{d-|J|}(\overline{M},
\underline{\omega}^{-\epsilon_J(n+t)+t}\otimes \underline{\Det}^{p(J)}).$$

D'apr{\`e}s \cite{hmv-coh} Cor.2.7 la suite spectrale est Hecke {\'e}quivariante.
Il suffit donc de voir que pour tout $J\subset J_F$, le 
$\overline{\Q}_p$-espace vectoriel  
$$\left(\mathrm{H}^{d-|J|}(\overline{M},
\underline{\omega}^{-\epsilon_J(n+t)+t}\otimes \underline{\Det}^{p(J)})
\otimes \overline{\Q}_p\right)[f]$$ est de dimension $1$. 

Gr{\^a}ce {\`a} l'existence d'une structure $\overline{\Q}$-rationnelle 
du complexe BGG sous-jacent aux complexes BGG sur $\overline{\Q}_p$
et sur  $\overline{\C}$, en prenant un  plongement de $\Q_p$ dans $\C$, 
on a un  isomorphisme Hecke-{\'e}quivariant 
$$\mathrm{H}^{d-|J|}(\overline{M},
\underline{\omega}^{-\epsilon_J(n+t)+t}\otimes 
\underline{\Det}^{p(J)})  \otimes_{\Q_p} {\C}=
\mathrm{H}^{J,\overline{J}}(M^{\an }, \mathbb{V}_{n,\C}).$$

Or, par l'isomorphisme d'Eichler-Shimura-Harder  la
$f$-partie  $\mathrm{H}^{J,\overline{J}}(M^{\an }, \mathbb{V}_{n,\C})[f]$
est de dimension $1$ sur  $\C$, pour tout $J\subset J_F$ (voir \cite{hida74}).
\hfill$\square$

\begin{rque}
1) Le motif $W_f$ est pur de poids $(k_0-1)d$. 
L'ensemble de ses poids de Hodge-Tate est stable par la sym{\'e}trie 
$h \mapsto (k_0-1)d-h$, correspondant au passage au compl{\'e}mentaire
$|p(J)|\mapsto |p(\overline{J})|$. Cette sym{\'e}trie est induite par la 
dualit{\'e} de  Poincar{\'e} $W_f\times W_f\rightarrow \Q(-(k_0-1)d)$.

2) Si $F$ est un corps quadratique et $f$ une forme de Hilbert 
cuspidale propre de poids $\kappa$ sur $F$. 
En notant  $\tau$ le plongement 
non-trivial de $F$, on voit que   les sauts de la filtration de 
Hodge de $W_f$ sont $m_\tau, k_0-m_\tau-1, k_0+m_\tau-1, 2k_0-m_\tau-2$.
\end{rque}

\bigskip
\noindent {\bf Sur $\Z_p$ :}
On a {\'e}galement une version cristalline.
Soit $p$ un nombre premier ne divisant pas $\Delta$ et 
soit  $V$ un $\Z_p$-module libre de type fini muni
d'une action alg{\'e}brique de $G$ de plus haut poids $n\in\n[J_F]$. 
 On suppose que  $p-1>\sum_\tau(n_\tau+1)$. Comme dans le 
paragraphe pr{\'e}c{\'e}dent, on peut associer {\`a} $V$

1) un faisceau lisse $\mathbb{V}$ sur $M^1\otimes \Q_p$, et 

2) un fibr{\'e} $\overline{\mathcal{V}}$  sur $\overline{M^1}\otimes \Z_p$, 
filtr{\'e} et {\`a} connexion logarithmique int{\'e}grable quasi-nilpotente. 

\medskip
{\'E}tant donn{\'e}e des compactifications toro{\"\i}dales 
$\overline{f}:\overline{\mathcal{A}}\rightarrow \overline{M^1}$,
on a un un  th{\'e}or{\`e}me de comparaison de Faltings \cite{Fa-jami}
reliant $\mathrm{H}^\bullet(M^1\otimes\overline{\Q}_p,\mathbb{V})$
et $\mathrm{H}^\bullet_{\mathrm{log-cris}}
(\overline{M^1},  \overline{\mathcal{V}})$.

\medskip
Prenons pour $V$ la repr{\'e}sentation g{\'e}om{\'e}triquement irr{\'e}ductible 
de plus haut poids $n\geq 0$.  Un analogue  sur $\F_p$ de la suite spectrale
du th{\'e}or{\`e}me  \ref{hodge-pad} donne le th{\'e}or{\`e}me (voir  \cite{MoTi} et 
\cite{hmv-coh} Sect.5)

\medskip
\begin{theo} Supposons que $p$ ne divise pas $\Delta$ et que
$p-1>\sum_\tau(n_\tau+1)$. Alors pour tout $j$ 
compris entre $0$ et $d$, les poids de Hodge-Tate 
du module cristallin 
$\mathrm{H}^j(M\otimes \overline{\Q}_p, \mathbb{V}_n)$
sont les $|p(J)| $, o{\`u} $J$ parcourt l'ensemble 
des parties de $J_F$ telles que $|J|\leq j$.
\end{theo}

\section[Autres formes  de Hilbert  arithm{\'e}tiques.]{Autres formes  de Hilbert  arithm{\'e}tiques}

Dans toute cette partie on suppose que $D=\Gm$, c'est-{\`a}-dire 
$G=G^*$, de sorte que tous les sous-groupes de congruence
consid{\'e}r{\'e}s sont contenus dans $\SL_2(F)$, et en 
particulier $\Gamma=\Gamma^1$, $M^{\an}=M^{1,\an}$ et $M=M^1$.

\subsubsection*{Formes de Hilbert-Jacobi.}

Le but de ce paragraphe est de poser les d{\'e}finitions et les
propri{\'e}t{\'e}s de base  des formes
modulaires de Hilbert-Jacobi arithm{\'e}tiques. Le cas des
formes de Siegel-Jacobi a d{\'e}j{\`a} {\'e}t{\'e} trait{\'e} par J. Kramer \cite{kramer}.
Nous donnons  d'abord la d{\'e}finition des formes modulaires de Hilbert-Jacobi 
sur $\C$, inspir{\'e}e de \cite{EiZa}.

Dans le paragraphe sur la VAHB analytique universelle nous avons 
d{\'e}j{\`a} consid{\'e}r{\'e}  le groupe produit semi-direct 
$\Gamma^J=(\mathfrak{o}\oplus\mathfrak{c}^*)\rtimes\Gamma$
(pour $\gamma\cdot(m,n)=(m,n)\gamma^{-1}$).
Ce groupe  agit {\`a} gauche sur $\mathfrak{H}_F\times(F\otimes\C)$ par :
$$\begin{cases} \gamma(z,v)=(\gamma(z),j(\gamma,z)^{-1}v) \\
(m,n)(z,v)=(z,v+m\otimes z+n\otimes 1)\end{cases}.$$

Soient $\kappa\in\Z[J_F]=X(T)$ et $\mu\in\mathfrak{c}=X(\Gm\otimes \mathfrak{c}^*)$
et soit $e^{\mu}=\mu\circ q: F\otimes\C\rightarrow \C^\times$ 
la compos{\'e}e de l'application 
$q:F\otimes\C\rightarrow \Gm\otimes \mathfrak{c}^*$, d{\'e}finie dans (\ref{qexp}),
et du caract{\`e}re  $\mu$. Pour chaque    {\'e}l{\'e}ment $\gamma=
\begin{pmatrix} a & b \\ c & d\end{pmatrix}\in \Gamma$  et 
$(m,n)\in \mathfrak{o}\oplus\mathfrak{c}^*$ on d{\'e}finit une transformation 
 lin{\'e}aire de  l'espace des fonctions holomorphes 
$$f : \mathfrak{H}_F\times(F\otimes\C) \rightarrow \C\text{ , }
(z,v)\mapsto f(z,v),$$  en posant  :
\begin{equation}\label{jacobi}
\begin{cases} (f|_{\kappa,\mu}\gamma)(z,v)=
j(\gamma,z)^{-\kappa}e^{\mu}(-\frac{cv^2}{j(\gamma,z)})f(\gamma(z,v)) \\
(f|_{\kappa,\mu}(m,n))(z,v)=e^{\mu}(m^2z+2mv)f((m,n)(z,v))
\end{cases}
\end{equation}

 Les relations suivantes :

(i) $(f|\gamma)|\gamma'=f|(\gamma\gamma')$, pour tout $\gamma,\gamma'\in \Gamma$,

(ii) $(f|(m,n))|(m',n')=f|(m+m',n+n')$, pour tout $m,m'\in \mathfrak{o}$
et $n,n'\in\mathfrak{c}^*$,

(iii) $(f|(m,n))|\gamma=(f|\gamma)|((m,n)\gamma)$, pour tout $\gamma\in \Gamma$,
$m \in \mathfrak{o}$ et $n\in\mathfrak{c}^*$,

\noindent sont  faciles {\`a} d{\'e}montrer 
(le calcul r{\'e}v{\`e}le que (ii) et (iii) sont  {\'e}quivalentes 
respectivement {\`a} $e^{\mu}(2mn)=1$ et $e^{\mu}(cdn^2+abm^2+2bcmn)=1$).

Une fa{\c c}on {\'e}quivalente de formuler (i), (ii) et (iii) {\`a} la fois,  est de dire que
$(\ref{jacobi})$ d{\'e}finit une  action  du groupe produit semi-direct 
$\Gamma^J$ sur les  fonctions holomorphes sur  $\mathfrak{H}_F\times(F\otimes\C)$.

\begin{defin}
Une forme modulaire de Hilbert-Jacobi de poids $\kappa\in\Z[J_F]$, indice 
$\mu\in\mathfrak{c}$ et niveau $\Gamma$ est une fonction holomorphe 
$f : \mathfrak{H}_F\times(F\otimes\C) \rightarrow \C$ v{\'e}rifiant:

(i) $f|_{\kappa,\mu}\gamma=f$, pour tout $\gamma\in \Gamma$,

(ii) $f|_{\kappa,\mu}(m,n)=f$, pour tout $(m,n)\in\mathfrak{o}\oplus\mathfrak{c}^*$;

(iii) $f$ est ``holomorphe {\`a} l'infini'' : pour chaque pointe 
$\mathcal{C}=\gamma\infty\in
\mathbb{P}^1(F)$, avec $\gamma\in G_{\Q}$, 
la fonction $f_{\mathcal{C}}:=f|_{\kappa,\mu}\gamma$ 
admet  un d{\'e}veloppement en s{\'e}rie de Fourier
$$f_{\mathcal{C}}(z,v)=\sum_{\xi\in X, \alpha\in \mathfrak{a}} a_{\xi,\alpha} 
e^{2i\pi \Tr_{F/\mathbb{Q}}(\xi z+ \alpha v)},$$
similaire {\`a} celui du (\ref{dev-fourier}).
 La condition d'holomorphie en la pointe $\mathcal{C}$ se lit alors :
\begin{equation}\label{holo-jacobi}
a_{\xi,\alpha} \neq 0\Rightarrow 4\xi\mu-\alpha^2\in(X\mathfrak{c})_+\cup\{0\}.
\end{equation}

\end{defin}

\noindent{\bf Principe du $q$-d{\'e}veloppement : } Si
pour tout $\xi\in X$, $\alpha\in \mathfrak{a}$ on a 
$a_{\xi,\alpha}=0$, alors $f=0$.

\medskip

Pour tout  $u\in \mathfrak{o}_{\mathcal{C}}^\times$, il
existe $\xi_u^*\in (\mathfrak{ab})^*$, d{\'e}fini {\`a} $X^*$ pr{\`e}s, tel que
$\begin{pmatrix} u & \xi_u^* \\ 0 & u^{-1}\end{pmatrix}\in
\gamma^{-1}\Gamma\gamma\cap B_{\R}$. L'invariance de 
$f_{\mathcal{C}}$ par le groupe $\gamma^{-1}\Gamma\gamma\cap B_{\R}$
nous donne pour tout $\xi\in X$ la relation :
\begin{equation}\label{units-jacobi}
a_{u^2\xi,u\alpha}=u^\kappa  e^{2i\pi\Tr_{F/\mathbb{Q}}(u\xi\xi_u^*)}a_{\xi,\alpha}.
\end{equation}

En utilisant le diagramme commutatif suivant :
$$\xymatrix@R=15pt@C=30pt{0\ar[r]&\mathfrak{c}^*\ar^{\phi}[r]&
\C\otimes\mathfrak{c}^*
\ar^{e^{2i\pi\cdot}\otimes \id}[r]\ar_{\id\otimes \Tr(\mu\cdot)}[d]
 &\C^\times\otimes\mathfrak{c}^*\ar[r]\ar^{\mu}[d] & 0,\\
& &\C\ar^{e^{2i\pi\cdot}}[r] &\C^\times & }$$

on obtient pour tout $\beta\in \mathfrak{b}$ la relation :
\begin{equation}\label{relation-jacobi}
a_{\mu\beta^2+\alpha\beta+\xi,\alpha+2\mu\beta}=a_{\xi,\alpha}.
\end{equation}

\noindent{\bf Principe de Koecher : } Si $F\neq \mathbb{Q}$, alors 
la condition (\ref{holo-jacobi}) est toujours satisfaite.
Si $\kappa$ n'est pas parall{\`e}le, alors $a_{0,0}=0$ (pas de 
s{\'e}ries d'Eisenstein parmi les formes de Hilbert-Jacobi).
Enfin, il n'existe de formes de Hilbert-Jacobi
non-nulles que si $\mu\in\mathfrak{c}_+$.

\smallskip

On montre d'abord, par l'absurde, que si $a_{\xi,\alpha}$
est non nul, alors $\xi \in X_+$.
Soient $\alpha\in\mathfrak{a}$,
$\xi\in X$ et $\xi^* \in X_+^*$ 
tels que $a_{\xi,\alpha}\neq 0$ et  $\langle\xi,\xi^*\rangle<0$. 
D'apr{\`e}s (\ref{units-jacobi}),
pour tout $u\in \gamma^{-1}\Gamma\gamma\cap T_{\R}$, on 
a $a_{u^2\xi,u\alpha}=u^\kappa  a_{\xi,\alpha}\neq 0$.
En particulier, $a_{0,0}=u^\kappa  a_{0,0}$,
d'o{\`u} la deuxi{\`e}me propri{\'e}t{\'e}.
Comme $f$  est holomorphe au point $(\underline{i}\xi^*,0)$, 
la s{\'e}rie $\sum_{u\in \gamma^{-1}\Gamma\gamma\cap T_{\R}}
u^\kappa  e^{-2i\pi\langle u^2\xi,\xi^*\rangle}$ converge absolument, 
ce qui est impossible, par le th{\'e}or{\`e}me des unit{\'e}s de Dirichlet.

La relation (\ref{relation-jacobi}), nous dit alors que 
pour tout $\beta\in\mathfrak{b}$, on a 
$\mu\beta^2+\alpha\beta+\xi\in X_+$. On en d{\'e}duit que
$\mu\in\mathfrak{c}_+$ et $4\xi\mu-\alpha^2\in(X\mathfrak{c})_+\cup\{0\}$.
\hfill $\square$

\medskip
Tout comme les formes modulaires de Hilbert, les formes modulaires
 de Hilbert-Jacobi admettent elles aussi une d{\'e}finition
purement g{\'e}om{\'e}trique. 

Tout {\'e}l{\'e}ment $\mu\in\mathfrak{c}$ donne un morphisme 
$(\id,\lambda\circ (\id\otimes \mu)):\mathcal{A}\rightarrow \mathcal{A}\times \mathcal{A}^t$,
d'o{\`u} un faisceau inversible ample sur $A$, 
$\mathcal{L}_{\mu}=(\id,\lambda\circ 
(\id\otimes \mu))^*\mathcal{P}_{\mathcal{A}}$, o{\`u} $\mathcal{P}_{\mathcal{A}}$ 
d{\'e}signe le faisceau de Poincar{\'e} sur $\mathcal{A}\times \mathcal{A}^t$.

\begin{defin}
Soit $R$ une $\Z[\frac{1}{\Delta}]$-alg{\`e}bre contenant les valeurs de $\kappa$.
Une forme modulaire de Hilbert-Jacobi de poids $\kappa\in\Z[J_F]$, indice 
$\mu\in\mathfrak{c}$, niveau $\Gamma$ et {\`a} coefficients dans $R$, est 
une section globale de $f^*\underline{\omega}^\kappa\otimes\mathcal{L}_{\mu}$ sur
$\mathcal{A}\times_{\Spec(\Z[\frac{1}{\Delta}])} \Spec(R)$.
On note 
$J_{\kappa,\mu}(R)=J_{\kappa,\mu}(\mathfrak{c},\mathfrak{n};R):=
\mathrm{H}^0(\mathcal{A}\times_{\Spec(\Z[\frac{1}{\Delta}])} \Spec(R)
,f^*\underline{\omega}^{\kappa}\otimes\mathcal{L}_{\mu})$ l'espace  
de ces formes modulaire de Hilbert-Jacobi. 
\end{defin}

Soit une $(R,\mathfrak{n})$-pointe $\mathcal{C}$ et soit 
$R$ une $\Z[\frac{1}{\Delta},\zeta_{\mathcal{C}}]$-alg{\`e}bre.
En {\'e}valuant  une forme de Hilbert-Jacobi $f\in 
J_{\kappa,\mu}(R)$
sur les objets de Tate associ{\'e}s  {\`a}  une  pointe $\mathcal{C}$,
on obtient comme dans \cite{kramer} le $q$-d{\'e}veloppement de $f$ en $\mathcal{C}$
$$f_\mathcal{C}=\sum_{\xi\in X, \alpha\in\mathfrak{a}}
a_{\xi,\alpha}q^\xi\mathfrak{X}^\alpha.$$

En transposant la m{\'e}thode de \cite{kramer} au cas de Hilbert on obtient
alors les relations (\ref{units-jacobi}) et (\ref{relation-jacobi})
qui impliquent, comme plus haut le principe de Koecher (\ref{holo-jacobi}),
lorsque $F\neq\Q$.

{\bf Principe du $q$-d{\'e}veloppement : } Soient  $M_1\subset M_2$ sont des groupes
ab{\'e}liens et $f\in J_{\kappa,\mu}(M_2)$. Si le  $q$-d{\'e}veloppement
de $f$ en une $(R,\mathfrak{n})$-pointe  $\mathcal{C}$ est 
{\`a} coefficients dans $M_1$, alors $f\in J_{\kappa,\mu}(M_1)$. 

\medskip

{\'E}tant donn{\'e} une compactification $\overline{\mathcal{A}}=
\mathcal{A}_{\widetilde{\Sigma},\varphi}$ de la 
vari{\'e}t{\'e} de Hilbert-Blumenthal universelle $\mathcal{A}$, comme dans la partie
\ref{kugasato}, on d{\'e}finit l'espace  
des formes modulaire de Hilbert-Jacobi relatives {\`a} cette compactification :
$$J_{\kappa,\mu}(R;\widetilde{\Sigma},\varphi):=
\mathrm{H}^0(\overline{\mathcal{A}}\times_{\Spec(\Z[\frac{1}{\Delta}])} \Spec(R)
,f^*\underline{\omega}^{\kappa}\otimes\mathcal{L}_{\mu}).$$

{\`A} noter que le prolongement {\`a} $\overline{\mathcal{A}}$ du faisceau 
inversible ample $\mathcal{L}_{\mu}$ d{\'e}pend de la polarisation $\varphi$.

Contrairement au cas des formes modulaires de Hilbert, pour les formes 
de Hilbert-Jacobi on a juste une inclusion 
$$J_{\kappa,\mu}(R;\widetilde{\Sigma},\varphi)\hookrightarrow 
J_{\kappa,\mu}(R)$$

Dans \cite{kramer} Kramer d{\'e}montre que cette inclusion est stricte pour
les formes de Siegel-Jacobi. Il serait int{\'e}ressant d'{\'e}tudier cette
 question dans le cas de formes de  Hilbert-Jacobi.

Comme dans la partie pr{\'e}c{\'e}dente, 
on peut alors associer {\`a} toute $(R,\mathfrak{n})$-pointe  
$\mathcal{C}$ et tout $f\in J_{\kappa,\mu}(R;\widetilde{\Sigma},\varphi)$ 
son $q$-d{\'e}veloppement en $\mathcal{C}$
$$f_{\mathcal{C}}=\underset{(\xi,\alpha)\geq \varphi}{\sum_{\xi\in X, \alpha\in\mathfrak{a}}}
a_{\xi,\alpha}q^\xi\mathfrak{X}^\alpha.$$

Notons que 
$$(\xi,\alpha)\geq \varphi \iff \forall q\in X_+^*,\, \forall l \in \mathfrak{a}^*
\max_{\beta\in \mathfrak{b}} \Tr_{F\!/\!\Q}(-\mu q\beta^2-2\mu l \beta +\xi q+\alpha l)\geq 0,$$

 alors que
$$4\xi\mu-\alpha^2 \geq 0  \iff \forall q\in X_+^* ,\, \forall l \in \mathfrak{a}^*
\max_{\beta\in \mathfrak{b}_{\R}} \Tr_{F\!/\!\Q}(-\mu q\beta^2-2\mu l \beta +\xi q+\alpha l)\geq 0.$$

L'anneau o{\`u} vit le $q$-d{\'e}veloppement d'une forme de 
$J_{\kappa,\mu}(R;\widetilde{\Sigma},\varphi)$ est donc en g{\'e}n{\'e}ral strictement inclus
dans celui o{\`u} vit le $q$-d{\'e}veloppement d'une forme de 
$J_{\kappa,\mu}(R)$. 
Il serait int{\'e}ressant de savoir s'il existe des {\'e}ventails 
$\widetilde{\Sigma}$, munis d'une fonction de polarisation $\varphi$
pour lesquelles on a une {\'e}galit{\'e}.

\subsubsection*{S{\'e}ries th{\^e}ta et formes de Hilbert de poids demi-entier.}
${}$

\label{theta}R{\'e}f{\'e}rences:\cite{sh-half1},\cite{sh-half2},\cite{wu}

Dans cette section, on suppose que $\mathfrak{c}=\mathfrak{c}_0^2$ 
est un carr{\'e} et est premier {\`a} $2$.

\noindent{\bf Sur $\C$:}   
 On d{\'e}duit facilement de  \cite{sh-half1}
(Prop.1.1 et 1.2, en fait plut{\^o}t Prop.3.2 et Lemme 3.5) 
qu'il existe un facteur d'automorphie de poids $t/2$ pour 
$\Gamma_0(\mathfrak{c},4)$,
c'est-{\`a}-dire
une fonction  $h:\Gamma_0^1(\mathfrak{c},4)\times\mathfrak{H}_F\rightarrow \C^\times$
holomorphe en la seconde variable  telle que 
$$h(\gamma_1\gamma_2,z)=h(\gamma_1,\gamma_2(z))\cdot h(\gamma_2,z)$$

et pour tout $\gamma\in \Gamma_0^1(\mathfrak{c},4)$, en notant 
$t=\sum_{\tau\in J_F} \tau$,
$$(*)\quad h(\gamma,z)^2=\chi(\gamma)\cdot j(\gamma,z)^t$$
 o{\`u} $\chi$ est un caract{\`e}re quadratique de $\Gamma_0(4)$ qui ne d{\'e}pend que de l'image dans
$(\mathfrak{o}/4\mathfrak{o})^\times$ de $d_\gamma$.

Soit $\Gamma_0^1(\mathfrak{c},4)_+=\Ker(\chi)$. C'est un sous-groupe de congruences d'indice $2$ de 
$\Gamma_0^1(\mathfrak{c},4)$.
Par exemple, si $F=\Q$, on a $\Gamma_0(4)_+=\Gamma_1(4)$.

Soit $\mathfrak{n}$ un id{\'e}al entier de $\mathfrak{o}$ tel que $\Gamma$ soit sans torsion; 

\noindent{\bf Hypoth{\`e}se:} on suppose dans toute cette partie ainsi que dans celle
concernant les formes de Hilbert $p$-adiques de poids demi-entier que $4$ divise $\mathfrak{n}$ 
dans $\mathfrak{o}$. 

On a alors
$$\Gamma=\Gamma_1^1(\mathfrak{c},\mathfrak{n})\subset 
\Gamma_0^1(\mathfrak{c}, 4)_+.$$

Par la formule $(*)$ ci-dessus, on a donc pour tout $\gamma\in \Gamma$, 
$h(\gamma,z)^2=j(\gamma,z)^t$.
 
Pour tout poids demi-entier $\kappa={\frac{t}{2}}+\lambda$, $\lambda\in\Z[J_F]$, on pose
pour tout $\gamma\in \Gamma$, 
$j_\kappa(\gamma,z)=h(\gamma,z)\cdot j(\gamma,z)^\lambda$.

On note $G_\kappa(\Gamma)$ l'espace des fonctions
$f$
holomorphes sur $\mathfrak{H}_F$
satisfaisant $f(\gamma(z))=j_\kappa(\gamma,z)\cdot f(z)$ pour tout 
$\gamma\in \Gamma$.

On d{\'e}finit un fibr{\'e} inversible $\underline{\omega}^\kappa$ holomorphe sur $M^{\an}$, correspondant
 au facteur
d'automorphie $j_\kappa(\gamma,z)$; c'est le quotient de $\mathfrak{H}_F\otimes \C$ par 
$\Gamma$
agissant par $\gamma(z,u)=(\gamma(z),j_\kappa(\gamma,z)\cdot u)$.

On a un isomorphisme canonique
$G_\kappa(\Gamma)\cong \mathrm{H}^0(M^{\an},\underline{\omega}^\kappa)$.

Pour  $z\in \C$, soit $e(z)=\exp(2i\pi z)$ et 
${\bf e}:(F\otimes\C)\rightarrow \C^\times, z=(z_\tau)_\tau\mapsto e(\sum_\tau z_\tau)$.

Soit $\mathfrak{n}=4\mathfrak{n}_0$. Pour chaque fonction  
 $\eta:\mathfrak{c}_0 \rightarrow \C$,  constante modulo 
$\mathfrak{n}_0\mathfrak{c}_0$, on d{\'e}finit une s{\'e}rie th{\^e}ta par
$$\theta(z,\eta)=\sum_{\alpha\in 
 \mathfrak{c}_0}\eta(\alpha){\bf e}(\alpha^2z)$$

C'est la valeur en $v=0$ (Thetanullwert) de la fonction th{\^e}ta de 
$(z,v)\in \mathfrak{H}_F\times (F\times \C)$ 
et $\eta\in\mathcal{S}(F_f)$ d{\'e}finie comme suit : on fixe un {\'e}l{\'e}ment $c_0$ de $\mathfrak{c}_0$ 
premier {\`a} $2$ et qui engendre 
$\mathfrak{c}_0/\mathfrak{nc}_0$ et on pose
$$\theta(z,v;\eta)=\sum_{\alpha\in \mathfrak{c}_0}
\eta(\alpha){\bf e}\left(\alpha^2z+c_0\alpha v\right)$$

\begin{lemme} On a $\theta(z,\eta)\in G_{t/2}(\Gamma)$. Les fonctions
 $\theta(z,\eta)$ pour $\eta$ parcourant l'ensemble des 
fonctions sur $\mathfrak{c}_0/\mathfrak{n}_0\mathfrak{c}_0$,
d{\'e}finissent des sections 
globales de $\underline{\omega}^{t/2}$ sur $M^{\an}$.
\end{lemme}

\noindent{\bf D{\'e}monstration : }
La modularit{\'e} se d{\'e}duit de \cite{sh-half1} Prop.1.2, ou de \cite{sh-half2} (4.3).
Le point de la v{\'e}rification est que, au sens de la Prop.2.4 de \cite{sh-half2}, 
on a  ${}^g\eta=\eta$ pour tout
$g$ dans (le rel{\`e}vement de) l'adh{\'e}rence de $\Gamma$ dans $G_f$. 
Il suffit pour cela de voir que pour tout $x\in\mathfrak{c}_0\widehat{\mathfrak{o}}$,
la fonction caract{\'e}ristique $\eta_x$ de $x+\mathfrak{nc}_0\widehat{\mathfrak{o}}$ satisfait
${}^g\eta_x=\eta_x$. On se ram{\`e}ne {\`a} $g$ triangulaire sup{\'e}rieure ou bien triangulaire
inf{\'e}rieure, mais dans ce cas on exige que le coefficient $c$ engendre l'id{\'e}al 
$\mathfrak{c}_0^2\mathfrak{dn}$. On conclut alors {\`a} l'aide des formules (3.2) 
et (3.3) de \cite{sh-half1} et de la formule (i) de Prop.2.3 de \cite{sh-half2}.

\vskip 5mm

\noindent{\bf Sur $\Z[\frac{1}{2\Delta}]$ : } On va examiner la question de l'alg{\'e}brisation et du prolongement des fibr{\'e}s inversibles 
$\underline{\omega}^\kappa$ ($\kappa$ demi-entier) {\`a} une 
compactification toro{\"\i}dale $\overline{M}$. Pour cela on va rappeler des r{\'e}sultats classiques
dans le cas du groupe symplectique mais qui m{\'e}riteraient peut-{\^e}tre d'{\^e}tre d{\'e}taill{\'e}s davantage
dans le cas des vari{\'e}t{\'e}s de Hilbert.

Soit $N=N_1(\mathfrak{c},\mathfrak{n})$ l'espace de modules des $(A,\mathcal{L})$ o{\`u} 
$A$ est une VAHB $\mathfrak{c}$-polaris{\'e}e avec structure de niveau 
$\Gamma$ et o{\`u} $\mathcal{L}$ est un fibr{\'e}
inversible ample sym{\'e}trique, trivialis{\'e} le long de la section nulle,
d{\'e}finissant la polarisation sur $A$. 
Cet espace existe et est un sch{\'e}ma lisse sur $\Z[\frac{1}{2\Delta}]$ par \cite{FaCh} IV.7.
 Soit $f:\mathcal{A}\rightarrow M$ la vari{\'e}t{\'e} ab{\'e}lienne universelle sur $M$;
$N$ est un $\mathcal{A}[2]$-torseur sur $M$; il n'est pas connexe; on peut montrer que
son groupe des composantes connexes est naturellement isomorphe {\`a}
$\mathfrak{d}^{-1}/2\mathfrak{d}^{-1}$.
Soit $f_{ N}:\mathcal{A}_{ N}\rightarrow N$ le pull-back de $f$ {\`a} $N$.
Par d{\'e}finition, $\mathcal{A}_{N}$ est muni d'un fibr{\'e} inversible sym{\'e}trique
relativement ample universel $\mathcal{L}$.

Donnons une description  transcendante de $N$ 
(inspir{\'e}e de \cite{FaCh} V.3, p.161):
$$N^{\an}\cong \Gamma\bs \left(\mathfrak{H}_F\times 
\frac{1}{2}(\mathfrak{o}\oplus \mathfrak{c}^*)/
(\mathfrak{o}\oplus \mathfrak{c}^*)\right) $$
pour l'action $\gamma\cdot (z,x)=(\gamma(z),x 
\gamma^{-1}+(\frac{1}{2} cd,\frac{1}{2} ab))$,
pour $\gamma=\left(\begin{array}{cc}a&b\\c&d\end{array}\right)$.

\vskip 5mm
La construction de cet isomorphisme est un exercice utilisant les Thetanullwerte $\theta(z,z+c_0x,\eta_0)$,
 $x\in
\frac{1}{2}(\mathfrak{o}\oplus \mathfrak{c}^*)/(\mathfrak{o}\oplus \mathfrak{c}^*)$ et
o{\`u} $c_0\in\mathfrak{c}_0$; voir  p.160-161
de \cite{FaCh}.

\vskip 5mm

Le sous-groupe d'Igusa de $\Gamma_1^1(\mathfrak{c},\mathfrak{n}_0)$ donn{\'e} par les conditions $ab$ et $cd$ pairs
est le stabilisateur de $x=(0,0)$ dans 
$\Gamma$; ceci fournit par conjugaison par
la matrice $\left(\begin{array}{cc}2&0\\0&1\end{array}\right)$ une section $j^{\an}$
de  $M^{\an}$ dans $N^{\an}$. 
La description alg{\'e}brique de l'image de $M$ est le lieu o{\`u}
le fibr{\'e}  $\mathcal{L}$  est engendr{\'e} par ses sections paires. C'est une composante connexe de $N$.

\begin{rque}
Si $F=\Q$, seule une des quatre s{\'e}rie th{\^e}ta est paire et 
 seul le point $(\frac{1}{2}$, $\frac{1}{2})$,
parmi les  quatre points de $2$-torsion, est pr{\'e}serv{\'e} par le sous-groupe d'Igusa.
\end{rque}

Soit $j:M\rightarrow N$ la section alg{\'e}brique (sur $\Z[\frac{1}{2\Delta}]$)
ainsi d{\'e}finie.
Le fibr{\'e} inversible $\mathcal{L}$ vu comme $N$-sch{\'e}ma n'est pas un sch{\'e}ma en groupes,
mais il est muni d'une section nulle $0_{\mathcal{L}}$ (compos{\'e}e des sections nulles de 
$\mathcal{L}\rightarrow A$ et $A\rightarrow N$);
soit $\mathcal{L}_0=0_{\mathcal{L}}^*T_{\mathcal{L}/\mathcal{A}}$ le fibr{\'e} tangent relatif
de $\mathcal{L}\times_N M$
sur $M$.
On note pour abr{\'e}ger $\mathcal{V}=0_{\mathcal{L}}^* T_{\mathcal{L}/M}=
j^*0_{\mathcal{L}}^* T_{\mathcal{L}/N}$ le fibr{\'e} tangent de $\mathcal{L}\times_N M$ sur $M$.
Ce fibr{\'e}, muni de la filtration $\{0\}\subset \mathcal{L}_0\subset \mathcal{V}$, d{\'e}finit un torseur sous
le parabolique stabilisateur de $\mathcal{L}_0$ dont le quotient de Levi est $\mathcal{T}_1\times \Gm$.
 
Pour toute $\Z[\frac{1}{2\Delta}]$-alg{\`e}bre $R$ et tout caract{\`e}re $\xi$ du quotient de Levi
$\xi:\mathcal{T}_1\times \Gm
\rightarrow R^\times$, on d{\'e}finit un fibr{\'e} inversible alg{\'e}brique 
$\mathcal{V}_\xi=\mathcal{V}\stackrel{\xi}{\times}\mathbb{G}_a$  sur $M_R$ par contraction par $\xi$; 
c'est-{\`a}-dire par quotient par la relation d'{\'e}quivalence
$(v\theta,\alpha)\sim (v,\xi(\theta)\cdot \alpha)$ pour $\theta\in \mathcal{T}_1\times \Gm$, $v\in\mathcal{V}$
et $\alpha\in \mathbb{G}_a$.
\vskip 5mm

\begin{lemme} Pour tout poids demi-entier $\kappa={t/2}+\lambda$, soit $\xi_\kappa:
\mathcal{T}_1\times \Gm\rightarrow \C^\times$ donn{\'e} par $(x,y)\mapsto \lambda^{-1}(x)y$. On a un isomorphisme
canonique 
$$\underline{\omega}^{\kappa}\cong \mathcal{V}_{\xi_\kappa} $$
de fibr{\'e}s inversibles analytiques sur $M^{\an}$.
\end{lemme}

\noindent{\bf D{\'e}monstration : }
On a la description transcendante suivante de la vari{\'e}t{\'e} analytique 
$\mathcal{L}$ sur $M^{\an}$:
$$\mathcal{L}^{\an}\cong\Gamma\bs \left(\mathfrak{H}_F\times (F\otimes\C)
\times \C\right)/
(\mathfrak{o}\oplus \mathfrak{c}^*),$$
pour l'action $\gamma(z,v,u)=(\gamma (z), j(\gamma,z)^{-1}\cdot v,h(\gamma,z)
{\bf e}(\frac{c}{(cz+d)}v^2)\cdot u )$  

\noindent et  $(z,v,u)\cdot (m,n)=(z,v+m\cdot z+n,{\bf e}
(-\frac{1}{2}m^2v)\cdot u)$.

Pour voir que la formule donnant l'action de $\Gamma$ est correcte,
il suffit de donner le cocycle $\Gamma\times\mathfrak{H}_F\rightarrow (F\otimes\C)^\times\times \C^\times$
 qui d{\'e}finit le fibr{\'e} $0_{\mathcal{L}}^* T_{\mathcal{L}/M}^{\an}$. 

1) Sa $(F\otimes \C)^\times$-partie est impos{\'e}e: c'est $j(\gamma,z)^{-1}$;

2) pour la $\C^\times$-partie, il suffit de noter que les fonctions $v \mapsto \theta(z,v;\eta)$
 sont des sections
globales sur $(F\otimes \C)=\Lie(\mathcal{A}_z)$ 
de $\mathcal{L}_z$ par la th{\'e}orie analytique des diviseurs th{\^e}ta
(on laisse en exercice la d{\'e}monstration d{\'e}taill{\'e}e
de ce fait, disons seulement que l'on v{\'e}rifie que les invariants de Weil $(H,\psi,0,L)$ \cite{W} VI.3 de
$\theta(z,v,\eta)$ correspondent {\`a} ceux du
diviseur sym{\'e}trique ample $\Theta_z$ d{\'e}fini par $\mathcal{L}_z$). Le facteur d'automorphie cherch{\'e} est alors donn{\'e} par les relations
$\theta(\gamma(z,v);\eta)=h(\gamma,z)\cdot {\bf e}(\frac{c}{(cz+d)}v^2)\cdot\theta(z,v;\eta)$ pour 
$\gamma\in \Gamma$ (voir  Prop.1.1 de \cite{sh-half1}).

On en d{\'e}duit ais{\'e}ment que $\mathcal{V}_{\xi_\kappa}$ est d{\'e}fini par le m{\^e}me
facteur d'automorphie que $\underline{\omega}^{\kappa}$; d'o{\`u} le r{\'e}sultat.
L'argument des s{\'e}ries th{\^e}ta montre aussi que la formule pour l'action de $\mathfrak{o}\oplus 
\mathfrak{c}^*$ est {\'e}galement correcte.

\begin{cor} Pour tout poids demi-entier $\kappa=t/2+\lambda$, le fibr{\'e} inversible
complexe $\underline{\omega}^\kappa$ sur $M^{\an}$ est l'analytifi{\'e} du fibr{\'e} inversible alg{\'e}brique 
$\mathcal{V}_{\xi_\kappa}$,
d{\'e}fini sur $M$ au-dessus de $\Z[\frac{1}{2\Delta}]$. On notera encore $\underline{\omega}^\kappa$ ce fibr{\'e} 
alg{\'e}brique.
\end{cor}

\vskip 5mm

On construit, comme dans \cite{FaCh} IV.7, pour tout {\'e}ventail admissible $\Sigma$ une 
compactification toro{\"\i}dale lisse $\overline{\pi}:
\overline{N}\rightarrow \overline{M}$ compatible {\`a}  $N\rightarrow M$.  Rappelons que le sch{\'e}ma semi-ab{\'e}lien $\mathfrak{G}$ sur $\overline{M}$ a une action de $\mathfrak{o}$
et est lui-aussi muni d'une 
$\mathfrak{c}$-polarisation $\lambda: \mathfrak{c}\otimes_{\mathfrak{o}}\mathfrak{G}\cong \mathfrak{G}^t$
prolongeant celle de $\mathcal{A}$. Par d{\'e}finition
de $\overline{N}$, le pull-back $\mathfrak{G}_{N}$ de $\mathfrak{G}$ est muni d'un fibr{\'e} inversible sym{\'e}trique
relativement ample $\overline{\mathcal{L}}$ prolongeant $\mathcal{L}$.
Soit $0_{\overline{\mathcal{L}}}$ la section nulle de $\overline{\mathcal{L}}$ sur $\overline{N}$. 
On d{\'e}finit le prolongement de 
$\underline{\omega}^\kappa$ {\`a} $\overline{M}$ sur $\Z[\frac{1}{2\Delta}]$ comme $\overline{\mathcal{V}}_{\xi_\kappa}$
o{\`u} 
$\overline{\mathcal{V}}=\overline{j}^*0^*_{\overline{\mathcal{L}}}T_{\overline{\mathcal{L}}/\overline{N}}$.

Pour toute $\Z[\frac{1}{2\Delta}]$-alg{\`e}bre
$R$ contenant les valeurs de $\kappa$, on peut ainsi d{\'e}finir le $R$-module des formes arithm{\'e}tiques de poids demi-entier $\kappa$ par 
$$G\!D_\kappa(\Gamma;R)=\mathrm{H}^0(M\times \Spec(R),\underline{\omega}^\kappa)$$
et le principe de Koecher en poids demi-entier s'{\'e}nonce:
$$G\!D_\kappa(\Gamma;R)=\mathrm{H}^0(\overline{M}\times \Spec(R),\underline{\omega}^\kappa)
$$
Sa d{\'e}monstration est analogue au cas entier.

On peut m{\^e}me d{\'e}finir le module des formes de tout poids demi-entier. 
Notons  que $\underline{\omega}=f_*\Omega^1_{\mathcal{A}/M}=0^*\Omega^1_{\mathcal{A}/M}$ 
par propret{\'e} de $\mathcal{A}$ sur $M$ mais que
$0_{\mathcal{L}}^*\Omega^1_{\mathcal{L}/M}\neq f_{\mathcal{L},*}\Omega^1_{\mathcal{L}/M}$.

Consid{\'e}rons  le 
$\mathcal{T}_1$-torseur
$\mathfrak{M}=\Isom_{\mathcal{O}\otimes\mathcal{O}_{M}}
(\mathcal{O}\otimes\mathcal{O}_{M},\underline{\omega})$; 
formons un $\Gm$-torseur
sur $\mathfrak{M}$ donn{\'e} par
$$\mathfrak{M}^+=\Isom_{\mathfrak{M}}
(\mathcal{O}_{\mathfrak{M}}, 0^*_{\mathcal{L}}
\Omega^1_{\mathcal{L}/\mathfrak{M}}/0^*\Omega^1_{\mathcal{A}/\mathfrak{M}}))$$
ou encore $\mathfrak{M}^+=\Isom_{M}(\mathcal{O}_{M}, 
\mathcal{V}^\vee/\underline{\omega})$.

C'est un $\mathcal{T}_1\times\Gm$-torseur de Zariski sur $M$.

En tant que $M$-sch{\'e}ma, il classifie les syst{\`e}mes $(A,\lambda,\iota,\alpha,\mathcal{L},\omega,s)_S$ sur un 
$\Z[\frac{1}{2\Delta}]$-sch{\'e}ma $S$, o{\`u} $\omega$ est une 
$\mathcal{O}\otimes\mathcal{O}_{S}$-base de $0^*\Omega_{A/S}$, et $s$ une $\mathcal{O}_S$-base 
de $0^*\Omega_{\mathcal{L}/A}$.

On d{\'e}finit le module des formes modulaires de poids demi-entiers comme
$$G\!D(\Gamma,R)=\mathrm{H}^0(\mathfrak{M}^+\times 
\Spec(R),\mathcal{O}_{\mathfrak{M}^+})$$

C'est un module gradu{\'e} sur l'alg{\`e}bre $\Z[J_F]$-gradu{\'e}e 
$$G(\mathfrak{c},\mathfrak{n},R)=
\mathrm{H}^0(\mathfrak{M},\mathcal{O}_{\mathfrak{M}})$$
L'action $(f,g)\mapsto f\cdot g$ est induite par le pull-back {\`a} $\mathfrak{M}^+$ de 
$f\in \mathrm{H}^0(\mathfrak{M},\mathcal{O}_{\mathfrak{M}})$ par $\mathfrak{M}^+\rightarrow \mathfrak{M}$;
ce dernier morphisme provient de la projection de $\mathcal{L}$ 
sur $\mathcal{A}$.

\noindent{\bf Fait :}
Pour tout $\kappa={t/2}+\lambda$, on a
$\mathcal{O}_{\mathfrak{M}^+}[\xi_\kappa]=\mathcal{V}_{\xi_\kappa}=
\underline{\omega}^\kappa$ et donc 
donc  $G\!D_{\kappa}(\Gamma,R)=
G\!D(\Gamma,R)[\xi_\kappa]$.
En effet, une fonction $f(*,\omega,s)$ satisfaisant pour tout 
$$(x,y)\in\mathcal{T}_1\times\Gm,\quad f(*,x\omega,y\cdot s)=\lambda^{-1}(x)yf(*,\omega,s)$$
d{\'e}finit une section de $\mathcal{V}_{\xi_\kappa}$.

Par la d{\'e}finition alg{\'e}brique de $\underline{\omega}^\kappa$, on peut d{\'e}finir le $q$-d{\'e}veloppement
en une pointe $\mathcal{C}$ de $M$, 
d'une forme arithm{\'e}tique de poids demi-entier $\kappa$ comme
le compos{\'e} du d{\'e}veloppement de Fourier-Jacobi ({\it i.e.} 
le pull-back {\`a} la vari{\'e}t{\'e} de  Tate au-dessus
de $\overline{S}_{\Sigma^{\mathcal{C}}}^\wedge$, 
d{\'e}fini comme une s{\'e}rie gr{\^a}ce
{\`a} la trivialit{\'e} du pull-back du fibr{\'e} 
$e^*\Omega^1_{\overline{\mathcal{L}}/\G}$) 
et de l' {\'e}valuation $q_v\mapsto 1$  ({\it i.e.} $v=0$).  

Observons que sur le groupe des unit{\'e}s totalement positives de $F$, 
$u^\kappa  =u^{t/2}u^\lambda=u^\lambda$
est bien d{\'e}fini et est {\`a} valeurs dans l'anneau des valeurs de $\lambda$.

On a le ``{principe du $q$-d{\'e}veloppement}'' en une pointe $\mathcal{C}$ (disons
non-ramifi{\'e}e):

\begin{prop} Pour toute extension  $R\subset R'$
de $\Z[\frac{1}{2\Delta}]$-alg{\`e}bres contenant les valeurs de $\kappa$, 
si $f\in G\!D_\kappa(\Gamma;R')$ et 
$f\in R_{\mathcal{C}}^{(\kappa)}(R)$,
alors $f\in G\!D_\kappa(\Gamma;R)$.
\end{prop}

\noindent{\bf Id{\'e}e de d{\'e}monstration:} Par irr{\'e}ductibilit{\'e} g{\'e}om{\'e}trique de 
$\overline{M}$, une section du fibr{\'e} $\underline{\omega}^\kappa$
est uniquement d{\'e}termin{\'e}e par son pull-back
{\`a} la VAHB de Tate.

\noindent{\bf Application : } 
Soit $R\subset \C$ une $\Z[\frac{1}{2\Delta}]$-alg{\`e}bre contenant les valeurs de
$\eta:\mathfrak{c}_0/\mathfrak{nc}_0\rightarrow \C$; 
la s{\'e}ries th{\^e}ta $\theta(z,\eta)$ dont 
on a montr{\'e} qu'elle appartient {\`a} $G\!D_{t/2}(\Gamma)$ est en fait dans 
$G\!D_{t/2}(\Gamma;R)$.
De m{\^e}me pour les s{\'e}ries d'Eisenstein \cite{sh-half1}, 
une fois {\'e}tabli leur $q$-d{\'e}veloppement. 
Nous esp{\'e}rons revenir sur ce point ult{\'e}rieurement.

\vspace{-.4cm}
\section[Tour d'Igusa et formes modulaires de Hilbert $p$-adiques.]{Tour d'Igusa et formes modulaires de Hilbert $p$-adiques}

Une autre ``{presque-application}'' est le r{\'e}sultat fondamental 
d'irr{\'e}ductibilit{\'e} de la tour 
d'Igusa. C'est une presque-application au sens que la m{\'e}thode des 
compactifications toro{\"\i}dales
donne un r{\'e}sultat int{\'e}ressant mais insuffisant pour {\'e}tablir cette 
irr{\'e}ductibilit{\'e}.
Rappelons que, par ailleurs, l'irr{\'e}ductibilit{\'e} g{\'e}om{\'e}trique a {\'e}t{\'e} 
{\'e}tablie dans cette situation
par Ribet \cite{ribet29} et dans une situation plus g{\'e}n{\'e}rale 
(pour les groupes symplectiques ou unitaires
sur des corps totalement r{\'e}els) par Hida \cite{hida-igusa}. Ribet 
utilise la monodromie locale 
en un point ordinaire et Hida utilise la th{\'e}orie de Galois 
des tours de vari{\'e}t{\'e}s de Shimura, due {\`a} Shimura.

Dans toute cette partie, on se limite au cas $D=\Gm$ et donc $M=M^1$ 
de sorte qu'il
y a une VAHB universelle $\mathcal{A}\rightarrow M$.

\subsubsection*{L'invariant de Hasse.}
On suppose toujours $[F:\Q]>1$.

Soit $M^*$, resp. $\overline{M}$ la 
compactification minimale, resp. une compactification toro{\"\i}dale, de $M$.

Rappelons que le prolongement canonique
 $\underline{\omega}_{\mathfrak{G}/\overline{M}}^t$ de
 $\underline{\omega}^t$ {\`a} $\overline{M}$
descend {\`a} $M^*$ en un faisceau ample (on suppose comme d'habitude $\Gamma$ net).
 Ceci est d{\'e}montr{\'e} dans \cite{dimdg} Thm8.6(vi).

Soit $p$ premier, premier avec $\mathfrak{ndc}$. Pour tout entier $m\geq 1$, soit 
$M^*_m=M^*\times \Spec(\Z/p^m\Z)$. 
L'invariant de Hasse $H$ est la section globale de $\underline{\omega}^{(p-1)t}$ sur 
$M^*_1$ d{\'e}finie comme suit.
Soit $V:\mathcal{A}^{(p)}\rightarrow \mathcal{A}$ le morphisme de 
Verschiebung sur $M_1$. Formons 
$V^*:\underline{\omega}_{\mathcal{A}/M_1}\rightarrow\underline{\omega}_{\mathcal{A}^{(p)}/M_1}$ et prenons sa 
puissance ext{\'e}rieure maximale $\bigwedge^dV^*$. Via l'identification canonique 
$\bigwedge^d\underline{\omega}_{\mathcal{A}^{(p)}/M_1}=\left(\bigwedge^d\underline{\omega}_{\mathcal{A}/M_1}\right)
^{\otimes p}$,
on interpr{\`e}te  $\bigwedge^dV^*$ comme une section globale de 
$\underline{\omega}_{\mathcal{A}/M_1}^{(p-1)t}$.
Par le principe de Koecher, c'est une section globale 
de $\underline{\omega}^{(p-1)t}$ sur $M_1^*$. On d{\'e}finit pour tout $m\geq 1$
le lieu ordinaire $M_m^{*,\ord}$ de $M^*_m$ comme l'image inverse du lieu 
$M_1^{*,\ord}$ de $M^*_1$ o{\`u} $H\neq 0$.
$M_1^{*,\ord}$  est un ouvert affine de  $M^*_1$ car $H$ est une 
section globale d'un
faisceau ample; de plus,  une puissance $H^s$ se 
rel{\`e}ve en $E$ sur $\Z_p$, car une section  modulo $p$
d'un faisceau tr{\`e}s ample se rel{\`e}ve; cette section d{\'e}finit un ouvert
affine $M_m^{*,\ord}$ sur $\Z/p^m\Z$ d'{\'e}quation $E\not\equiv 0\pmod p$. 
Lorsque $m$ tend vers l'infini, les compl{\'e}mentaires de ces ouverts 
d{\'e}finissent un voisinage tubulaire ouvert d'{\'e}quation $|E|_p<1$
du diviseur non-ordinaire $E=0$ dans l'espace rigide 
$M^{\rig}$ associ{\'e} aux $M_m$.     

Soit $M_m^{\ord}=M_m^{*,\ord}\cap M_m$ et soit $\overline{M}_m^{*,\ord}=\pi^{-1}(M_m^{*,\ord})$.

\begin{rque}
Par le caract{\`e}re local du principe de Koecher 
(voir la d{\'e}monstration du Th{\'e}or{\`e}me \ref{koecher}),  on a lorsque $d>1$ :
$$\mathrm{H}^0(M_m^{\ord},\mathcal{O}_{M_m^{\ord}})=\mathrm{H}^0(\overline{M}_m^{\ord},\mathcal{O}_{\overline{M}_m^{\ord}})
=\mathrm{H}^0(M_m^{*,\ord},\mathcal{O}_{M_m^{*,\ord}})$$ 

Ceci entra{\^\i}ne que l'ouvert $M_m^{\ord}$ ne peut {\^e}tre affine, 
puisqu'il est distinct de l'ouvert affine $M_m^{*,\ord}$.
\end{rque}

\subsubsection*{Tour d'Igusa et monodromie aux pointes.}
Pour tout $n\geq 1$, la composante neutre 
du sch{\'e}ma fini et plat $\mathcal{A}[p^n]$ est purement
torique au-dessus de $M_m^{\ord}$. 
On d{\'e}finit la tour d'Igusa sur $M_m^{\ord}$ par
$$T_{m,n}=\Isom_{M_m^{\ord}}(\mathfrak{d}^{-1}\otimes \mu_{p^n},\mathcal{A}[p^n]^0)\rightarrow M_m^{\ord} ,\quad n\geq 1$$

Notons que $T_{m,n}\rightarrow M_m^{\ord}$ est 
une suite de $(\mathfrak{o}/p^n\mathfrak{o})^\times$-torseurs 
finis {\'e}tales  puisque
$$(\mathfrak{o}/p^n\mathfrak{o})^\times 
=\Aut(\mathfrak{d}^{-1}\otimes \mu_{p^n}).$$
Le morphismes $T_{m,n}\rightarrow M_m^{\ord}$ sont donc affines.

La question de l'irr{\'e}ductibilit{\'e} g{\'e}om{\'e}trique de la tour d'Igusa 
{\'e}quivaut {\`a} celle de la connexit{\'e} g{\'e}om{\'e}trique
 de ces rev{\^e}tements (puisque la base est g{\'e}om{\'e}triquement connexe).
Ceci {\'e}quivaut aussi {\`a} la surjectivit{\'e} des repr{\'e}sentations associ{\'e}es {\`a} ces rev{\^e}tements
(apr{\`e}s choix d'un point base g{\'e}om{\'e}trique):
$$\rho_n: \pi_1(M_1\otimes \overline{\F}_p)\rightarrow (\mathfrak{o}/p^n\mathfrak{o})^\times.$$

Cette surjectivit{\'e} a {\'e}t{\'e} {\'e}tablie par Ribet \cite{ribet29} 
et Hida \cite{hida-igusa} dans un cadre plus g{\'e}n{\'e}ral (voir aussi \cite{hida-padic}).
Chai, \cite{chai} Sect.4.6, a propos{\'e} une autre approche via les compactifications toro{\"\i}dales, 
inspir{\'e}e de \cite{FaCh} V.7.
Malheureusement, cette approche n'est pas concluante dans toutes les situations 
de compactifications toro{\"\i}dales,
et en particulier pas dans le cas des vari{\'e}t{\'e}s de Hilbert. Expliquons pourquoi. 
L'id{\'e}e est de calculer la monodromie locale au voisinage d'une pointe
non-ramifi{\'e}e dans la compactification minimale.
Elle ne permet cependant que d'obtenir un sous-groupe ferm{\'e} de $(\mathfrak{o}\otimes \Z_p)^\times$
qui est ouvert si la conjecture de Leopoldt est vraie.

Pour all{\'e}ger la notation dans l'argument ci-dessous, on omet (mais on sous-entend) l'extension des scalaires {\`a} 
$\overline{\F}_p$. Tout d'abord le cas $F=\Q$ est connu. Par cons{\'e}quent, la VAHB obtenue en formant le
produit tensoriel de la courbe elliptique universelle au-dessus de la courbe modulaire $Y$
par $\mathfrak{o}$ (puis si n{\'e}cessaire une $\mathfrak{c}$-isog{\'e}nie pour obtenir une $\mathfrak{c}$-polarisation)  fournit un morphisme de $Y$ dans la vari{\'e}t{\'e} de Hilbert qui permet d'inclure 
$(\Z/p^n\Z)^\times$ dans l'image de $\rho_n$.

Pour obtenir un gros sous-groupe du noyau $(\mathfrak{o}/p^n\mathfrak{o})^\times_1$ 
de la norme, Chai proc{\`e}de comme suit. Soit $\pi:\overline{M}_1\rightarrow M_1^*$.
Soit $\mathcal{C}$ une pointe non-ramifi{\'e}e de $M_1^*$ et 
soit $R$ la compl{\'e}tion de l'anneau local de $M_1^*$ en $\mathcal{C}$.
Posons $\overline{S}_{\Sigma^{\mathcal{C}}}=
\coprod_{\sigma\in\Sigma^{\mathcal{C}}} \overline{S}_{\sigma}$ et   
 $\overline{S}'_{\Sigma^{\mathcal{C}}}= 
\overline{S}_{\Sigma^{\mathcal{C}}}/\mathfrak{o}^\times_{\mathcal{C}}$. 
 Soit $(\overline{S}_{\Sigma^{\mathcal{C}}})^0$
(resp. $(\overline{S}'_{\Sigma^{\mathcal{C}}})^0$) le compl{\'e}mentaire du diviseur 
{\`a} l'infini $(\overline{S}_{\Sigma^{\mathcal{C}}})^\infty$ 
(resp. $(\overline{S}'_{\Sigma^{\mathcal{C}}})^\infty$). L'image de 
$\Spec(R)$ dans $M_1^*$ est contenue dans $M_1^{*\, \ord}$, 
 gr{\^a}ce {\`a} la description de la VAHB universelle sur $\Spec(R)$
 comme quotient de Mumford. On a des morphismes de groupes
$$\pi_1((\overline{S}'_{\Sigma^{\mathcal{C}}})^0)\rightarrow
 \pi_1(\Spec(R)\bs\{\mathcal{C}\})\rightarrow \pi_1(M_1^{\ord})$$

On observe que l'on a un rev{\^e}tement {\'e}tale 
$$\overline{\xi}_n:\overline{T}_{1,n}=\Isom_{\overline{M}_1}
(\mathfrak{d}^{-1}\otimes \mu_{p^n},\mathfrak{G}_1[p^n]^0)\rightarrow
\overline{M}_1^{\ord}=\pi^{-1}(M_1^{*\, \ord})$$
prolongeant le rev{\^e}tement {\'e}tale $T_{1,n}\rightarrow M_1^{\ord} $.

Le pull-back de $\overline{\xi}_n$ {\`a} $\overline{S}_{\Sigma^{\mathcal{C}}}$
est trivial car la description de 
$\mathfrak{G}_1|_{\overline{S}_{\Sigma^{\mathcal{C}}}}$ comme vari{\'e}t{\'e} de Tate 
montre que sur $\overline{S}_{\Sigma^{\mathcal{C}}}$, on a un isomorphisme 
canonique $\mathfrak{G}_1[p^n]^0\cong \mu_{p^n}\otimes \mathfrak{a}$ compatible 
avec l'action de $\mathfrak{o}^\times_{\mathcal{C}}$. Ainsi 
$\overline{S}_{\Sigma^{\mathcal{C}}}\times\overline{T}_{1,n}=
\overline{S}_{\Sigma^{\mathcal{C}}}\times (\mathfrak{o}/p^n)^\times$
et l'action de $\mathfrak{o}^\times_{\mathcal{C}}$ sur le membre de 
gauche correspond {\`a} l'action diagonale sur le membre de droite. 
L'image de $\rho_n$ contient donc l'image de $\mathfrak{o}^\times_{\mathcal{C}}$
dans $(\mathfrak{o}/p^n)^\times$. 

On notera que le  rev{\^e}tement connexe {\'e}tale 
$\overline{S}_{\Sigma^{\mathcal{C}}}^\infty\rightarrow 
(\overline{S}'_{\Sigma^{\mathcal{C}}})^\infty$ de groupe 
$\mathfrak{o}_{\mathcal{C}}^\times$  n'est pas de type fini 
(seulement localement de type fini). N{\'e}anmoins on peut former une tour 
de rev{\^e}tements finis connexes {\'e}tales en quotientant 
$\overline{S}_{\Sigma^{\mathcal{C}}}^\infty$
par $(\mathfrak{o}_{\mathcal{C}}^\times)^{ p^n}$.

Comme le pull-back de $\mathfrak{G}_1\rightarrow \overline{M}_1$ {\`a} 
$(\overline{S}'_{\Sigma^{\mathcal{C}}})^\infty$ est le quotient du tore 
de Tate sur $ (\overline{S}_{\Sigma^{\mathcal{C}}})^\infty$ 
par $\mathfrak{o}_{\mathcal{C}}^\times$, on voit que 
$\mathfrak{o}_{\mathcal{C}}^\times/(\mathfrak{o}_{\mathcal{C}}^\times)^{p^n}$ 
est contenue dans l'image de $\rho_n$ dans $(\mathfrak{o}/p^n\mathfrak{o})^\times$.
En passant {\`a} la limite projective, on obtient  l'image de $\pi_1(M_1^{\ord})$ dans 
$(\mathfrak{o}\otimes \Z_p)^\times$ contient bien $\Z_p^\times\cdot 
\overline{\mathfrak{o}^\times}$.

\subsubsection*{Formes modulaires $p$-adiques.}

Cette th{\'e}orie est une application directe de l'irr{\'e}ductibilit{\'e} de la tour d'Igusa.

On a d{\'e}fini $\overline{T}_{m,n}$ dans la section pr{\'e}c{\'e}dente; c'est un torseur fini {\'e}tale sur 
$\overline{M}_m^{\ord}$; il est donc affine sur ce sch{\'e}ma. Ce dernier n'est pas affine, mais on a vu
que 
$\mathrm{H}^0(\overline{M}_m^{\ord},\mathcal{O}_{\overline{M}_m^{\ord}})=\mathrm{H}^0(M_m^{\ord},\mathcal{O}_{M_m^{\ord}})$
et que le spectre de cet anneau est $M_m^{*,\ord}$. Soit 
$$V_{m,n}=\mathrm{H}^0(T_{m,n},\mathcal{O}_{T_{m,n}})=\mathrm{H}^0(\overline{T}_{m,n},\mathcal{O}_{\overline{T}_{m,n}})$$
l'alg{\`e}bre des fonctions r{\'e}guli{\`e}res sur $T_{m,n}$. 
En passant {\`a} la limite inductive sur $n$ puis projective sur $m$, on obtient l'anneau $\bV$ des 
fonctions sur le ind-pro-sch{\'e}ma
$$T_{\infty,\infty}\rightarrow M_\infty^{\ord}$$

Notons que comme les morphismes de transition $T_{m,n+1}\rightarrow T_{m,n}$ sont affines, $T_{m,\infty}$
est un sch{\'e}ma et $T_{\infty,\infty}$ est un sch{\'e}ma formel sur $M_\infty$.

$M_\infty^{\ord}$ est le sch{\'e}ma formel obtenu en retirant 
le lieu supersingulier $E\equiv 0\pmod p$ de la
compl{\'e}tion $M_\infty$ du sch{\'e}ma $M_{\Z_p}$ le long de sa fibre sp{\'e}ciale ($E$ d{\'e}signe un rel{\`e}vement quelconque de $H^s$ sur $\Z_p$ comme pr{\'e}c{\'e}demment). Dans le langage de la g{\'e}om{\'e}trie rigide, on retire le voisinage tubulaire ouvert de rayon $1$ autour du diviseur supersingulier.

Soit $\cA_p$ la cat{\'e}gorie des alg{\`e}bres $p$-adiquement compl{\`e}tes.
Le ind-pro-sch{\'e}ma $T_{\infty,\infty}$ repr{\'e}sente le foncteur qui associe {\`a} $R\in \cA_p$ l'ensemble des
classes d'isomorphisme de $(A,\iota,\lambda,\alpha,\phi)_{/R}$ o{\`u} $(A,\iota,\lambda,\alpha)$ d{\'e}finit un point de $M\times \Spec(R)$ et $\phi:\mu_{p^\infty}\otimes\mathfrak{d}^{-1}\rightarrow A[p^\infty]$
est une immersion ferm{\'e}e de groupes $p$-divisibles. De fa{\c c}on {\'e}quivalente, 
 on peut d{\'e}finir $\phi$ comme un isomorphisme de groupes formels
$\widehat{\mathbb{G}}_m\otimes\mathfrak{d}^{-1}\rightarrow \widehat{A}$ (compl{\'e}tions des sch{\'e}mas en groupes
le long de leur section unit{\'e}). On appelle un tel plongement $\phi$ une rigidification du groupe formel 
$\widehat{A}$.

 Pour toute $\Z_p$-alg{\`e}bre $p$-adiquement compl{\`e}te $R$, on pose
$$\bV(\mathfrak{c},\mathfrak{n};R)=\bV\widehat{\otimes} R$$
On l'appelle l'alg{\`e}bre des formes modulaires $p$-adiques de niveau 
auxiliaire $\mathfrak{n}$ {\`a} coefficients dans $R$. Elle est munie 
d'une action du groupe de Galois 
$\mathcal{T}_1(\Z_p)=(\mathfrak{o}\otimes\Z_p)^\times$ du rev{\^e}tement 
$T_{\infty,\infty}\rightarrow M^{\ord}_{\infty}$. On pose pour tout 
homomorphisme continu $\kappa:\mathcal{T}_1(\Z_p) \rightarrow R^\times$
$$\bV_\kappa(\mathfrak{c},\mathfrak{n};R)=\bV(\mathfrak{c},\mathfrak{n};R)[-\kappa]
$$

La donn{\'e}e d'une rigidification $\phi:\widehat{\Gm}\otimes\mathfrak{d}^{-1}\rightarrow \widehat{A}$ 
induit un isomorphisme 
$\phi^*:\mathfrak{o}\otimes\mathcal{O}_{T_{\infty,\infty}}\cong
\underline{\omega}$ par le principe G{\'e}om{\'e}trie Alg{\'e}brique-G{\'e}om{\'e}trie Formelle
 pour le sch{\'e}ma propre $\mathcal{A}$ sur $ M$.
Soit $\mathfrak{M}_\infty$ la compl{\'e}tion  de $\mathfrak{M}$
formelle le long de la fibre sp{\'e}ciale. On obtient un morphisme
 de $M_\infty$-sch{\'e}mas formels $T_{\infty,\infty}\rightarrow
\mathfrak{M}_\infty$; il
fournit un homomorphisme $j$ de 
l'alg{\`e}bre des formes classiques vers celle des formes $p$-adiques.
$$j:G(\mathfrak{c},\mathfrak{n})=\mathrm{H}^0(\mathfrak{M},
\mathcal{O}_{\mathfrak{M}})\rightarrow  \bV(\mathfrak{c},\mathfrak{n};R)$$

\vspace{-2mm}
\begin{rque}
Si $\kappa$ est un caract{\`e}re alg{\'e}brique de $\mathcal{T}_1$, $j$ envoie 
les formes classiques de poids $\kappa$ vers les formes 
$p$-adiques de poids $\kappa$. Cette application
est injective (voir ci-dessous) mais est loin d'{\^e}tre surjective, 
le module de droite {\'e}tant  de rang infini. La th{\'e}orie de Hida 
permet de montrer que sur les sous-modules des formes ordinaires,
on a un isomorphisme. Il serait int{\'e}ressant par une g{\'e}n{\'e}ralisation 
appropri{\'e}e de m{\'e}thodes de
 Mazur-Coleman d'{\'e}tendre aux formes de Hilbert
surconvergentes de pente fix{\'e}e un tel r{\'e}sultat de classicit{\'e}.     
\end{rque} 

\vspace{-4mm}
\begin{theo}
1)  Pour toute $(R,\mathfrak{n})$-pointe non-ramifi{\'e}e $\mathcal{C}$,
 on a un homomorphisme injectif 
$$\ev_{\mathcal{C}}:\bV(\mathfrak{c},\mathfrak{n};R)
\rightarrow R[[q^\xi;\xi\in X_+\cup\{0\}]],$$
o{\`u}  $X=\mathfrak{cb}^2$, qui est 
donn{\'e} par l'{\'e}valuation sur la VAHB de Tate sur 
$\overline{S}^0_{\mathcal{C},\sigma}$,
munie de sa rigidification canonique (pour tout 
$\sigma\in \Sigma^{\mathcal{C}}$).
De plus, pour toute $R\in \cA_p$, si $f\in \bV(\mathfrak{c},\mathfrak{n};R)$
et pour toute $(R,\mathfrak{n})$-pointe  non-ramifi{\'e}e $\mathcal{C}$ de $M$, on a $\ev_{\mathcal{C}}(f)\in pR[[q^\xi]]$, alors
$f\in p\bV(\mathfrak{c},\mathfrak{n};R)$.

2)  le morphisme $\ev_{\mathcal{C}}$ est compatible via l'homomorphisme $j$
avec celui d{\'e}fini pour les formes classiques.
\end{theo}

\noindent {\bf D{\'e}monstration : }
 1) Rappelons la propri{\'e}t{\'e} universelle de 
$\overline{T}_{\infty,\infty}=
\Isom_{\overline{M}_\infty^{\ord}}(\mu_{p^\infty}\otimes\mathfrak{d}^{-1},
\mathfrak{G}[p^\infty]^0)$:
Pour tout sch{\'e}ma semi-ab{\'e}lien $G\rightarrow S$ sur un sch{\'e}ma formel $S$ et pour 
tout morphisme $G\rightarrow \mathfrak{G}$ au-dessus d'un morphisme
$S\rightarrow \overline{M}_\infty$, il y a une bijection canonique entre 
l'ensemble  des 
rigidifications $\mu_{p^\infty}\otimes\mathfrak{d}^{-1}\cong G[p^\infty]^0$
et celui des rel{\`e}vements $S\rightarrow \overline{T}_{\infty,\infty}$ du
morphisme donn{\'e} $S\rightarrow \overline{M}_{\infty}$.

 La rigidification canonique
$\mu_{p^\infty}\otimes \mathfrak{d}^{-1}\cong G_\sigma[p^\infty]^0$
fournit donc un morphisme de sch{\'e}mas formels de la compl{\'e}tion 
$p$-adique de 
$\overline{S}_{\mathcal{C},\sigma}$ vers  $\overline{T}_{\infty,\infty}$.
Ce morphisme est {\'e}tale. L'irr{\'e}ductibilit{\'e} du sch{\'e}ma formel 
$\overline{T}_{\infty,\infty}$ entra{\^\i}ne qu'une section globale 
d'un faisceau localement libre sur $\overline{T}_{\infty,\infty}$ nulle sur
$\overline{S}_{\mathcal{C},\sigma}$ est identiquement nulle.  

2) Cet {\'e}nonc{\'e} pour une alg{\`e}bre $R$ r{\'e}sulte de 1) appliqu{\'e} {\`a} l'alg{\`e}bre $R/pR$.
\hfill$\square$

\begin{cor} Pour tout poids $\kappa\in \Z[J_F]$,
$j$ restreint {\`a} $G_\kappa(\mathfrak{c},\mathfrak{n};R)$ est injectif.  
\end{cor}

\begin{rque}
Si $R$ est plate sur $\Z_p$, par ind{\'e}pendance lin{\'e}aire des caract{\`e}res
alg{\'e}briques distincts, on d{\'e}duit du corollaire que l'application $j$ elle-m{\^e}me est injective. Par
contre, ce n'est pas le cas pour une $\Z_p$-alg{\`e}bre $p$-adique quelconque (c'est faux pour $\F_p$).
\end{rque}

\begin{cor}
L'homomorphisme $\ev_{\mathcal{C}}$ est d'image ferm{\'e}e dans
$R[[q^\xi;\xi\in X_+\cup\{0\}]]$ muni de la topologie de la convergence coefficient
par coefficient.
\end{cor}

\noindent{\bf D{\'e}monstration:} Si $f_n\in \bV(\mathfrak{c},\mathfrak{n};R)$ et 
$\ev_{\mathcal{C}}(f_{n+1})\equiv
\ev_{\mathcal{C}}(f_{n})\pmod{p^n}$, alors, $f_{n+1}\equiv f_n\pmod{p^n\bV(\mathfrak{c},\mathfrak{n};R)}$
et donc $f_n$ converge dans l'alg{\`e}bre $p$-adiquement compl{\`e}te $\bV(\mathfrak{c},\mathfrak{n};R)$.   

\medskip
Pour chaque $r\geq 1$; on va d{\'e}finir un homomorphisme
d'alg{\`e}bres 
$$j_r:G(\mathfrak{c},\mathfrak{n}p^r; R)\rightarrow \bV(\mathfrak{c},\mathfrak{n};R)$$ 
compatible avec les poids.

Le th{\'e}or{\`e}me de repr{\'e}sentabilit{\'e} de Mumford (GIT) entra{\^\i} ne qu'il existe un 
$\Z_p$-sch{\'e}ma quasi-projectif de type fini $M_0(\mathfrak{c},\mathfrak{n}p^r)$ classifiant les VAHB $\mathfrak{c}$-polaris{\'e}es 
$(A,\iota,\lambda,\alpha)_S$ 
munies d'un $\mathfrak{o}/p^r$-module cyclique $C$ (c'est-{\`a}-dire, localement libre sur 
$\mathfrak{o}\otimes\mathcal{O}_S$ de rang $p^r$
et annul{\'e} par $p^r$
exactement). La normalisation de ce sch{\'e}ma dans 
$M_1(\mathfrak{c},\mathfrak{n}p^r)\times \Spec\,\Q_p$ d{\'e}finit un mod{\`e}le entier sur $\Z_p$ 
de $M_1(\mathfrak{c},\mathfrak{n}p^r)$. Il est muni d'un morphisme propre et plat $\pi_r$ vers $M$.

Les faisceaux $\underline{\omega}^\kappa$ sur $M_1(\mathfrak{c},\mathfrak{n}p^r)$ sont les pull-back par $\pi_r$
des $\underline{\omega}^\kappa$ sur $M$.
Soit $M_{\infty,r}$ la compl{\'e}tion formelle de ce 
sch{\'e}ma le long de sa fibre sp{\'e}ciale. 
La d{\'e}finition de $T_{\infty,r}$ comme solution du probl{\`e}me de modules des immersions ferm{\'e}es
$\mu_{p^r}\otimes\mathfrak{d}^{-1}\hookrightarrow A[p^r]$ montre qu'il y a une immersion ouverte 
canonique
$T_{\infty,r}\hookrightarrow M_{\infty,r}$. Cette immersion est d'image dense dans 
l'une des composantes irr{\'e}ductibles de $M_{\infty,r}$.   
Soit $\mathfrak{M}_{\infty,r}=\mathfrak{M}_{\infty}\times_{M_{\infty}}M_{\infty,r}$.
On a vu qu'il y a un $M_\infty$-morphisme canonique $T_{\infty,\infty}\rightarrow \mathfrak{M}_\infty$; 
on vient d'autre part de construire un $M_\infty$-morphisme $T_{\infty,\infty}\rightarrow T_{\infty,r}
\rightarrow M_{\infty,r}$; on obtient donc un morphisme
$T_{\infty,\infty}\rightarrow \mathfrak{M}_{\infty,r}$
qui fournit l'homomorphisme d'alg{\`e}bres $j_{r}$ par restriction au lieu ordinaire sur les fonctions.

Soit $M_{\infty,\infty}=\limproj M_{\infty,r}$.
En fait, en passant {\`a} la limite projective sur les morphismes 
$T_{\infty,\infty}\rightarrow \mathfrak{M}_{\infty,r}$, on voit facilement que
$$T_{\infty,\infty}\rightarrow \mathfrak{M}_{\infty,\infty}\quad {\rm et}\quad
T_{\infty,\infty}\rightarrow M_{\infty,\infty}^{\ord}$$

sont des immersions ouvertes.
On peut d{\'e}finir un fibr{\'e} inversible 

Si l'on ne tient pas {\`a} pr{\'e}server l'int{\'e}gralit{\'e} (c'est {\`a} dire que l'on
 se limite aux $\Q_p$-alg{\`e}bres),
on peut aussi formuler cette construction en consid{\'e}rant l'espace 
rigide $T_{\infty,r}^{\rig}$
associ{\'e} {\`a} $T_{\infty,r}$; c'est un ouvert de $M_1(\mathfrak{c},
\mathfrak{n}p^r)^{\rig}$;
l'avantage de cette approche est que ces espaces rigides sont lisses et connexes 
(alors que le sch{\'e}mas formel de $M_1(\mathfrak{c},\mathfrak{n}p^r)$ a un 
grand nombre de composantes irr{\'e}ductibles).
  
Soit $M_1(\mathfrak{c},\mathfrak{n}p^r)^{\ord}$ l'image inverse dans 
$M_1(\mathfrak{c},\mathfrak{n}p^r)^{\rig}$
de $M_1(\mathfrak{c},\mathfrak{n})^{\ord}$ par le morphisme
d'oubli du groupe cyclique d'ordre $p^r$.
On voit que $T_{\infty,r}^{\rig}$ est contenu dans $M_1(\mathfrak{c},
\mathfrak{n}p^r)^{\ord}$ .
Par cons{\'e}quent le $\mathcal{T}_1$-torseur rigide $\mathfrak{M}(\mathfrak{c},
\mathfrak{n}p^r)^{\rig}$
sur $M_1(\mathfrak{c},\mathfrak{n}p^r)^{\rig}$ classifiant les VAHB munies 
d' un $\mathfrak{o}$-groupe cyclique
d'ordre $p^r$ et d'une $\mathfrak{o}\otimes\mathcal{O}_S$-base de 
$ \underline{\omega}$, est image de 
$T_{\infty,\infty}^{\rig}$ par le m{\^e}me argument que pr{\'e}c{\'e}demment.

\begin{theo} Pour tout $r\geq 1$, pour tout $\kappa\in\Z[J_F]$,  
l'homomorphisme d'alg{\`e}bres des formes classiques de niveau $p^r$ 
vers les formes $p$-adiques
$$j_{r,\kappa}:G_\kappa(\mathfrak{c},\mathfrak{n}p^r; R)\rightarrow \bV(\mathfrak{c},\mathfrak{n};R)$$ 
est injectif.
\end{theo}

\noindent{\bf D{\'e}monstration : } 
On sait d'une part que l'int{\'e}gralit{\'e} est pr{\'e}serv{\'e}e; on consid{\`e}re d'autre
part les espaces rigides  $\mathfrak{M}_{\infty,r}^{\rig}$ resp. $T_{\infty,\infty}$,
qui sont connexes par l'irr{\'e}ductibilit{\'e} de la tour d'Igusa. Ceci entra{\^\i} ne
l'injectivit{\'e} du $q$-d{\'e}veloppement des formes classiques de niveau $\mathfrak{n}p^r$ ainsi que celui
des formes $p$-adiques aux $(R,\mathfrak{n})$-pointes de $M$ relev{\'e}es {\`a} l'aide de la rigidification 
canonique $\phi_{\can}$ de la VAHB de 
Tate associ{\'e}e {\`a} $\mathcal{C}$. On d{\'e}duit alors le th{\'e}or{\`e}me de 
la compatibilit{\'e} des morphismes de $q$-d{\'e}veloppement. \hfill $\square$

\begin{rque} 1) {\`A} la diff{\'e}rence du cas du niveau premier {\`a} $p$, l'utilisation de la g{\'e}om{\'e}trie
rigide simplifie l'argument car la fibre sp{\'e}ciale de $M_1(p^r)$ est 
compliqu{\'e}e.

2) L'homomorphisme $j_r$ somme des $j_{r,\kappa}$ est injectif si $R$ est $\Z_p$-plate.
\end{rque}

En fait, la compatibilit{\'e} des morphismes de $q$-d{\'e}veloppement est vraie sur un ensemble de 
pointes $p$-adiques plus grand que l'ensemble fini des rel{\`e}vements standards des $(R,\mathfrak{n})$-pointes 
$\mathcal{C}$ de $M$ donn{\'e}s par la rigidification canonique $\phi_{\can}$ de la VAHB de 
Tate associ{\'e}e {\`a} $\mathcal{C}$.
Cet ensemble est appel{\'e} l'ensemble des pointes non-ramifi{\'e}es. On va le d{\'e}finir et le 
caract{\'e}riser {\`a} 
l'aide de l'immersion ouverte (de sch{\'e}mas formels ou rigides)
$$T_{\infty,\infty}\rightarrow \mathfrak{M}_{\infty,\infty}$$

\begin{defin} Une pointe de $T_{\infty,\infty}$ est 
un couple constitu{\'e}
d'une $(R,\mathfrak{n})$-pointe $\mathcal{C}$ de $M$ et de la classe d'isomorphisme $i(a)$ 
($a\in(\mathfrak{o}\otimes\Z_p)^\times$)
du morphisme canonique
$$\Tate_{\mathcal{C}}(q)\rightarrow \overline{T}_{\infty,\infty}$$
donn{\'e} par la VAHB de Tate $(\Tate_{\mathcal{C}}(q),\iota_{\can},\lambda_{\can},\alpha_{\can})$
au-dessus de $\overline{S}_{\mathcal{C},\Sigma^{\mathcal{C}}}\rightarrow \overline{M}$,
munie de la rigidification $\phi_{\can}\circ[a])$. 
 
\end{defin}

Deux telles classes d'isomorphisme $i(a)$ et $i(a')$ co{\"\i}ncident si et seulement
si $a$ et $a'$ sont congrus modulo $\overline{\mathfrak{o}^\times}$ 
car ce groupe agit par automorphismes du carr{\'e} cart{\'e}sien

$$\begin{array}{ccc}
\Tate_{\mathcal{C}}(q)&\rightarrow& \overline{T}_{\infty,\infty}\\
\downarrow&&\downarrow \\
\overline{S}_{\mathcal{C}}&\rightarrow& \overline{M}
\end{array} $$

Soit $P^T_{\infty}$ 
l'ensemble des pointes de $T_{\infty,\infty}$. La fibre de $P^T_\infty$ au-dessus d'une pointe de $M$
est donc un $(\mathfrak{o}\otimes\Z_p)^\times/\overline{\mathfrak{o}^\times}$-torseur.  

Notons $P_\infty$ l'ensemble des pointes de $M_1(\mathfrak{c},\mathfrak{n}p^\infty)=\limproj M_1(\mathfrak{c},\mathfrak{n}p^r)$.
C'est un espace $p$-adique compact et
il r{\'e}sulte imm{\'e}diatement du th{\'e}or{\`e}me d'approximation forte 
et de la d{\'e}composition d'Iwasawa que 
$P_\infty$ est un fibr{\'e} au-dessus de l'ensemble des pointes de 
$\Gamma_1(\mathfrak{c},\mathfrak{n})$, ses fibres {\'e}tant des copies de 
$$U(\Z_p)\bs G(\Z_p)/
\overline{\mathcal{T}_1(\Z_p)}U(\Z_p)$$
o{\`u} $U$ est le radical unipotent du Borel sup{\'e}rieur.
Par action sur l'ouvert des vecteurs primitifs 
$((\mathfrak{o}\otimes\Z_p)^2)_{\mathrm{prim}}$ de $(\mathfrak{o}\otimes\Z_p)^2$,
on identifie ce quotient {\`a}
$$U(\Z_p)\bs
((\mathfrak{o}\otimes\Z_p)^2)_{\mathrm{prim}}/\overline{\mathfrak{o}^\times}
$$
 Ce quotient s'identifie {\`a} l'ensemble des classes d'homoth{\'e}tie par un scalaire
de $\overline{\mathfrak{o}^\times}$ de vecteurs primitifs 
$\begin{pmatrix}a\\c\end{pmatrix}$
dans $(\mathfrak{o}\otimes\Z_p)^2$, o{\`u}
$c\in (\mathfrak{o}\otimes\Z_p)$ et $a$ 
parcourt un syst{\`e}me de repr{\'e}sentants de $(\mathfrak{o}\otimes\Z_p)/(c)$.
Notons  $\begin{bmatrix}a\\c\end{bmatrix}$ une telle classe. 
Tout {\'e}l{\'e}ment de $P_\infty$ s'{\'e}crit donc comme un couple constitu{\'e} d'une 
$(R,\mathfrak{n})$-pointe de $M$ et d'une classe $\begin{bmatrix}a\\c\end{bmatrix}$.
On introduit la notion de profondeur d'une pointe:
$$\mathrm{prof}_p(\begin{bmatrix}a\\c\end{bmatrix})=\ord_p(c)$$

Les pointes de profondeur infinie sont appel{\'e}es les pointes non ramifi{\'e}es. 
Soit $P_\infty^{nr}$ le sous-espace des pointes non-ramifi{\'e}es; c'est un
torseur sous 
$(\mathfrak{o}\otimes\Z_p)^\times/\overline{\mathfrak{o}^\times}$ au-dessus de l'ensemble des pointes de 
$M$.

\begin{lemme} 
L'immersion ouverte  $T_{\infty,\infty}\hookrightarrow M_{\infty,\infty}$ identifie $P^T_\infty$ et
$P^{nr}_\infty$.
\end{lemme}

 Soit $\mathcal{C}(P_\infty^{nr};R)$ l'espace des fonctions continues
sur $P_\infty^{nr}$ {\`a} valeurs dans une alg{\`e}bre $p$-adique $R$. On peut rep{\'e}rer une telle pointe
comme un couple constitu{\'e} d'une $(R,\mathfrak{n})$-pointe $c$ et d'un {\'e}l{\'e}ment $a$ de 
$(\mathfrak{o}\otimes \Z_p)^\times$.
Pour tout $f\in\bV(\mathfrak{c},\mathfrak{n};R)$, le $q$-d{\'e}veloppement en une pointe non ramifi{\'e}e $(c,a)$ est 
bien d{\'e}fini: c'est l'{\'e}valuation sur la vari{\'e}t{\'e} de  $\Tate_c(q)$ munie de la rigidification
$$\phi_c\circ[a]:\widehat{\mathbb{G}}_m\otimes\mathfrak{d}^{-1}\rightarrow \widehat{\Tate}_c(q),
\quad t\mapsto \phi_c([a](t))$$
Ces fl{\`e}ches de $q$-d{\'e}veloppement induisent un homomorphisme d'alg{\`e}bres
$$\ev:\bV(\mathfrak{c},\mathfrak{n};R)\rightarrow \mathcal{C}(P_\infty^{nr};R)$$
Le noyau s'appelle l'espace des formes modulaires $p$-adiques cuspidales.
C'est un id{\'e}al de $\bV(\mathfrak{c},\mathfrak{n};R)$ not{\'e} $\bV^{\cusp}(\mathfrak{c},\mathfrak{n};R)$

Par ce qui pr{\'e}c{\`e}de, $\bV(\mathfrak{c},\mathfrak{n};R)$ contient pour chaque poids $\kappa\geq 2t$
$$G_{\kappa}(\mathfrak{c},\mathfrak{n}p^\infty; R)=
\bigcup_{r\geq 0} G_{\kappa}(\mathfrak{c},\mathfrak{n}p^r; R)$$

ainsi que 
$$G(\mathfrak{c},\mathfrak{n}; R)=\bigoplus_{\kappa\geq 0}
G_{\kappa}(\mathfrak{c},\mathfrak{n}; R)$$

Il est de plus {\'e}vident que les formes cuspidales classiques sont envoy{\'e}es dans 
les formes $p$-adiques cuspidales.
Hida a d{\'e}montr{\'e}  

\begin{theo} Soit $R$ une $\Z_p$-alg{\`e}bre plate et $p$-adiquement compl{\`e}te. Alors

{\rm (i)}  $G^{\cusp}(\mathfrak{c},\mathfrak{n}; R)[\frac{1}{p}]\cap 
\bV(\mathfrak{c},\mathfrak{n};R)$ est dense dans 
$\bV^{\cusp}(\mathfrak{c},\mathfrak{n};R)$,

{\rm (ii)} De m{\^e}me, pour tout $\kappa\geq 2t$, $G_{\kappa}^{\cusp}(\mathfrak{c},\mathfrak{n}p^\infty; R)$ est dense
dans $\bV^{\cusp}(\mathfrak{c},\mathfrak{n};R)$ sur les parties ordinaires. 

\end{theo}

Le premier {\'e}nonc{\'e} est d{\'e}montr{\'e} dans \cite{hida-pel}, 
essentiellement pour toutes les vari{\'e}t{\'e}s PEL  (voir Th.2.2 et Cor.3.4). 
Le second est d{\'e}montr{\'e} dans \cite{hida-iwasawa}.
L'{\'e}nonc{\'e} analogue sans limitation aux parties cuspidales n'est pas {\'e}tabli
mais il est conjectur{\'e}.

\medskip
\subsubsection*{Formes modulaires de Hilbert de poids demi-entier $p$-adiques.}

Nous allons donner une d{\'e}finition purement alg{\'e}bro-g{\'e}om{\'e}trique des formes modulaires de Hilbert 
$p$-adiques de poids demi-entier. L'avantage de ce point de vue est qu'il {\'e}vite de recourir
{\`a} la multiplication par une s{\'e}rie th{\^e}ta de poids $t/2$ pour passer du poids demi-entier au poids entier pour
{\'e}tablir les propri{\'e}t{\'e}s arithm{\'e}tiques des formes de poids demi-entier. De cette mani{\`e}re,
on peut {\'e}tablir directement diff{\'e}rents r{\'e}sultats obtenus par Hida dans 
\cite{hida-annarbor}
 Sect.2, h.1-h5.
Nous reviendrons sur ce point dans un travail ult{\'e}rieur.

Soit $p$ premier, premier avec  $2\N(\mathfrak{cdn})$.
Pour chaque couple $(m,n)$ d'entiers, $m,n\geq 1$, on a d{\'e}fini
 le $(\mathfrak{o}/p^n\mathfrak{o})^\times$-torseur $\widetilde{T}_{m,n}$ fini
{\'e}tale sur $M_{m}^{\ord}$.

  On consid{\`e}re le $T_{m,n}$-sch{\'e}ma donn{\'e} par  
$$T^+_{m,n}=\Isom_{T_{m,n}}(\mathcal{O}_{T_{m,n}},0^*\Omega_{\mathcal{L}/\mathcal{A}})$$
c'est un $\Gm$-torseur sur $T_{m,n}$. 

Chaque sch{\'e}ma $T^+_{m,n}$ est g{\'e}om{\'e}triquement connexe par irr{\'e}ductibilit{\'e} de 
la tour d'Igusa. Soit $V^+_{m,n}$ l'anneau des fonction sur $T^+_{m,n}$. Il est muni
d'une action de $\Gm$. La partie de $V^+_{m,n}$ sur laquelle $\Gm$ agit trivialement n'est autre 
que $V_{m,n}$.
Consid{\'e}rons celle o{\`u} $\Gm$ agit par $m_{1}:y\mapsto y$. Notons-la $V\!D_{m,n}$. On a une
action naturelle de $V_{m,n}$ sur $V\!D_{m,n}$ induite par la multiplication dans $V^+_{m,n}$.
Pour toute alg{\`e}bre $p$-adiquement compl{\`e}te $R$, on d{\'e}finit le module des formes modulaires 
$p$-adiques de Hilbert {\`a} coefficients dans $R$ comme
$$\bV\!\D(\mathfrak{c},\mathfrak{n};R)=\left(\underset{m}{\limproj}\,\underset{n}{\limind}\, 
V\!D_{m,n}\right)\widehat{\otimes} R$$
C'est naturellement un $\bV(\mathfrak{c},\mathfrak{n};R)$-module.

On a un morphisme naturel $T^+_{\infty,\infty}\rightarrow \mathfrak{M}^+$ d{\'e}fini par pull-back
{\`a} $\mathfrak{M}^+$ de $T_{\infty,\infty}\rightarrow \mathfrak{M}$. Il fournit un morphisme
$$j_{d}:G\!D(\mathfrak{c},\mathfrak{n};R)\rightarrow \bV\!\D(\mathfrak{c},\mathfrak{n};R)$$
compatible avec les actions des formes modulaires de poids entier classiques resp. $p$-adiques
pour l'homomorphisme
$$j:G(\mathfrak{c},\mathfrak{n};R)\hookrightarrow \bV(\mathfrak{c},\mathfrak{n};R)$$

En particulier, pour tout poids demi-entier $\kappa=t/2+\lambda$, si l'on pose
$$\bV\!\D_{\kappa}(\mathfrak{c},\mathfrak{n};R)=\bV^+(\mathfrak{c},\mathfrak{n};R)[\xi_\kappa]=
\bV\!\D(\mathfrak{c},\mathfrak{n};R)[\lambda^{-1}]$$
on a
$$j_{d}:G\!D_{\kappa}(\mathfrak{c},\mathfrak{n};R)\rightarrow \bV_{\kappa}(\mathfrak{c},\mathfrak{n};R)$$

On a une notion de  $q$-d{\'e}veloppement aux pointes non-ramifi{\'e}es. Consid{\'e}rons d'abord les pointes non-ramifi{\'e}es standards (donn{\'e}es par la rigidification canonique des VAHB de Tate).
Observons que pour toute $(R,\mathfrak{n})$-pointe $\mathcal{C}$ de $M$,
et tout $\sigma\in \Sigma^\mathcal{C}$,
le pull-back du fibr{\'e} $\overline{\mathcal{L}}$
sur $\mathfrak{G}$ {\`a} la VAHB de Tate $\Tate_{\mathcal{C}}(q)$ sur
 $\overline{S}_{\mathcal{C},\sigma}$ est trivial. 
La base $1$ fournit alors un morphisme de la compl{\'e}tion $p$-adique de 
$\overline{S}_{\mathcal{C},\sigma}$ vers $ T^+_{\infty,\infty}$,
qui fournit un homomorphisme d'alg{\`e}bres
$$\ev_{\mathcal{C}}:\bV\!\D(\mathfrak{c},\mathfrak{n};R)\rightarrow 
R[[q^\xi;\xi\in X_+\cup\{0\}]].$$

On a un principe du $q$-d{\'e}veloppement (cons{\'e}quence de l'irr{\'e}ductibilit{\'e} g{\'e}om{\'e}trique
de la tour des $T^+_{m,n}$) qui affirme l'injectivit{\'e} de $\ev_{\mathcal{C}}$ et m{\^e}me son ``{universelle
injectivit{\'e}}'': si $R\subset R'$, si $f$ est  d{\'e}finie sur $R'$ et si $\ev_{\mathcal{C}}(f)$ est 
{\`a} coefficients dans $R$, alors $f$ est d{\'e}finie sur $R$.

La compatibilit{\'e} de $\ev_{\mathcal{C}}$ avec la fl{\`e}che de $q$-d{\'e}veloppement des formes de Hilbert de poids demi-entier arithm{\'e}tiques
entra{\^\i}ne l'injectivit{\'e} de $j_{\kappa}$ pour chaque poids demi entier $\kappa$ et
l'injectivit{\'e} de $j_d$ sur toute alg{\`e}bre $p$-adique $\Z_p$-plate.
En outre, ceci entra{\^\i}ne que l'image de $\bV\!\D(\mathfrak{c},\mathfrak{n};R)$ 
est ferm{\'e}e dans
$R[[q^\xi;\xi\in \mathfrak{c}_+\cup\{0\}]]$ muni de la topologie de la convergence coefficient
par coefficient. La d{\'e}monstration est la m{\^e}me que pour les formes $p$-adiques de poids entier.

On d{\'e}finit le module des formes cuspidales de poids demi-entier comme le noyau de l'application
$$\bV\!\D(\mathfrak{c},\mathfrak{n};R)\rightarrow \mathcal{C}(P_\infty^{nr},R) $$
donn{\'e}e par
$f\mapsto \phi_f$ o{\`u} pour chaque $a\in (\mathfrak{o}\otimes\Z_p)^\times/\overline{\mathfrak{o}^\times}$,
$\phi_f(a)$ d{\'e}signe le terme constant du $q$-d{\'e}veloppement
$$f(\Tate_{\mathcal{C}}(q),\iota_{\can},\lambda_{\can},\alpha_{\can},\phi_{\can}\circ[a],s_{\can})\big)$$

\noindent{\bf Application:} On peut d{\'e}finir des formes modulaires de 
Hilbert $p$-adiques de poids demi-entier en formant les
$q$-d{\'e}veloppements de s{\'e}ries 
th{\^e}ta et de s{\'e}ries d'Eisenstein associ{\'e}es {\`a} des caract{\`e}res 
$p$-adiques continus de conducteur divisant
$\mathfrak{c}_0Np^\infty$. Ces caract{\`e}res sont limites
de caract{\`e}res alg{\'e}briques qui d{\'e}finissent des formes de poids demi-entier 
arithm{\'e}tiques.
Ces s{\'e}ries d{\'e}finissent donc des {\'e}l{\'e}ments de $\bV\!\D(\mathfrak{c},\mathfrak{n};R)$,
 o{\`u} $R$ est
une alg{\`e}bre $p$-adiquement compl{\`e}te contenant les coefficients des s{\'e}ries finie et plate sur $\Z_p$.

\bibliographystyle{siam}

\begin{thebibliography}{10}

\bibitem{AMRT}
{\sc A.~Ash, D.~Mumford, M.~Rapoport, and Y.~Tai}, {\em \rm Lie groups.
  History, frontiers and applications. Vol. IV : Smooth compactification of
  locally symmetric varieties}, Math Sci Press, 1975.

\bibitem{BL}
{\sc J.-L. Brylinski and J.-P. Labesse}, {\em {Cohomologie d'intersection et
  fonctions $L$ de certaines vari{\'e}t{\'e}s de Shimura}}, Ann. Sci. {\'E}c. Norm. Sup.,
  17 (1984), pp.~361--412.

\bibitem{chai}
{\sc C.-L. Chai}, {\em {Arithmetic minimal compactification of the
  Hilbert-Blumenthal moduli space}}, Ann. of Math., 131 (1990), pp.~541--554.

\bibitem{clo}
{\sc L.~Clozel}, {\em Motifs et formes automorphes: applications du principe de
  fonctorialit{\'e}}, in Automorphic forms, Shimura varieties, and $L$-functions,
  vol.~I, Academic Press, 1990, pp.~77--159.

\bibitem{De-hodge2}
{\sc P.~Deligne}, {\em {Th{\'e}orie de Hodge, II}}, Publ. Math. IHES, 40 (1971),
  pp.~5--58.

\bibitem{DePa}
{\sc P.~Deligne and G.~Pappas}, {\em {Singularit{\'e}s des espaces de modules de
  Hilbert, en les caract{\'e}ristiques divisant le discriminant}}, Compositio
  Math., 90 (1994), pp.~59--79.

\bibitem{dimdg}
{\sc M.~Dimitrov}, {\em {Compactifications arithm{\'e}tiques des vari{\'e}t{\'e}s de
  Hilbert et formes modulaires de Hilbert pour
  $\Gamma_1(\mathfrak{c},\mathfrak{n})$}}, in this volume : {Geometric Aspects of Dwork
  Theory}; A.~Adolphson, F.~Baldassarri,
  P.~Berthelot, N.~Katz and F.~Loeser, eds., Walter de Gruyter, Berlin
  2004, pp.~527--554.

\bibitem{hmv-coh}
{\sc M.~Dimitrov}, {\em {Galois representations modulo $p$ and cohomology of
  Hilbert modular varieties}},
\newblock preprint.

\bibitem{EiZa}
{\sc M.~Eichler and D.~Zagier}, {\em \rm On the theory of Jacobi forms},
  vol.~55 of Progr. Math., Birkh{\"a}user-Verlag, 1985.

\bibitem{Fa-jami}
{\sc G.~Faltings}, {\em {Crystalline cohomology and $p$-adic Galois
  representations}}, in Algebraic analysis, geometry, and number theory,
  J.~H.~U. Press, ed., 1989, pp.~25--80.

\bibitem{FaCh}
{\sc G.~Faltings and C.-L. Chai}, {\em \rm Degeneration of Abelian Varieties},
  Springer-Verlag, 1990.

\bibitem{freitag}
{\sc E.~Freitag}, {\em \rm Hilbert Modular Forms}, Springer-Verlag, 1990.

\bibitem{goren}
{\sc E.~Goren}, {\em \rm Lectures on Hilbert modular varieties and modular
  forms}, American Mathematical Society, 2002.

\bibitem{hida-igusa}
{\sc H.~Hida}, {\em {Irreducibility of generalized Igusa towers}}.
\newblock preprint.

\bibitem{hida-padic}
\leavevmode\vrule height 2pt depth -1.6pt width 23pt, {\em {\rm $p$-Adic
  Automorphic Forms on Shimura Varieties}}, Springer-Verlag, 2004.

\bibitem{hida-annarbor}
\leavevmode\vrule height 2pt depth -1.6pt width 23pt, {\em {$p$-adic
  $L$-functions for base change lifts of ${\rm GL}\sb 2$ to ${\rm GL}\sb 3$}},
  in Automorphic forms, Shimura varieties, and $L$-functions, Vol. II (Ann
  Arbor), Persp. in Math., Acad. Press, 1990, pp.~93--142.

\bibitem{hida-iwasawa}
\leavevmode\vrule height 2pt depth -1.6pt width 23pt, {\em {On nearly ordinary
  Hecke algebras for $\mathrm{GL}_2$ over totally real fields}}, in \rm
  Algebraic Number Theory - in honor of K. Iwasawa, vol.~17 of Advanced Studies
  in Pure Mathematics, 1989, pp.~139--169.

\bibitem{hida74}
\leavevmode\vrule height 2pt depth -1.6pt width 23pt, {\em {On the critical
  values of $L$-functions of $GL_2$ and $GL_2\times GL_2$}}, Duke Math. J., 74
  (1994), pp.~431--529.

\bibitem{hida-pel}
\leavevmode\vrule height 2pt depth -1.6pt width 23pt, {\em {Control theorems of
  coherent sheaves on Shimura varieties of PEL-type}}, J. Inst. 
  Math. Jussieu, 1 (2002), pp.~1--76.


\bibitem{illusie2}
{\sc L.~Illusie}, {\em {R{\'e}duction semi-stable et d{\'e}composition de complexes de
  de Rham {\`a} coefficients.}}, Duke Math. J., 60 (1990), pp.~139--185.

\bibitem{jantzen}
{\sc J.~Jantzen}, {\em \rm Representations of algebraic groups}, Academic
  Press, 1987.

\bibitem{KaOd}
{\sc N.~Katz and T.~Oda}, {\em {On the differentiation of de Rham cohomology
  classes with respect to parameters}}, J. Math. Kyoto Univ., 8-2 (1968),
  pp.~199--213.

\bibitem{KKMS}
{\sc G.~Kempf, F.~Knudsen, D.~Mumford, and B.~Saint-Donat}, {\em \rm Toroidal
  embeddings I}, vol.~339 of Lecture notes in mathematics, Springer, 1973.

\bibitem{kunn}
{\sc K.~K{\"u}nnemann}, {\em {Projective regular models for abelian varieties,
  semistable reduction, and the height pairing}}, Duke Math. J., 95 (1998),
  pp.~161--212.

\bibitem{kramer}
{\sc J.~Kramer}, {\em {An arithmetic theory of Jacobi forms in higher
  dimensions}}, J. Reine Angew. Math., 458 (1995), pp.~157--182.

\bibitem{milne}
{\sc J.~Milne}, {\em \rm {\'E}tale Cohomology}, Princeton University Press, 1980.

\bibitem{MoTi}
{\sc A.~Mokrane and J.~Tilouine}, {\em {Cohomology of Siegel varieties with
  $p$-adic integral coefficients and applications, {\rm in Cohomology of Siegel
  Varieties}}}, Ast{\'e}risque, 280 (2002), pp.~1--95.

\bibitem{mum}
{\sc D.~Mumford}, {\em {An analytic construction of degenerating abelian
  varieties over complete rings}}, Compositio Math., 24 (1972), pp.~239--272.

\bibitem{GIT}
{\sc D.~Mumford and J.~Fogarty}, {\em \rm Geometric Invariant Theory},
  Springer-Verlag, 1982.

\bibitem{rapoport}
{\sc M.~Rapoport}, {\em {Compactification de l'espace de modules de
  Hilbert-Blumenthal}}, Compositio Math., 36 (1978), pp.~255--335.

\bibitem{ribet29}
{\sc K.~Ribet}, {\em {$p$-adic interpolation via Hilbert modular forms}}, in
  Proceedings of Symposia of Pure Mathematics, vol.~29, 1975, pp.~581--592.

\bibitem{sh-half1}
{\sc G.~Shimura}, {\em {On Eisenstein series of half-integral weight}}, Duke
  Math. J., 52 (1985), pp.~281--314.

\bibitem{sh-half2}
\leavevmode\vrule height 2pt depth -1.6pt width 23pt, {\em {On Hilbert modular
  forms of half-integral weight}}, Duke Math. J., 55 (1987), pp.~765--838.

\bibitem{geer}
{\sc G.~van~der Geer}, {\em \rm Hilbert Modular Surfaces}, Springer-Verlag,
  1988.

\bibitem{W}
{\sc A.~Weil}, {\em \rm Introduction {\`a} l'{\'e}tude des vari{\'e}t{\'e}s k{\"a}hl{\'e}riennes},
  Hermann, 1958.

\bibitem{wu}
{\sc H.-T. Wu}, {\em On $p$-adic Hilbert modular adjoint $L$-functions}, PhD
  thesis, UCLA, 2001.

\bibitem{Za}
{\sc Yu.~Zarhin}, {\em On equations defining moduli of abelian
  varieties with endomorphisms in a totally real field}, Trans. Moskow Math.
  Soc., 42 (1981), pp.~1--46.

\end{thebibliography}

\bigskip

Universit{\'e} Paris 7 \hfill Universit{\'e} Paris 13

UFR de Math{\'e}matiques,  Case 7012 \hfill LAGA, Institut Galil{\'e}e

2 place Jussieu \hfill 99, avenue J.-B. Cl{\'e}ment

75251 Paris \hfill 93430 Villetaneuse

FRANCE \hfill FRANCE

\medskip
 \texttt{dimitrov@math.jussieu.fr}  \hfill \texttt{tilouine@math.univ-paris13.fr}

\end{document}